\documentclass[a4paper,12pt]{article}

\begin{document}

\title{The two-sided Laplace transformation over the Cayley-Dickson
algebras and its applications.}
\author{Ludkovsky S.V.}
\date{25 March 2006}
\maketitle
\begin{abstract}

This article is devoted to the noncommutative version of the Laplace
transformation. New types of the direct and inverse transformations
of the Laplace type over general Cayley-Dickson algebras, in
particular, also the skew field of quaternions and the  octonion
algebra are investigated. Examples are given. Theorems are proved
about properties of such transformations and also theorems about
images and originals in conjuction with operations of
multiplication, differentiation, integration, convolution, shift and
homothety. Applications are given to a solution of differential
equations over the Cayley-Dickson algebras.
\end{abstract}

\section{Introduction.}
The classical transformation of Laplace plays very important role in
complex analysis and has numerous applications, including
differential equations \cite{vladumf,lavrsch,polbremm}. In this
paper the noncommutative analog of the classical Laplace
transformation is defined and investigated, as well as new types of
the Laplace transforms over the Cayley-Dickson algebras. For this
there are used the preceding results of the author on holomorphic
and meromorphic functions of the Cayley-Dickson numbers
\cite{ludoyst,ludfov}. As it is well-known, quaternions and
operations over them were first defined and investigated by W.R.
Hamilton in the third quarter of the 19 century \cite{hamilt}. Later
on Cayley and Dickson had introduced generalizations of quaternions
known now as the Cayley-Dickson algebras
\cite{baez,kansol,kurosh,rothe}. They, in particular quaternions and
octonions, have found applications in physics, for example, they
were used by Maxwell, Yang and Mills while derivation of their
equations, which they then have rewritten in the real form because
of the insufficient development of mathematical analysis over such
algebras in their time  \cite{emch,guetze,lawmich}. This is actual,
since non commutative gauge fields are widely used in theoretical
physics \cite{solov}. In their proceedings Hamilton, Cayley and
Dickson have pointed out on the necessity of the development of
analysis over quaternions, octonions and so on for a subsequent
their using in mechanics and quantum field theory. In the latter
time it was studied mathematical analysis over graded algebras
\cite{berez,khren,oystaey}, but it frequently can be characterized
as super commutative, for example, over Grassmann algebras, but the
case of the  Cayley-Dickson algebras was not studied by another
authors, which correspond to the super noncommutative analysis.
\par Remind that the Cayley-Dickson algebra ${\cal A}_r$ over the
field $\bf R$ has $2^r$ generators $\{ i_0,i_1,$ \\
$...,i_{2^r-1} \} $ such that $i_0=1$, $i_j^2=-1$ for each
$j=1,2,...,2^r-1$, $i_ji_k=-i_ki_j$ for every $1\le k\ne j \le
2^r-1$, where $r\ge 1$. At the same time the algebra ${\cal
A}_{r+1}$ is formed from the algebra ${\cal A}_r$ with the help of
the doubling procedure by generator $i_{2^r}$, in particular, ${\cal
A}_1=\bf C$  coincides with the field of complex numbers, ${\cal
A}_2=\bf H$ is the skew field of quaternions, ${\cal A}_3$ is the
algebra of octonions, ${\cal A}_4$ is the algebra of sedenions. The
skew field of quaternions is associative, and the algebra of
octonions is alternative. The algebra ${\cal A}_r$ is power
associative, that is, $z^{n+m}=z^nz^m$ for each $n, m \in \bf N$ and
$z\in {\cal A}_r$, it is non associative and non alternative for
each $r\ge 4$.
\par At the beginning of this article the Laplace transformation is
defined. Then new types of the direct and inverse transformations of
the  Laplace type over the general Cayley-Dickson algebras are
investigated, in particular, also over the quaternion skew field and
the algebra of octonions. Examples are given. There are proved also
necessary  theorems about an integration of functions of ${\cal
A}_r$ variables. At the same time circumstances of the non
commutative Laplace transfornation are elucidated, for example,
related with the fact that in the Cayley-Dickson algebra ${\cal
A}_r$ there are $2^r-1$ imaginary generators $\{ i_1,...,i_{2^r-1}
\} $ apart from one in the field of complex numbers such that the
imaginary space in ${\cal A}_r$ has the dimension $2^r-1$. Theorems
about properties of images are proved, as well as theorems about
images and originals in conjunction with the operations of
multiplication, differentiation, integration, convolution, shift and
homothety. There is given also the extension of the non commutative
Laplace transform for generalized functions and there is considered
an application of the non commutative integral transformation for
solutions of differential equations. In the second section the
one-sided Laplace transformations are considered, then theorems
about their properties and examples are given. In the third section
the two-sided Laplace transformation over the Cayley-Dickson
algebras is considered and its applications to solutions of
differential equations are given.  All results of this paper are
obtained for the first time.
\section{Noncommutative integral transformations.}
\par {\bf 1. Definitions.} A function $f: {\bf R}\to {\cal A}_r$ we call
function-original, where ${\cal A}_r$ is the Cayley-Dickson algebra,
which may be, in particular, ${\bf K}=\bf H$ over the quaternion
skew field or ${\bf K}=\bf O$ the octonion algebra, if it satisfies
the following conditions $(1-3)$:
\par $(1)$ $f(t)$ satisfies the H\"older condition: $|f(t+h)-f(t)|
\le A |h|^{\alpha }$ for each $|h|<\delta $ (where $0<\alpha \le 1$,
$A=const >0$, $\delta >0$ are constants for a given $t$) everywhere
on $\bf R$ may be besides points of discontinuity of the first type.
On each finite interval in $\bf R$ the function $f$ may have only
the finite number of points of discontinuity of the first kind.
Remind, that a point $t_0$ is called the point of discontinuity of
the first type, if there exist finite left and right limits
$\lim_{t\to t_0, t<t_0} f(t) =: f(t_0-0)\in {\cal A}_r$ and
$\lim_{t\to t_0, t>t_0} f(t) =: f(t_0+0)\in {\cal A}_r$.
\par $(2)$ $f(t)=0$ for each $t<0$.
\par $(3)$ $f(t)$ increases not faster, than the exponential function,
that is there exist constants $C=const>0$, $s_0=s_0(f)\in \bf R$
such that $|f(t)|<C \exp (s_0t)$ for each $t\in \bf R$.
\par Mention, that for a bounded original it is possible to take
$s_0=0$.
\par If there exists an original \par $(4)$ $F(p;\zeta ):=\int_0^{\infty }
f(t)e^{-pt-\zeta }dt$, \\
then $F(p)$ is called the Laplace transformation at a point $p\in
\bf K$ of the function-original $f(t)$, where $\zeta \in {\cal A}_r$
is a parameter.
\par Let $\gamma : (-\infty ,\infty )\to {\cal A}_r$ be a path such that
the restriction $\gamma_l :=\gamma |_{[-l,l]}$ is rectifiable for
each $0<l\in \bf R$ and put by the definition \par $(5)$
$\int_{\gamma }f(z)dz=\lim_{l\to \infty }\int_{\gamma _l}f(z)dz$, \\
where ${\cal A}_r$ denotes the Cayley-Dickson algebra, $2\le r$,
${\cal A}_2=\bf H$ is the quaternion skew field, ${\cal A}_3=\bf O$
is the octonion algebra, and integrals by rectifiable paths were
defined in \S 2.5 \cite{ludfov}. So they are defined along curves
also, which may be classes of equivalence of paths relative to
increasing piecewise smooth mappings $\tau : [a,b]\to [a_1,b_1]$
realizing reparametrization of paths. Then we shall talk, that an
improper integral (5) converges. \par Consider now a function
$f(z,\zeta )$ defined for all $z$ from a domain $U$ and for each
$\zeta $ in a neighborhood $V$ of a curve $\gamma $ in ${\cal A}_r$.
The integral $G(z) := \int_{\gamma }f(z,\zeta )d\zeta $ converges
uniformly in a domain $U$, if for each $\epsilon
>0$ there exists $l_0>0$ such that
\par $(6)$ $|\int_{\gamma }f(z,\zeta )d\zeta - \int_{\gamma _l}
f(z,\zeta )d\zeta |<\epsilon $ for each $z\in U$ and $l>l_0$.
Analogously is considered the case of unbounded $\gamma $ in one
side with $[0,\infty )$ instead of $(-\infty ,\infty )$.

\par {\bf 2. Note.} Below there are considered theorems about an
existence of the Laplace transformation  $F(p)$ and its inverse. In
Definition 1 noncommutative analog of the usual Laplace transform
can be considered, but it is possible more general noncommutative
integral transformation, since $\bf H$ has 3 and $\bf O$ has 7
purely imaginary generators, in general ${\cal A}_r$ has $2^r-1$
purely imaginary generators, while in the complex field there is
only one purely imaginary unit. For the definition of $F(p)$ it was
used $\exp (- u)$, where \par $(1)$ $u:=pt+\zeta $ depends linearly
on $t$ and $p$. But it is possible to consider also the non linear
function $u=u(p,t)$ taking into account non commutativity of the
algebra ${\cal A}_r$.
\par {\bf 3. Definitions.} Let $\{ i_0, i_1,..., i_{2^r-1} \} $
be the standard generators of the algebra ${\cal A}_r$ with $r\ge
2$. Put \par $(1)$ $u(p,t) := p_0t + M(p,t)+\zeta _0$, where $p =
p_0i_0 + p_1i_1 + ... + p_{2^r-1}i_{2^r-1}\in {\cal A}_r$, $p_0,
p_1, ..., p_{2^r-1}\in \bf R$, \\ $(2)$ $M(p,t)=M(p,t;\zeta ) =
(p_1t+\zeta _1)[i_1 \cos (p_2t+\zeta _2) + i_2 \sin (p_2t+\zeta _2)
\cos (p_3t+\zeta _3) +i_3 \sin (p_2t+\zeta
_2) \sin (p_3t+\zeta _3)]$ for quaternions; \\
$(3)$ $M(p,t)=M(p,t;\zeta ) = (p_1t+\zeta _1)[ i_1 \cos
(p_2t +\zeta _2) + i_2 \sin (p_2t+\zeta _2)$ \\
$\cos (p_3t+\zeta _3) +...+ i_6 \sin (p_2t+\zeta _2) ...\sin
(p_6t+\zeta _6) \cos (p_7t+\zeta _7)$ \\ $+ i_7\sin (p_2t+\zeta
_2)...\sin (p_6t+\zeta _6) \sin (p_7t+\zeta _7)]$ \\ for octonions,
\par $(3')$ $M(p,t)=M(p,t;\zeta ) = (p_1t+\zeta _1)[ i_1 \cos
(p_2t +\zeta _2) + i_2 \sin (p_2t+\zeta _2)$ \\
$\cos (p_3t+\zeta _3) +...+ i_{2^r-2} \sin (p_2t+\zeta _2) ...\sin
(p_{2^r-2}t+\zeta _{2^r-2}) \cos (p_{2^r-1}t+\zeta _{2^r-1})$ \\ $+
i_{2^r-1}\sin (p_2t+\zeta _2)...\sin (p_{2^r-2}t+\zeta _{2^r-2})
\sin (p_{2^r-1}t+ \zeta _{2^r-1})]$ \\ for the general
Cayley-Dickson algebra with $2\le r<\infty $, where  $\zeta -\zeta
_0 = \zeta _1i_1+...+\zeta _{2^r-1}i_{2^r-1}\in {\cal A}_r$ is the
parameter of an initial phase, $\zeta _j\in \bf R$ for each
$j=0,1,...,2^r-1$. More generally, let \par $(4)$
$u(p,t)=u(p,t;\zeta )=(p_0t + v(p,t))+\zeta _0$, where $v(p,t)$ is a
locally analytic function, $Re (v(p,t))=0$ for each $p\in {\cal
A}_r$ and $t\in \bf R$, $Re (z) := (z+{\tilde z})/2$, ${\tilde
z}=z^*$ is the conjugated number for $z\in {\cal A}_r$. Then the
more general Laplace transformation over ${\cal A}_r$ is deined by
the formula:
\par $(5)$ $F_u(p;\zeta ):=\int_0^{\infty }f(t)\exp (-u(p,t))dt$ \\
for each numbers $p\in {\cal A}_r$, for which the integral exists.
At the same time the components $p_j$ of the number $p$ and $\zeta
_j$ for $\zeta $ in $v(p,t)$ we write in the $p$- and $\zeta
$-representations respectively such that
\par $(6)$ $h_j=(-hi_j+ i_j(2^r-2)^{-1} \{ -h
+\sum_{k=1}^{2^r-1}i_k(hi_k^*) \} )/2$ for each $j=1,2,...,2^r-1$, \\
$h_0=(h+ (2^r-2)^{-1} \{ -h +\sum_{k=1}^{2^r-1}i_k(hi_k^*) \} )/2$, \\
where $2\le r\in \bf N$, $h=h_0i_0+...+h_{2^r-1}i_{2^r-1}\in {\cal
A}_r$, $h_j\in \bf R$ for each $j$, $h^*={\tilde h}$ denotes the
conjugated number for $h$, $h\in {\cal A}_r$. Denote $F_u(p;\zeta )$
by ${\cal F}(f,u;p;\zeta )$.
\par Henceforth, there are used $u(p,t)$ given by $2(1)$
or $3(1-3,3')$, if another form $3(4)$ is not elucidated.
\par {\bf 4. Note.} The using of the function $M(p,t)$, given by
Formulas $(2,3)$ from Definition 3 appears naturally while
consideration of the spherical coordinates in the base of generators
of the algebra ${\cal A}_r$, as well as while consideration of
iterated exponential functions. Let $ \{ i_0, i_1, i_2, i_3 \} $ be
the standard generators of the quaternion algebra $\bf H$, where
$i_0=1$, $i_1^2=i_2^2=i_3^2=-1$, $i_1i_2= -i_2i_1=i_3$,
$i_2i_3=-i_3i_2=i_1$, $i_3i_1=-i_1i_3=i_2$, then
\par $\exp (i_1(p_1t+\zeta _1)\exp (-i_3(p_2t
+\zeta _2) \exp (-i_1(p_3t+\zeta _3)))= \exp (i_1(p_1t+\zeta _1)\exp
(- (p_2t+\zeta _2)(i_3\cos (p_3t+\zeta _3) - i_2\sin (p_3t+\zeta
_3))))$
\par $= \exp (i_1(p_1t+\zeta _1)(\cos (p_2t+\zeta _2) -
\sin (p_2t+\zeta _2)(i_3\cos (p_3t+\zeta _3) - i_2\sin (p_3t+\zeta
_3))))$
\par $= \exp ((p_1t+\zeta _1)(i_1\cos (p_2t+\zeta _2) + i_2
\sin (p_2t+\zeta _2)\cos (p_3t+\zeta _3) + i_3\sin (p_2t+\zeta
_2)\sin (p_3t+\zeta _3)))$.
\par Further by induction the equality is accomplished:
\par $\exp (\mbox{ }_{r+1}M(p,t;\zeta ))=
\exp \{\mbox{ }_rM((i_1p_1+_...+i_{2^r-1}p_{2^r-1}),t; (i_1\zeta
_1+...+i_{2^r-1} \zeta _{2^r-1})\exp (-i_{2^r+1}($ \\ $p_{2^r}t
+\zeta _{2^r}) \exp (-\mbox{
}_rM((i_1p_{2^r+1}+...+i_{2^r-1}p_{2^{r+1}-1}),t; (i_1\zeta
_{2^r+1}+...+i_{2^r-1} \zeta _{2^{r+1}-1}))) \} $, \\
where $i_{2^r}$ is the generator of the doubling of the algebra
${\cal A}_{r+1}$ from the algebra ${\cal A}_r$, such that
$i_ji_{2r}=i_{2^r+j}$ for each $j=0,...,2^r-1$, the function
$M(p,t;\zeta )$ from Definition 3 over ${\cal A}_r$ is written with
the lower index $\mbox{ }_rM$.
\par Write an image in the form $F_u(p;\zeta ):=
\sum_{s=0}^{2^r-1} F_{u,s}(p;\zeta )i_s$, where a function $f$ is
decomposed in the form $f(t)=\sum_{s=0}^{2^r-1} f_s(t)i_s$, $f_s:
{\bf R}\to \bf R$ for each $s=0,1,...,2^r-1$, $F_{u,s}(p;\zeta )$
denotes the image of the function-original $f_s$.
\par It can be taken the automorphism of the algebra ${\cal A}_r$
and instead of the standard generators $ \{ i_0,...,i_{2^r-1} \} $
to use new generators $ \{ N_0,...,N_{2^r-1} \} $, and to provide
also $M(p,t;\zeta )=M_N(p,t;\zeta )$ relative to a new basic
generators, where $2\le r\in \bf N$. In this more general case an
image we denote by $\mbox{ }_NF_u(p;\zeta )$ for the original $f(t)$
or in more details we denote it by $\mbox{ }_N {\cal F}(f,u;p;\zeta
)$.
\par If $p=(p_0,p_1,0,...,0)$, then the Laplace transformation
1(4) and 2(3,5) reduces to the complex case, that is, given above
definitions for quaternions, octonions and general Cayley-Dickson
algebras are justified.
\par {\bf 5. Theorem.} {\it Let $V$ be a bounded neighborhood
of a rectifiable curve $\gamma $ in ${\cal A}_r$, and a sequence of
functions $f_n: V\to {\cal A}_r$ be uniformly convergent on $V$,
where $2\le r<\infty $, then there exists the limit
\par $(i)$ $\lim_{n\to \infty }\int_{\gamma } f_n(z)dz=\int_{\gamma
}\lim_{n\to \infty }f_n(z)dz$.}
\par {\bf Proof.} For a given $\epsilon >0$ in view of
the uniform convergence of the sequence $f_n$ there exists $n_0\in
\bf N$ such that $|f_n(z)-f(z)|<\epsilon /l$ for each $n>n_0$, where
$0<l<\infty $ is the length of the rectifiable curve $\gamma $,
$f(z) := \lim_{n\to \infty }f_n(z)$. In view of the Inequality
2.7(4) \cite{ludfov} there are only two positive constants $C_1>0$
and $C_2>0$ such that $|\int_{\gamma }f(z)dz - \int_{\gamma
}f_n(z)dz|< (\epsilon /l) lC_1 \exp (C_2R^s)=\epsilon C_1\exp
(C_2R^s)$, where $s=2^r+2$, $0<R<\infty $, such that $V$ is
contained in the ball $B({\cal A}_r,z_0,R)$ in ${\cal A}_r$ of the
radius $R$ with the center at some point $z_0\in \bf K$. This means
the validity of Equality $(i)$.
\par {\bf 6. Theorem.} {\it If a function $f(z,\zeta )$ holomorphic by
$z$ is piecewise continuous by $\zeta $ for each $z$ from a simply
connected (open) domain $U$ in ${\cal A}_r$ with $2\le r<\infty $
and for each $\zeta $ from a neighborhood $V$ of the path $\gamma $,
where $\gamma _l$ is rectifiable for each $0<l<\infty $, and the
integral $G(z):=\int_{\gamma }f(z,\zeta )d\zeta $ converges
uniformly in the domain $U$, then it is  the holomorphic function in
$U$.}
\par {\bf Proof.} For each $0<l<\infty $ the function
$\int_{\gamma _l} f(z,\zeta )d\zeta =: G_l(z)$ is continuous by $z$
in view of Theorem 5, that together with 1(6) in view of the
triangle inequality gives the continuous function $G(z)$ on $U$. In
view of Theorems 2.16 and 3.10 \cite{ludfov} the integral
holomorphicity of the function $G(z)$ implies its holomorphicity.
But the integral holomorphicity is sufficient to establish in the
interior $Int (B({\cal A}_r,z_0,R))$ of each ball $B({\cal
A}_r,z_0,R)$ contained in $U$.  Let $\psi $ be a rectifiable path
such that $\psi \subset Int (B({\cal A}_r,z_0,R))$. Therefore,
$\int_{\psi } G(z)dz = \int_{\psi } (\int_{\gamma }f(z,\zeta )d\zeta
)dz$.  With the help of the  proof of Theorem 2.7 \cite{ludfov}
these integrals can be rewritten in the real можно переписать в
действительных coordinates and with the generators
$i_0,...,i_{2^r-1}$ of the Cayley-Dickson algebra ${\cal A}_r$. In
view of the uniform convergence $G(z)$ and the Fubini Theorem it is
possible to change the order of the integration and then $\int_{\psi
} G(z)dz = \int_{\gamma } (\int_{\psi }f(z,\zeta )dz )d\zeta =0$,
since $\int_{\psi }f(z,\zeta )dz=0$.
\par {\bf 7. Theorem.} {\it For each original $f(t)$ its
image $F_u(p;\zeta )$ is defined in the half space $\{ p\in {\cal
A}_r: Re (p)>s_0 \} $, where $2\le r\in \bf N$, and the indicator of
the growth of the function $f(t)$ is not greater, than $s_0$,
moreover, $F_u(p;\zeta )$ is holomorphic by $p$ in this half space,
and it is also holomorphic by $\zeta \in {\cal A}_r$.}
\par {\bf Proof.} Integral 1(4) or 3(5) is absolutely convergent
for $ Re (p)>s_0$, since it is majorized by the converging integral
\par $|\int_0^{\infty }f(t)\exp (-u(p,t))dt| \le \int_0^{\infty }
C\exp (-(s-s_0)t-\zeta _0)dt = C(s-s_0)^{-1}e^{-\zeta _0}$, \\
since $|e^z|=\exp (Re (z))$ for each $z\in {\cal A}_r$ in view of
Corollary 3.3 \cite{ludfov}. While an integral, produced from the
integral 1(4) or 3(5) differentiating by $p$ converges also
uniformly:
\par $(i)$ $|\int_0^{\infty }f(t)[\partial \exp (-u(p,t))/
\partial p].hdt|
\le |h| \int_0^{\infty } C t \exp (-(s-s_0)t-\zeta _0)dt =
|h|C(s-s_0)^{-2}e^{-\zeta _0}$ \\
for each $h\in {\cal A}_r$, since each $z\in {\cal A}_r$ can be
written in the form $z=|z|\exp (M)$, where $|z|^2=z{\tilde z}\in
[0,\infty )\subset \bf R$, $M\in {\cal A}_r$, $Re (M):= (M+{\tilde
M})/2=0$ in accordance with Proposition 3.2 \cite{ludfov}. In view
of Equations 3(6)
\par $\partial (\int_0^{\infty }f(t)\exp (- u(p,t;\zeta
))dt)/\partial {\tilde p}=0$ and \par $\partial (\int_0^{\infty
}f(t)\exp (- u(p,t;\zeta ))dt)/\partial {\tilde \zeta }=0$, while
\par $|\int_0^{\infty
}f(t)[\partial \exp (- u(p,t;\zeta ))/\partial \zeta ].hdt| \le |h|
\int_0^{\infty } C \exp (-(s-s_0)t-\zeta _0)dt = |h|
C(s-s_0)^{-1}e^{-\zeta _0}$ \\
for each $h\in {\cal A}_r$. In view of convergence of integrals
given above $F_u(p;\zeta )$ is (super)differentiable by $p$ and
$\zeta $, moreover, $\partial F_u(p;\zeta )/\partial {\tilde p}=0$
and $\partial F_u(p;\zeta )/{\tilde \zeta }=0$ in the considered
$(p,\zeta )$-representation, consequently, $F_u(p;\zeta )$ is
holomorphic by $p\in \bf K$ with $Re (p)>s_0$ and $\zeta \in {\cal
A}_r$ in view of Theorem 6.
\par {\bf 8. Note.} In view of Formula $7(i)$, if $p$ tends to
the infinity such that $Re (p)=:s$ converges to $+ \infty $, then
$\lim_{s\to +\infty } F_u(p;\zeta )=0$. Write $p$ in the polar form
$p=\rho e^{M\theta }$, $M\in {\cal A}_r$, $|M|=1$, $|\theta |<\pi
/2-\delta $, where $\rho \ge 0$, where a constant $\delta $ is such
that $0<\delta <\pi /2$, then $F_u(p;\zeta )$ tends to zero
uniformly by such $\theta $ and $M$ while $p$ tends to the infinity,
that is, for $\rho $ tending to the infinity, since $e^{M\theta
}=\cos (\theta ) + \sin (\theta )M$.
\par {\bf 9. Definitions and Notations.} 1. Let $\bf 1$ denotes
the unit operator on ${\cal A}_r$, that is, ${\bf 1}(h)=h$ for each
$h\in {\cal A}_r$, where ${\cal A}_r$ denotes the Cayley-Dickson
algebra with $2^r$ standard generators $ \{ i_0,i_1,...,i_{2^r-1} \}
$, $i_0=1$, $i_p^2=-1$ for each $p\ge 1$, $i_pi_s=-i_si_p$ for each
$1\le s\ne p\ge 1$. Consider also the operator of conjugation
${\tilde {\bf 1}}(h)={\tilde h}$ for each $h\in {\cal A}_r$.
\par 2. We shall distinguish symbols: $e = {\bf 1}(1)$ и ${\tilde e} =
{\tilde {\bf 1}}(1)$, where $1\in {\cal A}_r$.
\par 3. $\bf 1$ and ${\tilde {\bf 1}}$ we shall not identify
or mutually suppress with each other or with symbols of the type
${\bf 1} {\tilde {\bf 1}}$, or with others basic symbols; while
symbols $e$ and ${\tilde e}$ we shall not identify or suppress with
each other or with expressions of the form $e{\tilde e}$ or with
others basic symbols.
\par 4. $(z^pz^q)$ is identified with $z^{p+q}$, while
$({\tilde z}^p{\tilde z}^q)$ is identified with ${\tilde z}^{p+q}$,
where $p$ and $q$ are natural numbers, $p, q=1,2,3,...$. In each
word we suppose that in each fragment of the form $...(z^p{\tilde
z}^q)...$ the symbol $z^p$ is displayed on the left of the symbol
${\tilde z}^q$.
\par 5. If $w$ is some word or
a phrase, $c\in \bf R$, then the phrases $cw$ and $wc$ are
identified; if $b, c\in \bf R$, then the phrases $b(cw)$ and $(bc)w$
are identified.
\par 6. While the conjugation we shall use the identifications
of words: $(a({\bf 1}b))^*=({\tilde b}{\tilde {\bf 1}}){\tilde a}$,
$(a(eb))^*=({\tilde b}{\tilde e}){\tilde a}$ ; while the inversion
gives the identifications: $(a({\bf 1}b))^{-1}=(b^{-1}{\bf
1})a^{-1}$, $(a({\tilde {\bf 1}}b))^{-1}=(b^{-1}{\tilde {\bf
1}})a^{-1}$, $(a(eb))^{-1}=(b^{-1}e)a^{-1}$, $(a({\tilde
e}b))^{-1}=(b^{-1}{\tilde e})a^{-1}$, where $ab\ne 0$, $a, b\in
{\cal A}_r$.
\par 7. A word $\{ w_1...w_{k-1} (\{ a_1...a_l \} )_{q(l)} w_{k+1}...w_j
\} _{q(j)}$ is identified with a word \\ $\{
w_1...w_{k-1}w_kw_{k+1}...w_j \} _{q(j)}$, as well as with words $ c
\{ w_1...w_{k-1}bw_{k+1}...w_j \} _{q(j)}$, \\ $ \{ (c
w_1)...w_{k-1}bw_{k+1}...w_j \} _{q(j)}$, where $w_k= (\{ a_1...a_l
\} )_{q(l)}$, $a_1,...,a_l\in {\cal A}_r$, $l, k, j\in \bf N$,
$cb=w_k$, $c\in \bf R$, $b\in {\cal A}_r$.
\par 8. Phrases are sums of words. Sums may be finite or infinite
countable. A phrase $aw+bw$ is identified with a word $(a+b)w$, a
phrase $wa+wb$ is identified with a word $w(a+b)$ for each words $w$
and constants $a, b \in {\cal A}_r$, since $c=a+b\in {\cal A}_r$.
\par 9. Symbols of functions, that is, corresponding to them phrases
finite or infinite, are defined with the help of the following
initial symbols: constants from ${\cal A}_r$, $e$, ${\tilde e}$,
$z^p$, ${\tilde z}^p$, where $p\in \bf N$. Symbols of operator
valued functions are composed from the symbols of the set: $\{ {\bf
1}, {\tilde {\bf 1}}; e, {\tilde e}; z^p, {\tilde z}^p: p\in {\bf N}
\} $, where the symbol ${\bf 1}$ or ${\tilde {\bf 1}}$ is present in
a phrase. At the same time a phrase corresponding to a function or
an operator valued function can not contain words consisting of
constants only, that is, each word $w_1w_2...w_k$ of a phrase must
contain not less than one symbol $w_1, w_2,...,w_k$ from the set $\{
e, {\tilde e}; z^p, {\tilde z}^p: p\in {\bf N} \} $. For a phrase
$\nu $ and a void word $\emptyset $ put $\nu \emptyset =\emptyset
\nu =\nu $.
\par 10. For symbols from the set $\{ {\bf 1}, {\tilde
{\bf 1}}; e, {\tilde e}; z^p, {\tilde z}^p: p\in {\bf N} \} $ define
their lengths: $l(0)=0$, $l(a)=1$ for each $a\in {\cal A}_r\setminus
\{ 0 \} $, $l({\bf 1})=1$, $l({\tilde {\bf 1}})=1$, $l(e)=1$,
$l({\tilde e})=1$, $l(z^p)=p+1$, $l({\tilde z}^p)=p+1$ for each
natural number $p$. A length of a word is the sum of lengths of
composing it symbols. Analogously there are considered words and
phrases for several variables $\mbox{ }_1z,...,\mbox{ }_nz$, also
the function of a length of a word, where symbols corresponding to
different indices $v=1,...,n$ are different. \par Henceforth, we
consider phrases subordinated to the following conditions. Either
all words $w$ of a given phrase $\nu $ contain each of the symbols
$\mbox{ }_ve$, $\mbox{ }_v{\tilde e}$, $\mbox{ }_v\bf 1$ or $\mbox{
}_v{\tilde {\bf 1}}$ with the same finite constant multiplicity,
which may be dependant on the index $v=1,...,n$ or on a symbol
itself; or each word $w$ of $\nu $ contains neither $\mbox{ }_ve$
nor $\mbox{ }_v{\tilde e}$ besides words $w$ with one $\mbox{ }_ve$
and without any $\mbox{ }_vz^l$ with $l\ne 0$ or $w$ with one
$\mbox{ }_v{\tilde e}$ without any $\mbox{ }_v{\tilde z}^l$ with
$l\ne 0$. If a phrase $\nu $ in each its word $w$ contains $\mbox{
}_v\bf 1$ with multiplicity $p$ or $\mbox{ }_v{\tilde {\bf 1}}$ with
multiplicity $k$, then it is supposed to be a result of partial
differentiation $\partial ^{p+k}f(z,{\tilde z})/\partial \mbox{
}_vz^p \partial \mbox{ }_v{\tilde z}^k$; while the case of $\mbox{
}_ve$ with multiplicity $p$ or $\mbox{ }_v{\tilde e}$ with
multiplicity $k$ is supposed to be arising from $(\partial
^{p+k}f(z,{\tilde z})/\partial \mbox{ }_vz^p \partial \mbox{
}_v{\tilde z}^k).(1^{\otimes (p+k)})$. A presence of at least one of
the symbols $\mbox{ }_v\bf 1$ or $\mbox{ }_v{\tilde {\bf 1}}$ in a
phrase characterize an operator valued function apart from a
function.
\par 11. Supply the space $\cal P$ of all phrases subordinated
to the restrictions given above over ${\cal A}_r$ with the metric by
the formula:
$$ d(\nu ,\mu )=\sum_{j=0,1,2,...} l(\nu _j,\mu _j)b^j,$$
where $b$ is a fixed number, $0<b<1$, $\nu , \mu \in \cal P$,
\par $\nu =\sum_{j=0,1,2,...} \nu _j$, a function or an operator
valued function $f_j(z,{\tilde z})$ corresponding to $\nu _j$ is
homogeneous of a degree $j$, that is, $f_j(tz,t{\tilde
z})=t^jf_j(z,{\tilde z})$ for each $t\in \bf R$ and $z\in {\cal
A}_r$, $\nu _j=\sum_k \nu _{j,k}$, where $\nu _{j,k}$ are words with
corresponding homogeneous functions or operator valued functions
$f_{j,k}$ of degree $j$, $l(\nu _j,\mu _j):=\max_{k,s} l(\nu _{j,k};
\mu _{j,s})$, where $l(\nu _{j,k}; \mu _{j,s})=0$ for $\nu
_{j,k}=\mu _{j,s}$, $l(\nu _{j,k}; \mu _{j,s})=\max [l(\nu _{j,k});
l(\mu _{j,s})]$ for $\nu _{j,k}$ not equal to $\mu _{j,s}$.
\par 12. To each phrase $\nu $ there corresponds a continuous
function $f$, if a phrase is understood as a sequence of finite sums
of composing it words, such that the corresponding sequence of
functions converges. Consider the space $C^0_{\cal P}(U,X)$,
consisting of pairs $(f, \nu )$, $f\in C^0(U,X)$, where $X={\cal
A}_r^n$ or $X=X_t$ is the space of $\bf R$-polyhomogeneous
polyadditive operators on ${\cal A}_r^m$, $m\in \bf N$, $\nu \in
\cal P$ such that to a phrase $\nu $ there corresponds a continuous
function or operator valued function $f$ on a domain $U$ in ${\cal
A}_r^n$ with values in $X$, where $t=(t_1,t_2)$, $0\le t_1, t_2 \in
\bf Z$, $X_t$ is the ${\cal A}_r$-vector space generated by
operators of the form $A= \mbox{ }_1A\otimes ... \otimes \mbox{
}_{t_1+t_2}A$, $\mbox{ }_jA\in span_{{\cal A}_r}\{ \mbox{ }_1{\bf
1},...,\mbox{ }_n{\bf 1} \} $ for $j=1,...,t_1$; $\mbox{ }_jA\in
span_{{\cal A}_r} \{ \mbox{ }_1{\tilde {\bf 1}},...,\mbox{
}_n{\tilde {\bf 1}} \} $ for $j=t_1+1,...,t_1+t_2$. Introduce on
$C^0_{\cal P}(U,X)$ the metric:
$${\cal D}((f,\nu ); (g,\mu )):=
\sup_{z\in U} \| f(z,{\tilde z})-g(z,{\tilde z}) \| +d(\nu ,\mu ),$$
where to a phrase $\nu $ there corresponds a function $f$, while to
a phrase $\mu $ there corresponds a function $g$, where the norm of
a $m$-times $\bf R$-polyhomogeneous polyadditive operator is given
by the formula $\| S \| := \sup_{(\| h_1 \| =1,...,\| h_m \| =1)} \|
S(h_1,...,h_m) \| /[\| h_1 \| ... \| h_m \| ]$, where for a
function, that is, for $m=0$, there is taken the absolute value of
the function $|f|=\| f \| $ instead of the operator norm.
\par {\bf 10. Note.} For $r=\infty $ or $r=\Lambda $
with $card (\Lambda )\ge \aleph _0$ the variables $z$ and $\tilde z$
are algebraically independent over ${\cal A}_r$ due to Theorem 3.6.2
\cite{ludfov}. Since $A_m$ with $m<r$ has an algebraic embedding
into ${\cal A}_r$, then it is possible to consider functions
$f(z,{\tilde z})$ on domains $U$ in ${\cal A}_m$, which are
restrictions $f=g|_U$ of functions $g_f(z,{\tilde z})$ on domains
$W$ in ${\cal A}_r$, for which $W\cap {\cal A}_m=U$. \par In view of
Definition 9 symbols of functions, that is, corresponding to them
phrases finite or infinite are defined with the help of the initial
symbols. From these phrases are defined more general analytic
functions, in particular, $\exp (z)=e^z$, $Ln (z)$, $z^a$, $a^z$,
etc. For this local decompositions into power series are used for
locally analytic functions, as well as the Stone-Weierstrass
theorem, stating that the polynomials are dense in the space of
continuous functions on a compact canonical closed subset in $\bf
R^n$, hence also in ${\cal A}_r$ with $0<r<\infty $ (see
\cite{ludoyst,ludfov}).
\par For nonassociative algebras with $m\ge 3$ an order of
multiplication is essential, which is prescribed by an order of
brackets or by a vector $q(s)$ indicating on an order of
multiplications in a word. For example, words $(az^p)(z^qb)$ and
$a(z^{p+q}b)$ or $(az^{p+q})b$ for $a, b \in {\cal A}_m\setminus \bf
R$, $m\ge 3$, are different, where $p\ne 0$ and $q\ne 0$ are natural
numbers.
\par {\bf 11. Lemma.} {\it If two functions or operator valued functions
$f$ and $g$ are bounded on $U$ (see Definition 9.12), then ${\cal
D}((f,\nu ); (g,\mu ))<\infty $. In particular, if $U$ is compact,
then ${\cal D}((f,\nu );(g,\mu ))<\infty $ for each $(f,\nu ),
(g,\mu )\in C^0_{\cal P}(U,X)$.}
\par {\bf Proof.} For each degree of homogeneouity $j$ each word
$\nu _{j,k}$ consists no more, than of $b(j+\omega )$ symbols, where
$b=2$ for $\bf H$, $b=3$ for ${\cal A}_r$ with $r\ge 3$,
consequently, a number of such words is finite and $l(\nu _{j,k})\le
s(j+\omega )<\infty $, since a multiplicity $s_v$ or $u_v$, or
$p_v$, or $y_v$ of an appearance of $\mbox{ }_ve$ or $\mbox{
}_v{\tilde e}$, or $\mbox{ }_v\bf 1$, or $\mbox{ }_v{\tilde {\bf
1}}$ is finite and constant in all words of a given phrase may be
besides finite number of minor words associated with constants,
though it may depend on the index $v=1,...,n$ or on a symbol itself,
while each symbol may be surrounded on both sides by two constants
from ${\cal A}_r$, moreover, $\bf H$ is associative apart from
${\cal A}_r$ with $r\ge 3$, where $\omega =\omega (\nu
):=s_1+...+s_n+u_1+...+u_n+p_1+...+p_n+y_1+...+y_n$. Then $l(\nu
_{j,k};\mu _{j,s})\le b(j+y)$, where $y=\max (\omega (\nu ),\omega
(\mu ))$. Since the series $\sum_{j=0}^{\infty } [b(j+y)]b^j$
converges, then $d(\mu ,\nu )<\infty $. On the other hand,
$\sup_{z\in U}|f(z,{\tilde z})-g(z,{\tilde z})|<\infty $ for bounded
functions $f$ and $g$ on $U$, consequently, ${\cal D}((f,\nu );
(g,\mu ))<\infty $. At the same time it is known that a continuous
function on a compact set is bounded, that finishes the proof of the
second statement of this lemma.
\par {\bf 12. Theorem.} {\it Let $U$ be an open region in ${\cal
A}_r$, $r\ge 2$. Then there exists a continuous operator $S:
C^0_{\cal P}(U,X)\times \Gamma \to X$, where $\Gamma $ is a space of
all rectifiable paths in $U$ such that $S((f,\nu );\gamma
):=\int_{\gamma }(f,\nu )(z,{\tilde z})dz$.}
\par {\bf Proof.} On a space of rectifiable paths there exists a
natural metric $d_r(\gamma ,\eta )$ induced by a metric between
arbitrary subsets $A$ and $B$ in ${\cal A}_r$: $d_r(A,B):=\max (\psi
(A,B), \psi (B,A))$, $\psi (A,B):= \sup_{z\in A}\inf_{\xi \in
B}|z-\xi | $. To a function $f$ and its phrase representation $\nu $
there corresponds a unique function $g$ and its phrase
representation $\mu $, which is constructed by the following. A
function $g$ is characterized by the conditions.
\par Let $f$ be defined by a continuous function
$\xi : U^2\to {\cal A}_r$ such that
\par $(i)$ $\xi (\mbox{ }_1z,\mbox{ }_2z)|_{\mbox{ }_1z=z,
\mbox{ }_2z=\tilde z}=f(z,{\tilde z})$, \\
where $\mbox{ }_1z$ and $\mbox{ }_2z\in U$. Let also $g: U^2\to
{\cal A}_r$ be a continuous function, which is $\mbox{
}_1z$-superdifferentiable such that
\par $(ii)$ $(\partial g(\mbox{ }_1z,\mbox{ }_2z)/\partial \mbox{ }_1z).1=
\xi (\mbox{ }_1z,\mbox{ }_2z)$ on $U^2$. Then put
\par $(iii)$ ${\hat f}(z,{\tilde z}).h:={\hat f}_z(z,{\tilde z}).h:=
[(\partial g(\mbox{ }_1z,\mbox{ }_2z)/
\partial \mbox{ }_1z).h]|_{\mbox{ }_1z=z,\mbox{ }_2z=\tilde z}$
for each $h\in {\cal A}_r$. Shortly it can be written as $(\partial
g(z,{\tilde z})/ \partial z).1=f(z,{\tilde z})$ and ${\hat
f}_z(z,{\tilde z}).h:={\hat f}(z).h :=(\partial g(z,{\tilde z})/
\partial z).h$.
\par A phrase $\mu $ is constructed by the algorithm: at first in each
word $\nu _{j,k}$ of the phrase $\nu $ substitute each $e$ one time
on $z$, that gives a sum of words $\lambda _{j+1,\alpha (k,i)}$,
$\alpha =\alpha (k,i)\in \bf N$, $j, k, i \in \bf N$, that is, $\nu
_{j,k}\mapsto \sum_i\lambda _{j+1,\alpha (k,i)}$, where $i$
enumerates positions of $e$ in the word $\nu _{j,k}$, $i=1,...,s$,
where $s=s_1$. Therefore, to $\nu $ there corresponds the phrase
$\lambda =\sum_{j,\alpha }\lambda _{j,\alpha }$. The partial
differentiation $\partial f/\partial z$ was defined for functions
and their phrase representations. Consider the remainder $(\partial
\lambda /\partial z).1 - \nu =:\zeta $. If $\zeta =0$, then put $\mu
=\lambda $. If $\zeta \ne 0$, then to each word $\zeta _{j,k}$ of
this phrase $\zeta $ apply a left or right algorithm from \S 2.6
\cite{ludfov}. The first aforementioned step is common for both
algorithms. For this consider two $z$-locally analytic functions
$f_1$ and $q$ on $U$ such that $f_1$ and $q$ not commute as well as
corresponding to them phrases $\psi $ and $\sigma $. Let $\psi
^0:=\psi $, $\sigma ^0:=\sigma $, $\sigma ^{-n}:=\sigma ^{(n)}$,
$(\partial (\sigma ^n)/\partial z).1 =: \sigma ^{n-1}$ and $\sigma
^{-k-1}=0$ for some $k\in \bf N$ using the same notation $\sigma
_{j,k}^n$ for each word of the phrase $\sigma $. Then
\par $(iv)$ $(\psi \sigma )^1=\psi ^1\sigma -\psi ^2\sigma ^{-1}+
\psi ^3\sigma ^{-2}+...+(-1)^k\psi ^{k+1}\sigma ^{-k}$. In
particular, if $\psi =(a_1z^na_2)$, $\sigma =(b_1z^kb_2)$, with
$n>0$, $k>0$, $a_1, a_2, b_1, b_2\in {\cal A}_r\setminus {\bf R}I$,
then $\psi ^p=[(n+1)...(n+p)]^{-1}(a_1z^{n+p}a_2)$ for each $p\in
{\bf N}$, $\sigma ^{-s}=k(k-1)...(k-s+1)(b_1z^{k-s}b_2)$ for each
$s\in \bf N$. Apply this left algorithm by induction to appearing
neighbor subwords of a given word going from the left to right. Then
apply this to each word of a given phrase. Symmetrically the right
algorithm is:
\par $(v)$ $(\psi \sigma )^1=\psi \sigma ^1-\psi ^{-1}\sigma ^2+
\psi ^{-2}\sigma ^3+...+(-1)^n\psi ^{-n}\sigma ^{n+1}$, when $\psi
^{-n-1}=0$ for some $n\in \bf N$. Then apply this right algorithm
going from the right to the left by neighborhood subwords in a given
word.  For each word both these algorithms after final number of
iterations terminate, since a length of each word is finite. These
algorithms apply to solve the equation $(\partial \xi (z,{\tilde
z})/\partial z).1= \zeta (z,{\tilde z})$ for each $z\in U$ in terms
of phrases. Use one of these two algorithms to each word of the
phrase $\zeta $ that to get a unique phrase $\xi $ and then put $\mu
=\lambda +\xi $.
\par If $f_1$ and $q$ have series converging in $Int (B({\cal A}_r,0,R))$,
then these formulas demonstrate that there exists a $z$-analytic
function $(f_1q)^1$ with series converging in $Int (B({\cal
A}_r,z_0,R))$, since $\lim_{n\to \infty }(nR^n)^{1/n}=R$, where
$0<R<\infty $.  Since $f$ is locally analytic, then $g$ is also
locally analytic. Therefore, for each locally $z$-analytic function
$f$ and its phrase $\nu $ there exists the operator $\hat f$ and its
phrase $\partial \mu /\partial z$. Considering a function $G$ of
real variables corresponding to $g$ we get that in view of to Lemma
2.5.1 \cite{ludfov} all solutions $(g,\mu )$ differ on constants in
${\cal A}_m$, since $\partial g/\partial w_s+(\partial g/\partial
w_p).(s^*p)=0$ for each $s=i_{2j}$, $p=i_{2j+1}$,
$j=0,1,...,2^{r-1}-1$ and $\partial g/\partial w_1$ is unique, hence
$\hat f$ is unique for $f$.
\par Denote the described above mapping ${\cal P}\ni \nu \mapsto
\mu \in {\cal P}$ by $\phi (\nu )=\mu $. Thus each chosen algorithm
between two these algorithms gives $\phi (\nu _1+\nu _2)=\phi (\nu
_1)+\phi (\nu _2)=\mu =\mu _1+\mu _2$, if $\nu =\nu _1+\nu _2$,
moreover, $\omega (0)=0$. Therefore, this procedure gives the unique
branch of $S$. Recall, that if the following limit exists
$$(vi)\quad \int_{\gamma }(f,\nu )(z,{\tilde z})dz:=
\lim_P I((f,\nu );\gamma ;P),\mbox{ where}$$
$$(vii)\quad I((f,\nu );\gamma ;P):=\sum_{k=0}^{t-1}
({\hat f},{\hat {\nu }})(z_{k+1},{\tilde z}_{k+1}).(\Delta z_k),$$
where $\Delta z_k:=z_{k+1}-z_k$, $z_k:=\gamma (c_k)$ for each
$k=0,...,t$, then we say that $(f,\nu )$ is line integrable along
$\gamma $ by the variable $z$.
\par From Equation $2.7.(4)$ \cite{ludfov} it follows, that
$$(viii)\quad |S((f,\nu ) - (y,\psi );\gamma )|\le {\cal D}((f,\nu );
(y,\psi ))V(\gamma )C_1 \exp (C_2 R^{2^m+2})$$ for each $(f,\nu ),
(y,\psi )\in C^0_{\cal P}(U,X)$ for $U\subset {\cal A}_m$ with
finite $m\in \bf N$, where $C_1$ and $C_2$ are positive constants
independent from $R$, $(f,\nu )$ and $(y,\psi )$, $0<R<\infty $ is
such that $\gamma \subset B({\cal A}_m,z_0,R)$. The space $X_t$ is
complete relative to the norm of polyhomogeneous operators as well
as the Cayley-Dickson algebra ${\cal A}_r$ relative to its norm.
Therefore, $\{ S((f^v,\nu ^v);\gamma ): v\in {\bf N} \} $ is the
Cauchy sequence in $X$ for each Cauchy sequence $\{ (f^v, \nu ^v):
v\in {\bf N} \} \subset C^0_{\cal P}(U,X)$, hence $\{ S((f^v,\nu
^v);\gamma ): v\in {\bf N} \} $ converges.
\par For each rectifiable path $\gamma $ in $U$ there corresponds a
canonical closed bounded subset $V$ in ${\cal A}_r$ such that
$\gamma \subset Int (V)\subset V\subset U$. Then $V_m:=V\cap {\cal
A}_m$ is compact for each $m\in \bf N$ and a continuous function
$f|_V$ is uniformly continuous. Such that $\lim_{d(\eta ,0e)\to
0}{\cal D}((y|_{V_m},\eta ),(0,0e))=0$ for each $m\in \bf N$. For
finite $r$ take $m=r$. For infinite $r$ there exists a sequence of
rectifiable paths $\gamma _m$ in $V_m$ for suitable algebraic
embeddings of ${\cal A}_m$ into ${\cal A}_r$ such that $\lim_{m\to
\infty } d_r(\gamma _m,\gamma )=0$. In the algebra ${\cal A}_r$ the
union of all algebraically embedded subalgebras ${\cal A}_m$, $m\in
\bf N$ is dense. Therefore, for each $\epsilon _k=1/k$ there exists
a continuous functions $f_m$ on $V_m$ and a phrase $\nu _m$ over
${\cal A}_m$ for $m=m(k)$ with $m(k)< m(k+1)$ for each $k\in \bf N$
such that $\lim_{m\to \infty }{\cal D}((f_m,\nu _m),(f|_{V_m},\mu
))=0$, where each $\nu _m$ has a natural extension over ${\cal A}_m$
and hence $f_m$ as the corresponding to this phrase function has the
natural extension on a neighborhood of $V_m$ in $U$. This gives
$S((f,\nu ),\gamma )=\lim_{m\to \infty } S((f_m,\nu _m),\gamma _m)$.
\par In accordance with the proof of Theorem 2.7 \cite{ludfov} the
operator $S$ is continuous for functions and analogously for
operator valued functions and their phrases relative to the metric
${\cal D}$.
\par {\bf 13. Note.}
For a subset $U$ in ${\cal A}_r$ let $\pi _{s,p,t}(U):= \{ u: z\in
U, z=\sum_{v\in \bf b}w_vv,$ $u=w_ss+w_pp \} $ for each $s\ne p\in
\bf b$, where $t:=\sum_{v\in {\bf b}\setminus \{ s, p \} } w_vv \in
{\cal A}_{r,s,p}:= \{ z\in {\cal A}_r:$ $z=\sum_{v\in \bf b} w_vv,$
$w_s=w_p=0 ,$ $w_v\in \bf R$ $\forall v\in {\bf b} \} $, where ${\bf
b} := \{ i_0,i_1,...,i_{2^r-1} \} $ is the family of standard
generators of the algebra ${\cal A}_r$. That is, geometrically $\pi
_{s,p,t}(U)$ the projection on the complex plane ${\bf C}_{s,p}$ of
the intersection $U$ with the plane ${\tilde \pi }_{s,p,t}\ni t$,
${\bf C}_{s,p} := \{ as+bp:$ $a, b \in {\bf R} \} $, since $sp^*\in
{\hat b}:={\bf b}\setminus \{ 1 \} $. Recall that in \S \S 2.5-7
\cite{ludfov} for each continuous function $f: U\to {\cal A}_r$ it
was defined the operator $\mbox{ }_j{\hat f}$ by each variable
$\mbox{ }_jz\in {\cal A}_)r$.
\par  It is supposed further, that a domain $U$ in ${\cal A}_r^n$
has the property that each projection $p_j(U)=:U_j$ is
$2^r-1$-connected (see \cite{span}); $\pi _{s,p,t}(U_j)$ is simply
connected in $\bf C$ for each $k=0,1,...,2^r-1$, $s=i_{2k}$,
$p=i_{2k+1}$, $t\in {\cal A}_{r,s,p}$ and $u\in {\bf C}_{s,p}$, for
which there exists $z=u+t\in U_j$, where $e_j =
(0,...,0,1,0,...,0)\in {\cal A}_r^n$ is the vector with $1$ on the
$j$-th place, $p_j(z) = \mbox{ }_jz$ for each $z\in {\cal A}_r$,
$z=\sum_{j=1}^n\mbox{ }_jz e_j$, $\mbox{ }_jz\in {\cal A}_r$ for
each $j=1,...,n$, $n\in {\bf N} := \{ 1,2,3,... \} $.
\par {\bf 14. Theorem.} {\it Let a domain $U$ in ${\cal
A}_r$ satisfies the condition of Note 13 and $(f,\nu )$ be integral
holomorphic in $U$. Then $\partial S((f,\nu );\gamma )/\partial
z=({\hat f},{\hat {\nu }})(z)$ for each $z_0, z\in Int (U)$ and
rectifiable paths $\gamma $ in $Int (U)$ such that $\gamma (0)=z_0$
and $\gamma (1)=z$. In particular, if $f=(\partial g/\partial z).1$
and $\nu =(\partial \mu /\partial z).1$, where $\partial f/\partial
{\tilde z}=0$ and $\partial \nu /\partial {\tilde z}=0$ on $U$, then
$S((f,\nu );\gamma )=(g,\mu )(z)-(g,\mu )(z_0)$.}
\par {\bf Proof.} In view of the condition that $(f,\nu )$ is
integral holomorphic, that is, $S((f,\nu );\eta )=0$ for each
rectifiable loop $\eta $ in $U$, it follows, that $S((f,\nu );\gamma
)$ depends only on $z_0=\gamma (0)$ and $z=\gamma (1)$ and is
independent from rectifiable paths with the same ends. Since
$(\partial g(z,{\tilde z})/\partial z).1=f(z)$ and $(\partial \mu
(z,{\tilde z})/\partial z).1=\nu $, such that ${\hat f}(z,{\tilde
z})=\partial g(z,{\tilde z})/\partial z$ and ${\hat {\nu
}}(z,{\tilde z})=\partial \mu (z,{\tilde z})/\partial z$ exist in
the sense of distributions (see also \cite{luladfcdv}), then from
Formulas 12(vi,vii) the first statement of this theorem follows.
\par If $(f,\nu )=(\partial (g,\mu )/\partial z).1$, then
there exists $(y,\psi )\in C^0_{\cal P}(U,X)$ such that $(\partial
(y,\psi )/\partial z).1=(g,\mu )$, hence $(\partial ^2(y,\psi
)/\partial z^2)(h_1,h_2)=(\partial ^2(y,\psi )/\partial
z^2)(h_2,h_1)$ for each $h_1, h_2\in {\cal A}_r$, since $(f,\nu )$
and consequently $(g,\mu )$ and $(y,\psi )$ are holomorphic,
particularly, infinite continuous superdifferentiable in $z$ on $U$
(see Theorems 2.11, 2.16, 3.10 and Corollary 2.13 \cite{ludfov}).
Mention, that the first step of the algorithm of the proof of
Theorem 12 gives $\mu =\lambda $. For $h_1=1$ and $h_2=\Delta z$
this gives $({\hat f},{\hat {\nu }}).\Delta z=(\partial (g,\mu
)/\partial z).\Delta z$ and we can substitute the latter expression
in the integral sums instead of $({\hat f},{\hat {\nu }}).\Delta z$.
The second statement therefore follows from the fact that integral
holomorphicity of $(f,\nu )$ is equivalent to $\partial (f,\nu
)/\partial {\tilde z}=0$ on $U$, where $U$ satisfies conditions of
Note 13.
\par {\bf 15. Remark.} In the particular case of ${\bf C}={\cal A}_1$ we
simply have ${\hat f}=f$, but here the case of arbitrary $r\in \bf
N$ or $r=\Lambda $ with $card (\Lambda )\ge \aleph _0$ is
considered. For example, the pair of the function and its phrase
$((ab)(c{\tilde z}d),(aeb)(c{\tilde z}d))$ belong to $ C^0_{\cal
P}(U,{\cal A}_r)$ and there exists $(g,\mu )=((azb)(c{\tilde z}d),
(azb)(c{\tilde z}d))\in C^0_{\cal P}(U,{\cal A}_r)$ such that
$(\partial (g,\mu )/\partial {\tilde z}).1=(f,\nu )$, but $(f,\nu )$
is not integral holomorphic for each $r\ge 1$, where $a, b, c, d\in
{\cal A}_r$ are constants. In the commutative case of $\bf C$
conditions of section 9 simply reduce to the power series of complex
holomorphic functions and symbols $e$ and $\tilde e$ can be omitted,
but in the noncommutative case, when $r\ge 2$, the conditions of
section 9 are generally essential for Theorems 12 and 14.
\par {\bf 16. Lemma.} {\it Let $v: {\cal A}_r\to {\cal A}_r$ be an
automorphism of the algebra ${\cal A}_r$, and let $\{
i_0,i_1,...,i_{2^r-1} \} $ and $ \{ s_0,s_1,...,s_{2^r-1} \} $ be
generators of the algebra ${\cal A}_r$, where $v(i_j)=s_j$ for each
$j=0,1,...,2^r-1$, $2\le r< \infty $. Then the integrals
$\int_a^bf(t)dt$ and $\int_{\gamma }w(z)dz$ do not depend on a
choice of basic generators, where $f: [a,b]\to {\cal A}_r$ and $w:
U\to {\cal A}_r$ are piecewise continuous functions, $-\infty \le
a<b\le \infty $, $\gamma $ is a rectifiable curve in the domain
$U\subset {\cal A}_r$.}
\par {\bf Proof.} Each element of the Cayley-Dickson algebra
$z\in {\cal A}_r$ can be decomposed in the unique way by the base of
generators: $z = z_0i_0+...+z_{2^r-1}i_{2^r-1}$, $z = x_0s_0+...+
x_{2^r-1}s_{2^r-1}$, where $z_j, x_j\in \bf R$ for each
$j=0,...,2^r-1$. But the unit is unique in the algebra ${\cal A}_r$,
since $\zeta \zeta ^{-1}=1$ for each element $\zeta \ne 0$ in ${\cal
A}_r$, where $\zeta ^{-1} = {\tilde {\zeta }}|\zeta |^{-2}$.
Therefore, $i_0=v(i_0)=s_0$, where as usually $i_0=1$. Partition the
segment $[a,b]$ or a curve on segments $[a_k,b_k]$ or $\gamma _k$,
such that $f|_{[a_k,b_k]}$ and $w|_{U_k}$ is continuous, where $a\le
a_k<b_k=a_{k+1}\le b$ and $dim_{\bf R}(\partial (U_k\cap
U_{k+1}))\le 2^r-1$, $\gamma _k=\gamma \cap U$, $U_k\cap U_{k+1}\cap
\gamma $ is the one point set. At the same time
$\int_a^bf(t)dt=\sum_k\int_{a_k}^{b_k}f(t)dt$ and $\int_{\gamma }
w(z)dz=\sum_k \int_{\gamma _k}w(z)dz$.
\par The first integral $\int_a^bf(t)dt$ is taken along the segment
$[a,b]$ of the real line, consequently, it does not depend on the
basis of generators. In the second integral the curve $\gamma $ in
the domain $U$, as well as $\gamma _k$ with $U_k$, do not depend
naturally on the choice of the basis of generators.
\par Remind the definition of the integral along the rectifiable curve.
Let $w: U\to {\cal A}_r$ be a continuous function, where $U$ is open
in ${\cal A}_r$, $w$ defines the continuous function $\xi : U^2\to
{\cal A}_r$, such that
\par $(1)$ $\xi (\mbox{ }_1z,\mbox{ }_2z)|_{\mbox{ }_1z=z,
\mbox{ }_2z=\tilde z}=w(z,{\tilde z})$ \\
or shortly $w(z)$ instead of $w(z,{\tilde z})$, where $\mbox{ }_1z$
and $\mbox{ }_2z\in U$. Let also $g: U^2\to {\cal A}_r$ be a
continuous function, which is $\mbox{ }_1z$-superdifferentiable,
such that
\par $(2)$ $(\partial g(\mbox{ }_1z,\mbox{ }_2z)/\partial \mbox{ }_1z).1=
\xi (\mbox{ }_1z,\mbox{ }_2z)$ on $U^2$. Then put
\par $(3)$ ${\hat w}(z,{\tilde z}).h:={\hat w}_z(z,{\tilde z}).h:=
[(\partial g(\mbox{ }_1z,\mbox{ }_2z)/
\partial \mbox{ }_1z).h]|_{\mbox{ }_1z=z,\mbox{ }_2z=\tilde z}$
for each $h\in {\cal A}_r$. Shortly we can write this as $(\partial
g(z,{\tilde z})/
\partial z).1=w(z,{\tilde z})$ and ${\hat w}_z(z,{\tilde
z}).h:={\hat w}(z).h :=(\partial g(z,{\tilde z})/ \partial z).h$. If
there exists the following limit
$$(4)\quad \int_{\gamma }w(z,{\tilde z})dz:=
\lim_P I(w,\gamma ;P),\mbox{ where}$$
$$(5)\quad I(w,\gamma ;P):=\sum_{k=0}^{t-1}
{\hat w}(z_{k+1},{\tilde z}_{k+1}).(\Delta z_k),$$ where $\Delta
z_k:=z_{k+1}-z_k$, $z_k:=\gamma (c_k)$ for each $k=0,...,t$, then by
the definition $w$ is integrable along the path $\gamma $ by the
variable $z$. \par Since $i_0=s_0=1$, then ${\hat w}(z)$ also is not
dependant on the choice of the basis, since the operator ${\hat
w}(z)$ defined in the sense of distributions is given by Conditions
$(2,3)$ and the algorithm of the construction of ${\hat w}$ is
independent on the choice of the basis of generators in each
subdomain $U_k$ (see the proofs of Proposition 2.6 and Theorem 2.7
in \cite{ludoyst,ludfov}).
\par {\bf 17. Lemma.} {\it Let the straight line $a+S\theta $ be given ,
$\theta \in \bf R$ in ${\cal A}_r$ for $2\le r<\infty $, where $S\in
{\cal A}_r$, $|S|=1$, $Re (S)=0$. Then for the function $M(p,t)$
given by Equation 3(2) or 3(3,3') there exists an automorphism $v$,
$z\mapsto v(z)$, of the algebra ${\cal A}_r$, such that in the basis
of generators $s_j:=v(i_j)$ the restriction of the function
$v(M(p,t;\zeta ))$ on the straight line $p\in \{ z: z=a+S\theta ,
\theta \in \bf R \} $ is equal to
\par $(1)$ $v(M(p,t;\zeta ))= ({p'}_1t + {\zeta '}_1)K$, where \par
$(2)$ $K= [s_1 \cos ({\zeta '}_2) + s_2 \sin ({\zeta '}_2) \cos
({\zeta '}_3) +i_3 \sin ({\zeta '}_2) \sin ({\zeta '}_3)]$ for
quaternions;
\par $(3)$ $K= [ s_1 \cos ({\zeta '}_2) + s_2 \sin ({\zeta '}_2)$
$\cos ({\zeta '}_3) +...+ s_6 \sin ({\zeta '}_2) ...\sin ({\zeta
'}_6) \cos ({\zeta '}_7)$\\  $+ s_7\sin ({\zeta '}_2)...\sin ({\zeta
'}_6) \sin ({\zeta '}_7)]$ for octonions,
\par $(3')$ $K= [ s_1
\cos ({\zeta '}_2) + s_2 \sin ({\zeta '}_2)$  $\cos ({\zeta '}_3)
+...+ s_{2^r-2} \sin ({\zeta '}_2) ...\sin ({\zeta '}_{2^r-2}) \cos
({\zeta '}_{2^r-1})$\\  $+ s_{2^r-1}\sin ({\zeta '}_2)... \sin
({\zeta '}_{2^r-2}) \sin ({\zeta '}_{2^r-1})]$ for the general
Cayley-Dickson algebra, where ${p'}_1=v((p-{\tilde p})/2)$ and
${\zeta '}_n =v(\zeta _n) \in \bf R$ for each $n=1,...,2^r-1$.}
\par {\bf Proof.} On the given straight line the argument $p$
can be written in the form: $p=p_0+p_SS$, where $p_0 = Re (p) = a$
and $p_S = \theta $. Then the generators $s_0=1$ and $s_1=S$ we
complete to the basis of generators $\{ s_0,...,s_{2^r-1} \} $ in
${\cal A}_r$, using the standard procedure of subsequent doubling,
that $Re (s_{2^{n-1}}{\tilde s}_{2^n})=0$ for each $n=1,...,r-1$,
since ${\bf R}\oplus S{\bf R}$ is isomorphic with $\bf C$. In view
of $v(z+\zeta )=v(z)+v(\zeta )$ and $v(z\zeta )=v(z)v(\zeta )$ for
each $z, \zeta \in {\cal A}_r$, then $v(e^z)=e^{v(z)}$ for each
$z\in {\cal A}_r$. Therefore, $v(\exp (z\exp_n (\zeta ))=\exp
(v(z)\exp_n(v(\zeta ))$ for each $z, \zeta \in {\cal A}_r$ and $1\le
n\in \bf Z$, where $\exp _1(z):=\exp (z)$, $\exp_{n+1}(z):= \exp
(\exp_n(z))$ for each $1\le n\in \bf N$ and $z\in {\cal A}_r$. Using
the standard iterated exponents from Note 4 and Formulas
\par $h_j=(-hi_j+ i_j(2^r-2)^{-1} \{ -h
+\sum_{k=1}^{2^r-1}i_k(hi_k^*) \} )/2$ for each $j=1,2,...,2^r-1$, \\
$h_0=(h+ (2^r-2)^{-1} \{ -h +\sum_{k=1}^{2^r-1}i_k(hi_k^*) \} )/2$, \\
where $h=h_0i_0+...+h_{2^r-1}i_{2^r-1}$, $h_j\in \bf R$ for each
$j$, $h^* = {\tilde h}$ denotes the conjugated element of the
Cayley-Dickson algebra $h$, $h\in {\cal A}_r$,we get the statemenet
of this lemma, moreover, ${p'}_1=p_S=v((p-{\tilde p})/2)$.
\par {\bf 18. Proposition.} {\it Let $U$ be a domain in
the Cayley-Dickson algebra ${\cal A}_r$ with $2\le r$, $f: U^2\to
{\cal A}_r$ be a piecewise continuous function,  $\gamma :
(a_1,b_1)\to U$ and $\psi : (a_2,b_2)\to U$ be paths, where $-\infty
\le a_1<b_1\le \infty $, $-\infty \le a_2<b_2\le \infty $, while
$\gamma |_{(\alpha ,\beta )}$ and $\psi |_{(\alpha ,\beta )}$ are
rectifiable and for each $-\infty <\alpha <\beta <\infty $ the
integral $\int_{\gamma }(\int_{\psi }f(z,\zeta )d\zeta )dz$
converges. Then there exists
\par $(i)$ $\int_{\psi }(\int_{\gamma }f(z,\zeta )dz)d\zeta =
\int_{\gamma }(\int_{\psi }f(z,\zeta )d\zeta )dz$.}
\par {\bf Proof.} Since $f$ is piecewise continuous, then
there exists a countable partition of the domain $U^2$ in the form
$U^2=\bigcup_kW_k$ by subdomains $W_k$, where the codimension over
$\bf R$ of the intersection $codim_{\bf R}(W_k\cap W_l)\ge 1$ for
each $k\ne l$. At the same time restrictions $f|_{W_k}$ are
continuous. It can be chosen each $W_k$ bounded such that $(\gamma
\times \psi )\cap W_k=:\gamma _k\times \psi _k$ were with the
rectifiable $\gamma _k=\gamma |_{(\alpha _k,\beta _k)}$ and $\psi
_k=\psi |_{({\alpha '}_k,{\beta '}_k)}$. Then in view of continuity
of the restriction $f|_{W_k}$ there exists the finite integrals
\par $(ii)$ $\int_{\psi _k}(\int_{\gamma _k}f(z,\zeta )dz)d\zeta =
\int_{\gamma _k}(\int_{\psi _k}f(z,\zeta )d\zeta )dz$ for each $k$.
Since
\par $\sum_k \int_{\gamma _k}(\int_{\psi _k}f(z,\zeta )d\zeta
)dz=\int_{\gamma }(\int_{\psi }f(z,\zeta )d\zeta )dz$, then the
summation of the both parts of the Equality $(ii)$ gives the
Equation $(i)$.
\par {\bf 19. Theorem.} {\it If a function $f(t)$ is an original
(see definition 1), such that $\mbox{ }_NF_u(p;\zeta
):=\sum_{s=0}^{2^r-1} \mbox{ }_NF_{u,s}(p;\zeta )N_s$ is its image,
where the function $f$ is written in the form \\
$f(t)=\sum_{s=0}^{2^r-1} f_s(t)N_s$, $f_s: {\bf R}\to \bf R$ for
each $s=0,1,...,2^r-1$ (see Note 4), $f({\bf R})\subset \bf K$ for
${\bf K}=\bf H$ or ${\bf K}=\bf O$, $f({\bf R})\subset \bf R$ over
the Cayley-Dickson algebra ${\cal A}_r$ with $4\le r\in \bf N$. Then
at each point $t$, where $f(t)$ satisfies the H\"older condition
there is accomplished the equality:
\par $(i)$ $f(t) = (2\pi N_1)^{-1} Re (S{\tilde N}_1) \sum_{s=0}^{2^r-1}
\int_{a-S\infty }^{a+S\infty }\mbox{ }_NF_{u,s}(p;\zeta )
\exp (u(p,t;\zeta ))dp)N_s$, \\
where $u(p,t;\zeta )=pt+\zeta $ or $u(p,t;\zeta )=p_0t+M_N(p,t;\zeta
)+\zeta _0$, the integral is taken along the straight line $p(\tau
)=a+S\tau \in \bf K$, $\tau \in \bf R$, $S\in {\cal A}_r$, $Re
(S)=0$, $|S|=1$, $Re (p) = a>s_0$ and the integral is understood in
the sense of the principal value.}
\par {\bf Proof.} In view of the decomposition of a function $f$
in the form $f(t)=\sum_{s=0}^{2^r-1} f_s(t)N_s$ it is sufficient to
consider the inverse transformation of the real valued function
$f_s$, which we denote for simplicity by $f$. Since $t\in \bf R$,
then $\int_0^{\infty }f(\tau )d\tau $ is the Riemann integral. If
$w$ is a holomorphic function of the Cayley-Dickson variable, then
locally in a simply connected domain $U$ in each ball $B({\cal
A}_r,z_0,R)$ with the center at $z_0$ of radius $R>0$ contained in
the interior $Int (U)$ of the domain $U$ there is accomplished the
equality
\par $(\partial \int_{z_0}^zw(\zeta )d\zeta /\partial z).1=w(z)$, \\
where the integral depends only on an initial $z_0$ and a final $z$
points of a rectifiable path in $B({\cal A}_r,z_0,R)$. On the other
hand, along the straight line $a+S\bf R$ the restriction of the
antiderivative has the form $\int_{\theta _0}^{\theta }w(a+S\tau
)d\tau $, since\par $\int_{z_0=a+S\theta _0}^{z=a+S\theta }w(\zeta
)d\zeta =\int_{\theta _0}^{\theta } {\hat
w}(a+S\tau ).Sd\tau $, \\
 $\partial f(z)/\partial \theta =(\partial f(z)/\partial z).S$ for
superdifferentiable by $z\in U$ function $f(z)$, moreover, the
antiderivative is unique up to a constant from ${\cal A}_r$ with the
given representation of the function and the algorithm
\cite{ludfov,ludoyst}.
\par The integral $g_b(t) := \int_{a-SB}^{a+SB} \mbox{ }_NF_{u,s}(p;\zeta
)\exp (p_0t+M_N(p,t;\zeta )+\zeta _0)dp$ for each $0<b<\infty $ with
the help of generators of the algebra ${\cal A}_r$ and the Fubini
Theorem for real valued components of the function can be written in
the form:
\par $g_b(t) = (2\pi N_1)^{-1} Re (S{\tilde N}_1)
\int_0^{\infty }f(\tau )d\tau \int_{a-Sb}^{a+Sb} \exp
(p_0t+M_N(p,t;\zeta )+\zeta _0)$ \par $\exp (-(p_0\tau +M_N(p,\tau
;\zeta ) +\zeta _0))dp$, \\
since the integral $\int_0^{\infty }f(\tau )\exp (-(p_0\tau
+M_N(p,\tau ;\zeta ) +\zeta _0))d\tau $ is uniformly converging
relative to $p$ in the half space $Re (p)>s_0$ in ${\cal A}_r$ (see
also Proposition 18). In view of the alternativity of the algebra
$\bf K$ and the power associativity of the algebra ${\cal A}_r$ it
can be considered the automorphism $v$ from Lemma 17 $z\mapsto v(z)$
for each $z\in \bf K$. Then with such $v$ we get the function
$M(p,t;\zeta )$ given by  Formulas $17(1-3)$, where $p_0,{p'}_1,
\zeta _0, {\zeta '}_1,...,{\zeta '}_{2^r-1}\in \bf R$ are constant,
$t\in \bf R$, $K\in {\cal A}_r$ is constant with $Re (K)=0$, where
${p'}_0=p_0$ and ${\zeta '}_0=\zeta _0$, since $v(1)=1$ and,
consequently, $v(t)=t$ for each $t\in \bf R$. Formula $(i)$ is
satisfied if and only if it is accomplished after application of the
automorphism $v$ to both parts of the Equality, since $v(z)=v(\zeta
)$ for $z, \zeta \in {\cal A}_r$ is equivalent to that $z=\zeta $.
\par Then in view of Lemmas 16 and 17
up to an automorphism of the algebra ${\cal A}_r$ the proof reduces
to the case $p=(p_0,p_1,0,...,0)$, $N=(N_0,N_1,N_2,...,N_{2^r-1})$,
where $N_0=1$, since $\bf R$ is the center of the algebra ${\cal
A}_r$. But this gives $p_1=p_1(t)=Re (S{\tilde N}_1)t$ for each
$t\in \bf R$. Consider the particular case $c := Re (S{\tilde
N}_1)\ne 0$, then the particular case $Re (S{\tilde N}_1)=0$ is
obtained by taking the limit when $Re (S{\tilde N}_1)\ne 0$ tends to
zero. Thus,
\par $g_b(t) = (2\pi N_1)^{-1} c \int_0^{\infty }f(\tau )d\tau
\int_{a-Sb}^{a+Sb} \exp (
at+c(\zeta _1+ t)K) \exp (-(a\tau +c(\zeta _1+ \tau )K)dp$, \\
since $\zeta _0, a\in \bf R$, where $K$ is given by Formulas
17(2,3). Then
\par $g_b(t) = (\pi N_1)^{-1} c
\int_0^{\infty }f(\tau )e^{a(t-\tau )}[\sin (bc(t-\tau ))] (ct-c\tau
)^{-1}$ \par  $ =(\pi )^{-1} e^{at} \int_{-t}^{\infty } f(\zeta +t)
e^{-a(\zeta +t)}[\sin (b\zeta )]\zeta ^{-1}d\zeta $, \\
where it can be used the substitution $\tau -t= \zeta $. Put
$w(t):=f(t)e^{-at}$, where $w(t)=0$ for each $t<0$. Therefore,
\par $g_b(t)=(\pi )^{-1}e^{at}\int_{-\infty }^{\infty }[w(\zeta +t)-
w(t)]\zeta ^{-1}\sin (b\zeta )d\zeta +(\pi )^{-1}f(t)\int_{-\infty
}^{\infty }{\zeta }^{-1}\sin (b\zeta )d\zeta $. The integral in the
second term is known as the Euler integral: $\int_{-\infty }^{\infty
}{\zeta }^{-1}\sin (b\zeta )d\zeta =\pi $ for each $b>0$,
consequently, the second term is equal to $f(t)$. It remains to
prove, that $\lim_{b\to \infty }\int_{-\infty }^{\infty }[w(\zeta
+t)- w(t)]\zeta ^{-1}\sin (b\zeta )d\zeta =0$, that follows from the
subsequent lemma.
\par {\bf 20. Lemma.} {\it If a function $\psi (\zeta )$ with values in
the Cayley-Dickson algebra ${\cal A}_r$ is integrable on the segment
$[\alpha ,\beta ]$, then \par $\lim_{b\to \infty } \int_{\alpha
}^{\beta }\psi (\zeta )\sin (b\zeta )d\zeta =0$.}
\par {\bf Proof.} If $\psi $ is continuously differentiable on the
segment $[\alpha ,\beta ]$, then the result of the integration by
parts is:
\par $\int_{\alpha }^{\beta }\psi (\zeta )\sin (b\zeta )d\zeta
=-\psi (\zeta )\cos (b\zeta )b^{-1}|_{\alpha }^{\beta }+\int_{\alpha
}^{\beta }\psi '(\zeta )\cos(b\zeta )b^{-1}d\zeta $ and,
consequently,
\par $\lim_{b\to \infty } \int_{\alpha }^{\beta }\psi (\zeta )\sin
(b\zeta )d\zeta =0$. If $\psi (\zeta )$ is an arbitrary integrable
function, then for each $\epsilon >0$ there exists a continuous
differentiable function $\psi _{\epsilon }(\zeta )$ such that
$\int_{\alpha }^{\beta }|\psi (\zeta )- \psi _{\epsilon }(\zeta
)|d\zeta <\epsilon /2$. Then $\int_{\alpha }^{\beta }\psi (\zeta
)\sin (b\zeta )d\zeta =\int_{\alpha }^{\beta } [\psi (\zeta )-\psi
_{\epsilon }(\zeta )]\sin (b\zeta )d\zeta +\int_{\alpha }^{\beta
}\psi _{\epsilon }(\zeta )\sin (b\zeta )d\zeta $, where
$|\int_{\alpha }^{\beta } [\psi (\zeta )-\psi _{\epsilon }(\zeta
)]\sin (b\zeta )d\zeta |<\epsilon /2$ for each $b$, since $|\sin
(b\zeta )|\le 1$ and the second term tends to zero: $\lim_{b\to
\infty }\int_{\alpha }^{\beta }\psi _{\epsilon }(\zeta )\sin (b\zeta
)d\zeta =0$ by the one proved above.
\par The final part of the {\bf Proof} of Theorem 19. For a fixed
$\epsilon >0$ there is an equality:
\par $\int_{-\infty }^{\infty }[w(\zeta +t)-w(t)]\zeta ^{-1}\sin
(b\zeta )d\zeta =\int_{-B}^B [w(\zeta +t)-w(t)]\zeta ^{-1}\sin
(b\zeta )d\zeta + \int_{|\zeta |>B}w(\zeta +t)\zeta ^{-1}\sin
(b\zeta )d\zeta - w(t)\int_{|\zeta |>B}\sin (b\zeta )\zeta
^{-1}d\zeta $. The second and the third terms are converging
integrals and therefore for sufficiently large $B>0$ by the absolute
value they are smaller than $\epsilon /3$. In view of the H\"older
condition $|[w(\zeta +t)-w(t)] \zeta ^{-1}|\le A|\zeta |^{1-c}$,
where $c>0$, $A>0$. Then in view of Lemma 10 there exists $b_0>0$
such that \par $|\int_{-B}^B [w(\zeta +t)-w(t)]\zeta ^{-1}\sin
(b\zeta )d\zeta |<\epsilon /3$ for each $b>b_0$. Thus,
\par $\lim_{b\to \infty } \int_{-\infty }^{\infty } [w(\zeta
+t)-w(t)]{\zeta }^{-1}\sin (b\zeta )d\zeta =0$.
\par This theorem in the case of a general function $M_N(p,t;\zeta )$
can be proved directly by the calculation of appearing integrals by
real variables $t$ and $\tau $ using Lemma 20.
\par {\bf 21. Theorem.} {\it An original $f(t)$ with $f({\bf R})\subset
\bf K$ for ${\bf K}=\bf H$ or ${\bf K}=\bf O$, or $f({\bf R})
\subset \bf R$ over the Cayley-Dickson algebra ${\cal A}_r$ with
$4\le r \in \bf N$ is completely defined by its image $\mbox{
}_NF_u(p;\zeta )$  up to values at points of discontinuity.}
\par {\bf Proof.} In view of theorem 19 the value $f(t)$
at each point $t$ of continuity of $f(t)$ is expressible throughout
$\mbox{ }_NF_u(p;\zeta )$ by Formula $19(i)$. At the same time
values of the original at points of discontinuity do not influence
on the image $\mbox{ }_NF_u(p;\zeta )$, since on each bounded
interval a number of points of discontinuity is finite.
\par {\bf 22. Theorem.} {\it If a function $\mbox{
}_NF_u(p)$ is analytic by the variable $p\in {\cal A}_r$ in the half
space $W := \{ p\in {\cal A}_r: Re (p)>s_0 \} $, where $2\le r\in
\bf N$, $f({\bf R})\subset \bf K$ for ${\bf K}=\bf H$ or ${\bf
K}=\bf O$ and also an arbitrary algebra ${\cal A}_r$ with $f({\bf
R})\subset \bf R$ for $4\le r<\infty $, $u(p,t;\zeta )=pt+\zeta $ or
$u(p,t;\zeta ) := p_0t + M(p,t;\zeta )+\zeta _0$, in the half space
$Re (p)>s_0$, moreover, for each $a>s_0$ there exist constants $C_a
>0$ and $\epsilon _a >0$ such that \par $(i)$ $|\mbox{ }_NF_u(p)|\le
C_a\exp (-\epsilon _a |p|)$ for each $p\in {\cal A}_r$ with $Re
(p)\ge a$, where $s_0$ is fixed, the integral \par $(ii)$
$\int_{a-S\infty }^{a+S\infty }\mbox{ }_NF_u(p)dp$ is absolutely
converging, then $\mbox{ }_NF_u(p)$ is the image of the function
\par $(iii)$ $f(t)=(2\pi )^{-1}{\tilde S}\int_{a-S\infty
}^{a+S\infty }\mbox{ }_NF_u(p)\exp (u(p,t;0))dp$.}
\par {\bf Proof.} The case $u(p,t;\zeta ) = pt + \zeta $ follows from
$u(p,t;\zeta ) := p_0t + M(p,t)+\zeta _0$ when
$p=(p_0,p_1,0,...,0)$, but the integral along the straight line
$a+St$, $t\in \bf R$, with such $p$ in the basis of generators
$(N_0,...,N_{2^r-1})$ can be obtained from the general integral by
an automorphism $v$, $z\mapsto v(z)$, of the algebra ${\cal A}_r$,
$r=2$ for ${\bf K}=\bf H$, $r=3$ for ${\bf K}=\bf O$. That is, in
view of Lemmas 16 and 17 it is sufficient to prove the equality of
the type $(iii)$ after the action of the automorphism $v$ on the
left and the right its parts.
\par Let $Re (p)=a>s_0$, then \par $|\int_{a-S\infty }^{a+S\infty }
\exp (M(p,t;0))\mbox{ }_NF_u(p)dp|$  $\le  \int_{-
\infty }^{\infty } |\mbox{ }_NF_u(a+S\theta )|d\theta $. \\
In view of the supposition of this theorem this integral converges
uniformly relative to $t\in \bf R$. For $f(t)$ given by the Formula
$(iii)$ for $Re (\eta ) =: \eta _0>s_0$, we get
\par $\int_0^{\infty }f(t)\exp (-\eta _0 t)dt $ \\  $= (2\pi
)^{-1}{\tilde S}\sum_{j=0}^{2^r-1} \int_0^{\infty }(\int_{a-S\infty
}^{a+S\infty }\mbox{
}_NF_{u,j}(p) \exp (u(p,t;0))dp)\exp (-\eta _0 t)(dt)N_j$, \\
in which it is possible to change the order of the integration,
since $t\in \bf R$. Then
\par $\int_0^{\infty }f(t)\exp
(-\eta _0 t)dt = (2\pi )^{-1}{\tilde S}\sum_{j=0}^{2^r-1}
\int_{a-S\infty }^{a+S\infty }(\int_0^{\infty }\mbox{ }_NF_{u,j}
(p)\exp ((a-\eta _0)t)dt)(dp)N_j$, \\
since $e^v\in \bf R$ for each $v\in \bf R$,
$e^{aM}e^{bM}=e^{(a+b)M}$ for each $a, b\in \bf R$. In view of
$a<\eta _0$ and
\par $\int_0^{\infty }e^{(a-\eta _0)t}dt= - (a-\eta _0)^{-1}$, \\
then \par $\int_0^{\infty }f(t)\exp (-\eta _0 t)dt = - (2\pi
)^{-1}{\tilde S}\sum_{j=0}^{2^r-1} (\int_{a-S\infty }^{a+S\infty
}\mbox{ }_NF_{u,j}(p)(a-\eta _0)^{-1}dp)N_j$ \par $= - (2\pi
)^{-1}{\tilde S} (\int_{a-S\infty }^{a+S\infty }\mbox{ }_NF_u(p
)(a-\eta _0)^{-1}dp$.
\par To finish the proof it is necessary the following analog
of the Jordan lemma.
\par {\bf 23. Lemma.} {\it Let a function $F$ of the variable $p$ from
the Cayley-Dickson algebra ${\cal A}_r$ with $2\le r\in \bf N$
satisfy Conditions $(1-3)$: \par $(1)$ the function $F(p)$ is
continuous by the variable $p\in {\cal A}_r$ in an open domain $W$
of the half space $\{ p\in {\cal A}_r: Re (p)>s_0 \} $, moreover for
each $a>s_0$ there exist constants $C_a
>0$ and $\epsilon _a >0$ such that
\par $(2)$ $|F(p)|\le C_a\exp (-\epsilon _a |p|)$ for each
$p\in C_{R(n)}$, $C_R := \{ z\in {\cal A}_r: Re (z)\ge a \} $,
$0<R(n)<R(n+1)$ for each $n\in \bf N$, $\lim_{n\to \infty
}R(n)=\infty $, where $s_0$ is fixed, the integral
\par $(3)$ $\int_{a-S\infty }^{a+S\infty }F(p)dp$ is absolutely
converging. Then
\par $(4)$ $\lim_{n\to \infty }\int_{\gamma _n} F(p)\exp
(-u(p,t;\zeta ))dp=0$ \\
for each $t>0$ and each sequence of rectifiable curves $\gamma _n$
contained in $C_{R(n)}\cap W$, moreover either $F(p)$ is holomorphic
in $W$, which is $(2^r-1)$-connected open domain with ${\cal A}_r$
(see \cite{span}), such that the projection $\pi _{s,p,t}(W)$ is
simply connected in $s{\bf R}\oplus p\bf R $ for each $s=i_{2k}$,
$p=i_{2k+1}$, $k=0,1,...,2^{r-1}-1$ for each $t\in {\cal
A}_{r,s,p}:={\cal A}_r\ominus s{\bf R}\ominus p\bf R$ and $u\in
s{\bf R}\oplus p\bf R$, for which there exists $z=t+u\in {\cal
A}_r$; or there exists a constant $C_V>0$ such that the variations
(lengths) of curves are bounded $V(\gamma _n)\le C_V R_n$ for each
$n$, where $u(p,t;\zeta )$ is the function given by the Equations
3(1-3,3'), $n\in \bf N$.}
\par {\bf Proof.} If $0< \epsilon \le \min
(a-s_0,\epsilon _a)$, then in view of Condition $(ii)$ there exists
a constant  $C>0$ such that
\par $(iv)$  $|F(p)|\le Ce^{-\epsilon |p|}$, \\
for each $p\in {\cal A}_r$ with $Re (p)\ge a>s_0$. If $U$ is a
domain in ${\cal A}_r$ of the diameter not greater than $\rho $,
then in view of (4) from the proof of Theorem 7 \cite{ludfov} there
is accomplished the inequality:
\par $\sup_{p\in U} \| {\hat F}(p) \| \le \sup_{p\in U} |F(p)|
C_1\exp (C_2\rho ^n)$, \\
where $C_1$ and $C_2$ are positive constants independent from $F$,
$n=2^r+2$, $2\le r\in \bf N$. In particular, as $U$ it is possible
to take the interior of the parallelepiped with ribs of lengths not
greater than $\rho /2^{r/2}$. Then the path of integration can be
covered by a finite number of such parallelepipeds. In the case of
the circle of radius $R$ a number of necessary parallelepipeds is
not greater, than $2\pi R/\rho +1$. There exists $R_0>0$ such that
for each $R>R_0$ there is accomplished the inequality $2\pi R/\rho
<\exp (C_2\rho ^{n-1}(R-\rho ))$. Therefore, in $\rho $ neighborhood
$C_R^{\rho }$ of the circle of radius $R$ with $R>R_0$ there is
accomplished the inequality:
\par $\sup_{p\in C_R^{\rho }} \| {\hat F}(p) \| \le
\sup_{p\in C_R^{\rho }} |F(p)| C_1\exp (C_2\rho ^{n-1}R)$.\\
Since $\rho >0$ can be taken arbitrary small, then there exists
$\rho _0>0$ such that for each $0<\rho <\rho _0$ there is
accomplished the inequality $C_2\rho ^{n-1}<\epsilon $,
consequently,
\par $\sup_{p\in C_R^{\rho }, Re (p)\ge a} \| {\hat F}(p) \| \le
C C_1\exp ((C_2\rho ^{n-1}-\epsilon )R)\le C C_1\exp (-\delta R)$ \\
in view of the condition imposed on $F$, where $C$ is the positive
constant for the given $F$, $\delta = C_2\rho ^{n-1}-\epsilon >0$.
With this the length of the path of integration does not exceed
$2\pi R$ and $\lim_{R\to \infty } CC_1 2\pi R \exp (-\delta R)=0$.
The function $F(p)$ is continuous by $p$, hence it is integrable
along each rectifiable curve in the domain $W$ of the half space $\{
p\in {\cal A}_r: Re (p)>s_0 \} $. \par If $F(p)$ is holomorphic,
then in view of Theorem 2.11 \cite{ludfov,ludoyst} $\int_{\gamma _n}
F(p)\exp (-u(p,t;\zeta ))dp$ is independent from the type of the
curve and it is defined only by the initial and final its points. If
$V(\gamma _n)\le C_V R_n$ for each $n$, then it is sufficient to
prove the statement of this Lemma for each subsequence $R_{n(k)}$
with $R_{n(k+1)}\ge R_{n(k)}+1$ for each $k\in \bf N$.Denote for the
simplicity such subsequence by $R_n$. Each rectifiable curve can be
approximated by the converging sequence of rectifiable of polygonal
type line composed of arcs of circles. If a curve is displayed on
the sphere, then these circles can be taken with the common center
with the sphere. Condition $(2)$ in each plane ${\bf R}\oplus N\bf
R$, where $N\in {\cal A}_r$, $Re (N)=0$, $|N|=1$, is accomplished,
moreover, uniformly relative to a directrix $N$ and it can be
accomplished a diffeomorphism $g$ in ${\cal A}_r$, such that
$g(W)=W$, $g(S_n)=S_n$ for each $n\in \bf N$, and an image of a
$C^1$ curve from $W$ is an arc of a circle, since $0<R_n+1<R_{n+1}$
for each $n\in \bf N$ and $\lim_{n\to \infty }R_n=\infty $. The
function $(F,\gamma )\mapsto \int_{\gamma }F(p)dp$ is continuous
from $C^0(V,{\cal A}_r)\times \Gamma $ into ${\cal A}_r$, where $V$
is the compact domain in ${\cal A}_r$, $\Gamma $ is the family of
rectifiable curves in $V$ with the metric $\rho (v,w):= \max
(\sup_{z\in v}\inf_{\zeta \in w}|z-\zeta |, \sup_{z\in w}\inf_{\zeta
\in v}|z-\zeta |)$ (see Theorem 2.7 \cite{ludfov,ludoyst}). The
space $C^1$ of all continuously differentiable functions of the real
variable is dense in the space of continuous functions $C^0$ in the
compact-open topology in the case of a finite number of variables.
In addition a rectifiable curve is an uniform limit of $C^1$ curves,
since each rectifiable curve is continuous. Therefore, consider
$\gamma _n=\psi _n\cap \{ p\in {\cal A}_r: Re (p)>a \} \cap W$.
Consequently,
\par $\lim_{n\to \infty }\int_{\gamma _n} F(p)\exp
(-u(p,t;\zeta ))dp=0$, \\
since this is accomplished for $\gamma _n=\pi _n\cap C_{R(n)}$ and
hence for general $\gamma _n$ with the same initial and final
points, where $\pi _n$ are two dimensional over $\bf R$ planes in
${\cal A}_r$.
\par The continuation of the {\bf Proof} of Theorem 22.
In view of Lemma 23
\par $ |\int_{\psi _R} F(p)(p-\eta _0)^{-1}dp|\le u(R)\pi R/(R-|\eta
_ 0|)$, \\
where $0<u(R)<\infty $ and there exists $\lim_{R\to \infty }
u(R)=0$, while $\psi _R$ is the arc of the circle $|p|=R$ in the
plane ${\bf R}\oplus S{\bf R}$ with $Re (p)>a$, consequently,
\par $ \lim_{R\to \infty } \int_{\psi _R} F(p)(p-\eta _0)^{-1}dp=0$, \\
since $u(R)\le u_0\exp (-\delta R)$ for each $R>R_0$, where
$u_0=const >0$. \par Then the straight line $a+S\theta $ with
$\theta \in \bf R$ can be substituted by the closed contour $\phi
_R$ composed from $\psi _R$ and the segment $[a+Sb,a-Sb]$ passed
from the above to the bottom. Thus,
\par $\int_0^{\infty }\exp (- \eta _0t)f(t)dt = (2\pi )^{-1}
{\tilde S}\int_{\phi _R}F(p)(p-\eta _0)^{-1}dp,$ \\
where the sign in front of the integral is changed due to the change
of the pass direction of the loop $\phi _R$. Recall, that in the
case of the  Cayley-Dickson algebra ${\cal A}_r$ the residue of a
function is the operator $\bf R$-homogeneous and ${\cal
A}_r$-additive by the argument $L\in {\cal A}_r$ with $Re (L)=0$,
where the residue is naturally dependent on a function and a point.
In the domain $ \{ p\in {\cal A}_r: Re (p)\ge a, |p| \le R \} $ the
analytic function $F(p)$ has only one point of singularity $p=\eta
_0$, which is the pole of the first order with the residue $res
(\eta _0; (p-\eta _0)^{-1}F(p)).L = L F (\eta _0)$ for each $L\in
{\cal A}_r$ with $Re (L)=0$, consequently,
\par $\int_0^{\infty }  f(t)\exp (-\eta _0t)dt = F(\eta _0)$, since
$L=S$ in the given case and $S{\tilde S}=1$.
\par For $t<0$ in view of the aforementioned ${\cal A}_r$ Lemma 23
we get, that \par $\lim_{R\to \infty }\int_{\phi _R}
F(p)e^{pt}dp=0$, \\
since $Re (p)=a>0$, consequently, the straight line $a+S\theta $,
$\theta \in \bf R$, can be substituted by the loop $\phi _R$ as
above. Then for $t<0$ we get:
\par $f(t) = (2\pi )^{-1}\int_{\phi _R}F(p)e^{pt}dp=0$, \\
since $F(p)$ is analytic by $p$ together with $e^{pt}$ in the
interior of the domain $\{ p: p\in {\cal A}_r; |p|\le R', Re (p)>s_0
\} $, $a>s_0$, $0<R<R'\le \infty $. Then the condition 2 for the
original is accomplished. On the other hand,
\par $|f(t)| \le (2\pi )^{-1} e^{at} \int_{-\infty }^{\infty }
|F(a+S\theta )|d\theta = Ce^{at}$, \\
where $C=(2\pi )^{-1}\int_{-\infty }^{\infty }|F(a+S\theta )|d\theta
<\infty $, consequently, Condition (3) is satisfied. As well as
$f(t)$ is continuous, since the function $F(p)$ in the integral is
continuous satisfying Conditions $(i,ii)$ and
\par $\lim_{R\to \infty } \int_{\gamma (\theta ): |\theta |\ge R }
F(p)dp=0$. Moreover, the integral \par $\int_{a-S\infty }^{a+S\infty
}\mbox{ }_NF_u(p) [\partial \exp
(u(p,t;0))/\partial t]dp$ \\
converges due to Conditions $(i,ii)$ and the proof above,
consequently, the function $f(t)$ is differentiable and hence
satisfies the Lipschitz condition.
\par {\bf 24. Note.} In Theorem 22 Condition (i) can be replaced on
\par $\lim_{n\to \infty }\sup_{p\in C_{R(n)}} \| {\hat F}(p) \|
=0$ \\
by the sequence $C_{R(n)} := \{ z\in {\cal A}_r: |z| =R(n), Re
(z)>s_0 \} $, where $R(n)<R(n+1)$ for each $n$, $\lim_{n\to \infty }
R(n)=\infty $, since this leads to the accomplishment of the ${\cal
A}_r$ analog of the Jordan Lemma for each $r\ge 2$ (see also Note
36). But in Theorem 22 itself it is essential the alternativity of
the algebra, that is, in it in general are possible only $r=2$ or
$r=3$ for $f({\bf R})\subset \bf K$, or under the condition $f({\bf
R})\subset \bf R$ for $4\le r\in \bf N$.
\par Further properties of quaternion, octonion and
general ${\cal A}_r$ analogs of the Laplace transform are
considered.
\par {\bf 25. Proposition.} {\it If there are $\mbox{
}_NF_u(p;\zeta )$ and $\mbox{ }_NG_u(p;\zeta )$ images \\
of functions-originals $f(t)$ and $g(t)$ in the half space $Re
(p)>s_0(f)$ and $Re (p)>s_0(g)$ with values in ${\cal A}_r$, then
for each $\alpha , \beta \in {\cal A}_r$ in the case ${\bf K}=\bf
H$; as well as $f$ and $g$ with values in $\bf R$ and each $\alpha ,
\beta \in {\cal A}_r$ or $f$ and $g$ with values in ${\cal A}_r$ and
each $\alpha , \beta \in \bf R$ in the case of ${\cal A}_r$; the
function $\alpha \mbox{ }_NF_u(p;\zeta ) + \beta \mbox{
}_NG_u(p;\zeta ) $ is the image of the function $\alpha f(t) +\beta
g(t)$ in the half space $Re (p)>\max (s_0(f),s_0(g))$.}
\par {\bf Proof.} Since there exist
$\mbox{ }_NF_u(p;\zeta )$ and $\mbox{ }_NG_u(p;\zeta )$, then the
integral
\par $\int_0^{\infty }(\alpha f(t)+\beta g(t)) \exp (-u(p,t))dt=
\int_0^{\infty }\alpha f(t) \exp (-u(p,t))dt+ \int_0^{\infty }\beta
g(t) \exp (-u(p,t))dt$ \\
converges in the half space $Re (p)>\max (s_0(f),s_0(g))$. At the
same time $t\in \bf R$, while $\bf R$ is the center of the algebra
${\cal A}_r$. The algebra $\bf H$ is associative. Therefore under
the indicated conditions  the constants $\alpha , \beta $ can be
carried out outside integrals.
\par {\bf 26. Examples.}  Let at first $u(p,t)=pt$, $p\in
{\cal A}_r$, $2\le r$.
\par 1. If  $f(t)=Ch_{[0,\infty )}(t)$
is the characteristic function of the ray $[0,\infty )$, then
$F(p)=p^{-1}$, that follows from the definition, since
$\int_0^{\infty }e^{-pt}dt=-p^{-1}e^{-pt}|_0^{\infty }=1/p$.
\par 2. If $f(t)=\sin (\omega t)Ch_{[0,\infty )}(t)$, then
$F(p)=\omega (p^2+\omega ^2 )^{-1}$, where $\omega \in \bf R$. \par
3. If $f(t)=\cos (\omega t)Ch_{[0,\infty )}(t)$, then
$F(p)=p(p^2+\omega ^2)^{-1}$. \par 4. For $f(t)=sh(\omega
t)Ch_{[0,\infty )}(t)$ the image is $F(p)=\omega (p^2-\omega
^2)^{-1}$.
\par 5. For $f(t)=ch (\omega t)Ch_{[0,\infty )}(t)$ the image is
$F(p)= p(p^2-\omega ^2)^{-1}$. The Laplace transforms in Examples
2-5 follow from the proposition 25 and the formula:
\par $\int_0^{\infty }\cos (\omega t)e^{-pt}dt = \int_0^{\infty
}[e^{-(p_0+S(p_S-\omega ))t}+ e^{-(p_0+S(p_S+\omega )t)}]dt/2=
[(p-\omega S)^{-1}+(p+\omega S)^{-1}]/2$, \\
where $p=p_0+Sp_S$, $p_0=Re (p)$, $S\in {\cal A}_r$, $Re (S)=0$,
$|S|=1$, $p_s\in \bf R$;
\par $\int_0^{\infty }\sin (\omega t)e^{-pt}dt =
\int_0^{\infty }[e^{-(p_0+S(p_S-\omega ))t} - e^{-(p_0+S(p_S+\omega
)t)}]{\tilde S}dt/2= [(p-\omega S)^{-1} - (p+\omega
S)^{-1}]{\tilde S}/2$, \\
\par $\int_0^{\infty }ch
(\omega t)e^{-pt}dt = \int_0^{\infty }[e^{-((p_0-\omega )+Sp_S)t}+
e^{-((p_0+\omega )+Sp_St)}]dt/2= [(p-\omega )^{-1}+(p+\omega
)^{-1}]/2$, \\
\par $\int_0^{\infty }sh
(\omega t)e^{-pt}dt = \int_0^{\infty }[e^{-((p_0-\omega )+Sp_S)t}-
e^{-((p_0+\omega )+Sp_St)}]dt/2= [(p-\omega )^{-1}-(p+\omega
)^{-1}]/2$, \\
since $\cos (\omega t)= [e^{\omega tS} + e^{-\omega tS}]/2$, $\sin
(\omega t)= [e^{\omega tS} - e^{-\omega tS}]{\tilde S}/2$, $ch
(\omega t)= [e^{\omega t} + e^{-\omega t}]/2$, $sh (\omega t)=
[e^{\omega t} - e^{-\omega t}]/2$, where the considered
functions-originals take only real values.
\par 6. Consider now $f(t)=\exp (\zeta t)Ch_{[0,\infty )}(t)$, then
\par $F(p)=\int_0^{\infty } \exp (\zeta t)\exp (-pt)dt$, \\
where  $\exp (\zeta t)\exp (-pt)=e^{(\zeta _0-p_0)t}\exp (\zeta
't)\exp (-p't)$, $\zeta _0 := Re (\zeta )$, $p_0 := Re (p)$, $\zeta
' := \zeta -\zeta _0$, $p' = p-p_0$. Let the element $N_1\in {\cal
A}_r$ be such that $\zeta ' = \zeta _1N_1$, where $Re (N_1)=0$,
$\zeta _1\in \bf R$, $|N_1|=1$, then $p' = p_1N_1+p_2N_2$, where
$N_2\perp N_1$ relative to the scalar product $Re (N_1{\tilde
N}_2)$, $p_1, p_2\in \bf R$, that is, $Re (N_1{\tilde N}_2)=0$, and
such that $Re (N_2)=0$ и $|N_2|=1$. Therefore, $\exp (\zeta 't)\exp
(-p't)= (\cos (\zeta _1t) +\sin (\zeta _1t)N_1) (\cos (|p'|t) - \sin
(|p'|t) (p_1N_1+p_2N_2)/|p'|)=[\cos (\zeta _1t)\cos (|p'|t)+ \sin
(\zeta _1 t)\sin (|p'|t)/|p'|] + [-\cos (\zeta _1t)\sin
(|p'|t)p_1/|p'|+ \sin (\zeta _1t)\cos (|p'|t))]N_1- [\cos(\zeta
_1t)\sin
(|p'|t)p_2/|p'|)N_2 - [\sin (|p'|t)\sin (\zeta _1t)p_2/|p'|]N_1N_2 $ \\
in view of Proposition 3.2 \cite{ludfov,ludoyst}, since $|\zeta
'|=|\zeta _1|$, the function $\cos (\phi )$ is even by $\phi \in \bf
R$, $\sin (|\phi |)\phi N/|\phi | =\sin (\phi )N$ for each $\phi \in
\bf R$ and $N\in {\cal A}_r$ with $Re (N)=0$, where $\lim_{\phi \to
0} \sin (\phi )/\phi =1$, and the function  $\sin (\phi )$ is odd;
with $N_1=\zeta '/|\zeta '|$ for $\zeta '\ne 0$, $N_2=[p' -
p_1N_1]/|p'-p_1N_1|$ for $p'\ne p_1N_1$ for some $p_1\in \bf R$,
$p_2N_2=p'-p_1N_1$. Since $Re ((p'-p_1N_1){\tilde {\zeta }} ')=0$,
then $p_1=Re (p'{\tilde {\zeta }}')/|\zeta '|$, because of ${\zeta
'}^2= - |\zeta '|^2$ and $N_1\zeta '= -|\zeta '|$. Then $F(p)=
2^{-1}\{ c[(c^2+a^2)^{-1}(1-p_1/|p'|) +(c^2+b^2)^{-1}(1+p_1 /|p'|)+
[a(c^2+a^2)^{-1}+b(c^2+b^2)^{-1})(1-p_1/|p'|)N_1 -
(a(c^2+a^2)^{-1}+b(c^2+b^2)^{-1})p_2N_2/|p'| - c((c^2+a^2)^{-1}-
(c^2+b^2)^{-1})(p_2/|p'|)N_1N_2 \} $,
\\ where $c:=\zeta _0-p_0$, $a:=|p'|+|\zeta _1|$, $b:=|p'| - |\zeta _1|$.
\par In this and the following examples there are used expressions
of $p_j$ throughout $p$ with the help of generators of the algebra
${\cal A}_r$ (see Formulas 3(6)). With $|z|= [z
(2^r-2)^{-1}(-z+\sum_{j=1}^{2^r-1}i_j(zi_j^*))]^{1/2}$ for each
$z\in {\cal A}_r$ this gives the local analyticity of $F(p;\zeta )$
by $p$ and $\zeta $.
\par {\bf 27. Theorem.} {\it Let $\alpha =const >0$, $F(p)$
is an image of a function $f(t)$ with $u=pt+\zeta $ or there is
given $u$ by Formulas 3(1-3,3') over ${\cal A}_r$ with $2\le
r<\infty $, then there is an image $F(p/ \alpha ;\zeta )/ \alpha $
of a function $f(\alpha t)$.}
\par {\bf Proof.} Since $p_jt+\zeta _j=p_j(\tau /\alpha )+
\zeta _j$ for each $j$, where $t\alpha =\tau $, then substituting
the variable gives: $\int_0^{\infty }f(\alpha t)e^{-u(p,t;\zeta
)}dt= \int_0^{\infty }f(\tau )e^{-u(p,\tau /\alpha ;\zeta )} d\tau
=F(p/\alpha ;\zeta )/\alpha $, since $\bf R$ is the center of the
algebra ${\cal A}_r$.
\par {\bf 28. Theorem.} {\it If there is a function-original $f(t)$,
$F(p)$ is its image, where $u=pt$, $f^{(n)}(t)$ is an original for
$n\ge 1$, then \par $(i)$ ${\cal F}(f^{(n)}(t)Ch_{[0,\infty
)}(t),pt;p;0) =
F(p)p^n-f(0)p^{n-1}-...-f^{(n-1)}(0)$ \\
is the image of the function $f^{(n)}(t)$ over ${\cal A}_r$ in the
domain $s_0<Re (p)$ for \\ $s_0=\max
(s_0(f),s_0(f'),...,s_0(f^{(n)}))$, where $f^{(k)}(0) := \lim_{t\to
+0}f^{(k)}(t)$, $f$ is real-valued for $r\ge 4$, $f$ is with values
in $\bf K$ for ${\bf K}=\bf H$ or ${\bf K}=\bf O$.}
\par {\bf Proof.} The integration by parts gives:
\par $(ii)$ $\int_0^{\infty }f'(t)e^{-pt}dt=[f(t)e^{-pt}]|_0^{\infty }
+(\int_0^{\infty }f(t)e^{-pt}dt)p$ \\
due to the alternativity of $\bf H$ and $\bf O$, for $r\ge 4$ due to
taking real values by $f$. Since $Re (p) = s>s_0$, then
$|f(t)e^{-pt}|\le Ce^{-(s-s_0)t}$, therefore, from Formula $(ii)$ it
follows Formula $(i)$ for $n=1$. The application of the latter
Formula $n$ times gives the general Formula $(i)$.
\par {\bf 29. Theorem.} {\it If there is a function-original $f(t)$,
$F(p)$ is its image over ${\cal A}_r$ with $2\le r<\infty $, where
$u=pt$, then $F^{(n)}(p).(\mbox{ }_
1h,...,\mbox{ }_nh)$ is the image of the original \\
$(-t)^nf(t)\mbox{ }_1h...\mbox{ }_nh$ for each $n\in \bf N$ and
$\mbox{ }_1h,...,\mbox{ }_nh\in {\bf R}\oplus p'{\bf R}\subset {\cal
A}_r$, where $p' := p-Re (p)$.}
\par {\bf Proof.} A function $F(p)$ in the half space
$Re (p)>s_0$ is holomorphic. Then $F'(p).h=-[\int_0^{\infty
}f(t)(-t)e^{-pt}dt]h$, since $h$ commutes with $p$ and with $t$. The
multiple differentiation gives $F^{(n)}(p).(\mbox{ }_1h,...,\mbox{
}_nh)=[\int_0^{\infty } (-t)^nf(t)e^{-pt}dt]\mbox{ }_1h,...,\mbox{
}_nh$.
\par {\bf 30. Examples.} 1. The image of the function-original
$t^nCh_{[0,\infty )}(t)$ is $F(p) = n! p^{-n-1}$; for $f(t)=t\sin
(\omega t)Ch_{[0,\infty )}(t)$ the image is equal to $F(p)=2p\omega
(p^2+\omega ^2)^{-2}$; for $f(t)= t\cos (\omega t)Ch_{[0,\infty
)}(t)$ the image is equal to $F(p)=(p^2- \omega ^2) (p^2+\omega
^2)^{-2}$ for $u=pt$ in view of Theorem 29 and examples 26.1-3,
since $\omega \in \bf R$, where $p\in {\cal A}_r$, $r\ge 2$.
\par 2.  For the function-original $f(t)=t^n\exp (\zeta t)Ch_{[0,\infty )}
(t)$ its image is equal to \\
$(-1)^n(\partial ^nF(p)/\partial p^n).(1,...,1)=
\partial ^n{\sf F}(a,b,c)/\partial c^n$, where in the latter equality
the function $F$ is expressed through variables $(a,b,c)$ as in
Example 15.6 and it is denotes by $\sf F$. In particular, for $n=1$,
the image is:
\par ${\cal F}(te^{\zeta t}Ch_{[0,\infty )}(t),pt;
p;0)=2^{-1} \{ (a^2-c^2)(a^2+c^2)^{-2}(1-p_1/|p'|) +
(b^2-c^2)(b^2+c^2)^{-1}(1 + p_1/|p'|) - 2c(a(a^2+c^2)^{-2} +b
(b^2+c^2)^{-2})(1-p_1/|p'|)N_1 +
2c(a(a^2+c^2)^{-2}+b(b^2+c^2)^{-2})(p_2/|p'|)N_2 -
((a^2-c^2)(a^2+c^2)^{-2}-(b^2-c^2)(b^2+c^2)^{-2})(p_2/|p'|)
N_1N_2 \} $, \\
since $(a^2+c^2)^{-1} = (2a)^{-1}S((c+aS)^{-1} - (c-aS)^{-1}),$
$c(a^2+c^2)^{-1}=2^{-1}((c+aS)^{-1}+(c-aS)^{-1})$ for each $a, c\in
\bf R$, $a^2+c^2>0$, $S\in {\cal A}_r$, $|S|=1$, $Re (S)=0$, where
$p, \zeta \in {\cal A}_r$.
\par In the general case
\par $\partial ^n [c(c^2+b^2)^{-1}]/\partial c^n=
n!(-1)^n[(c+Sb)^{-n-1}+(c-Sb)^{-n-1}]/2$  \\
$=n!(-1)^n(c^2+b^2)^{-n-1}\sum_{0\le k\le
[(n+1)/2]} {{n+1}\choose {2k}} (-1)^k c^{n+1-2k}b^{2k}$, \\
where $S\in {\cal A}_r$, $|S|=1$, $b\ne 0$, $c, b\in \bf R$;
\par $\partial ^n(c^2+b^2)^{-1}/\partial c^n=n!(-1)^n
Sb^{-1} [(c+Sb)^{-n-1}-(c-Sb)^{-n-1}]/2$ \\
$=n!(-1)^n(c^2+b^2)^{-n-1}\sum_{0\le k\le [n/2]} {{n+1}\choose
{2k+1}} (-1)^k c^{n-2k}b^{2k}$. \\
Therefore, the image of the function $f(t)=t^ne^{\zeta
t}Ch_{[0,\infty )}(t)$ is the function
\par $F(p)=n!2^{-1} (-1)^n \{
\sum_{0\le k\le [(n+1)/2]} {{n+1}\choose {2k}} (-1)^k
c^{n+1-2k}[a^{2k}(c^2+a^2)^{-n-1}(1-p_1/|p'|)$ \\  $+
b^{2k}(c^2+b^2)^{-n-1}(1+p_1|p'|)] + [N_1(1-p_1/|p'|) - N_2p_2/|p'|]
\sum_{0\le k\le [n/2]} {{n+1}\choose {2k+1}} (-1)^k$ \\ $
c^{n-2k} [a^{2k+1}(c^2+a^2)^{-n-1} + b^{2k+1}(c^2+b^2)^{-n-1}]$ \\
$ - [N_1N_2p_2/|p'|] [ \sum_{0\le k\le [(n+1)/2]} {{n+1}\choose
{2k}} (-1)^k c^{n+1-2k}[a^{2k}(c^2+a^2)^{-n-1}-
b^{2k}(c^2+b^2)^{-n-1}] \} $.
\par 3. Let now $f(t)=\exp (u(\zeta ,t;\alpha ))Ch_{[0,\infty )}(t)$
and find its image $F(p;\beta )$ for $u(p,t;\beta ) := p_0t +
M(p,t;\beta )+\beta _0$, where $p, \zeta, \alpha , \beta \in {\cal
A}_r$ as in Definition 3. For this consider the trigonometric
Formulas:
\par $(i)$  $\cos (a-b)+\cos (a+b)=2\cos (a)\cos (b)$;
\par $(ii)$  $\cos (a-b) - \cos (a+b)=2\sin (a)\sin (b)$;
\par $(iii)$  $\sin (a+b)+\sin (a-b)=2\sin (a)\cos (b)$;
\par $(iv)$  $\sin (a+b)-\sin (a-b) = 2\sin (b)\cos (a)$ \\
for each $a, b\in \bf R$. Then
\par $(1)$ $\prod_{p=1}^n\cos (a_p)=2^{1-n}\sum_{(v_2,...,v_n
\in \{ 1, 2 \} )} \cos (a_1+(-1)^{v_2}a_2+(-1)^{v_3}a_3+...
+(-1)^{v_n}a_n)$ \\
for each $a_1,...,a_n\in \bf R$, $2\le n\in \bf N$;
\par $(2)$ $\prod_{p=1}^{2m}\sin (b_p)=2^{-m}\prod_{p=1}^m[\cos
(b_{2p-1}-b_{2p})-\cos (b_{2p-1}+b_{2p})]=$ \\
$(-2)^m\sum_{(w_1,...,w_m\in \{ 1,2 \} )}(-1)^{w_1+...+w_m}
(\prod_{p=1}^m\cos (b_{2p-1}+(-1)^{w_p}b_{2p}))$ \\
$=2^{1-2m}\sum_{(w_1,...,w_m\in \{ 1,2 \} )}(-1)^{w_1+...+w_m+m}
\sum_{(v_2,...,v_m\in \{ 1,2 \} )} \cos
(b_1+(-1)^{w_1}b_2+(-1)^{v_2}(b_3+$ \\ $(-1)^{w_2}b_4)+(-1)^{v_3}
(b_5+(-1)^{w_3}b_6)+...+(-1)^{v_m}(b_{2m-1}+(-1)^{w_m}b_{2m}))$ \\
for each $b_1,b_2,...\in \bf R$, $1\le m\in \bf N$;
\par $(3)$ $\prod_{p=1}^{2m+1} \sin (b_p)$
$=2^{-2m}\sum_{(w_1,...,w_m\in \{ 1,2 \} )}(-1)^{w_1+...+w_m+m}
\sum_{(v_1,...,v_m\in \{ 1,2 \} )} \sin (b_{2m+1}+$ \\
$(-1)^{v_1}(b_1+(-1)^{w_1}b_2)+(-1)^{v_2}(b_3+(-1)^{w_2}b_4)+
...+(-1)^{v_m}(b_{2m-1}+(-1)^{w_m}b_{2m}))$;
\par $(4)$ $(\prod_{p=1}^n\cos (a_p)) (\prod_{q=1}^{2m}\sin (b_q))$
$=2^{1-n-2m} \sum_{(v_2,...,v_n\in \{ 1,2 \} )}
\sum_{(w_1,...,w_m\in \{ 1,2 \} )} (-1)^{w_1+...+w_m+m}$ \\
$\sum_{(u_1,...,u_m\in \{ 1,2 \} )} \cos (a_1+(-1)^{v_2}a_2+...+
(-1)^{v_n}a_n+(-1)^{u_1}(b_1+(-1)^{w_1}b_2)+
(-1)^{u_2}(b_3+(-1)^{w_2}b_4)+...+(-1)^{u_m}(b_{2m-1}+(-1)^{w_m}b_{2m}))$;
\par $(5)$ $(\prod_{p=1}^n\cos (a_p)) (\prod_{q=1}^{2m+1}\sin (b_q))$
$=2^{-n-2m} \sum_{(v_1,...,v_n\in \{ 1,2 \} )} \sum_{(w_1,...,w_m\in
\{ 1,2 \} )} (-1)^{w_1+...+w_m+m}$ \\ $\sum_{(u_1,...,u_m\in \{ 1,2
\} )} \sin (b_{2m+1}+(-1)^{v_1}a_1+(-1)^{v_2}a_2+...+(-1)^{v_n}a_n+
(-1)^{u_1}(b_1+(-1)^{w_1}b_2)+
(-1)^{u_2}(b_3+(-1)^{w_2}b_4)+...+(-1)^{u_m}(b_{2m-1}+(-1)^{w_m}b_{2m}))$.
Moreover,
\par $(6)$  $\exp (M(p,t;\alpha )) = \cos (p_1t+\alpha _1) + i_1\sin
(p_1t+\alpha _1)\cos (p_2t+\alpha _2) + i_2\sin (p_1t+\alpha _1)
\sin (p_2t+\alpha _2)\cos (p_3t+\alpha _3) + i_3\sin (p_1t+\alpha
_1)\sin (p_2t+\alpha _2)\sin (p_3t+\alpha _3)$ for quaternions;
\par $(7)$ $\exp (M(p,t;\alpha ))=
\cos (p_1t+\alpha _1) + i_1\sin (p_1t+\alpha _1)\cos (p_2t+\alpha
_2) + i_2\sin (p_1t+\alpha _1) \sin (p_2t+\alpha _2)\cos
(p_3t+\alpha _3) +...+ i_6\sin (p_1t+\alpha _1)...\sin (p_6t+\alpha
_6)\cos (p_7t+\alpha _7) + i_7\sin (p_1t+\alpha _1)...\sin
(p_6t+\alpha _6)\sin (p_7t+\alpha _7)$ for octonions. Then \\
$(8)$  $\exp (M(\zeta ,t;\alpha )) \exp (- M(p,t;\beta )) = [\cos
(\zeta _1t+\alpha _1)\cos (p_1t+\beta _1) + \sin (\zeta _1t+\alpha
_1)\cos (\zeta _2t+\alpha _2)\sin (p_1t+\beta _1)\cos (p_2t+\beta
_2) + \sin (\zeta _1t+\alpha _1) \sin (\zeta _2t+\alpha _2)\cos
(\zeta _3t+\alpha _3)\sin (p_1t+\beta _1) \sin (p_2t+\beta _2)\cos
(p_3t+\beta _3) + \sin (\zeta _1t+\alpha _1)\sin (\zeta _2t+\alpha
_2)\sin (\zeta _3t+\alpha _3) \sin (p_1t+\beta _1)\sin
(p_2t+\beta _2)\sin (p_3t+\beta _3)]$ \\
$i_1[- \cos (\zeta _1t + \alpha _1)\sin (p_1t+\beta _1)\cos
(p_2t+\beta _2) + \sin (\zeta _1t+\alpha _1)\cos (\zeta _2t+\alpha
_2)\cos (p_1t + \beta _1) - \sin (\zeta _1t+\alpha _1) \sin (\zeta
_2t+\alpha _2)\cos (\zeta _3t+\alpha _3)\sin (p_1t+\beta _1)\sin
(p_2t+\beta _2)\sin (p_3t+\beta _3) + \sin (\zeta _1t+\alpha _1)\sin
(\zeta _2t+\alpha _2)\sin (\zeta _3t+\alpha _3)\sin (p_1t+\beta _1)
\sin (p_2t+\beta _2)\cos (p_3t+\beta _3)]$ \\
$i_2[- \cos (\zeta _1t+\alpha _1)\sin (p_1t+\beta _1) \sin
(p_2t+\beta _2)\cos (p_3t+\beta _3) + \sin (\zeta _1t+\alpha _1)
\sin (\zeta _2t+\alpha _2)\cos (\zeta _3t+\alpha _3)\cos (p_1t+\beta
_1) - \sin (\zeta _1t+\alpha _1)\sin (\zeta _2t+\alpha _2)\sin
(\zeta _3t+\alpha _3)\sin (p_1t+\beta _1)\cos (p_2t+\beta _2) + \sin
(\zeta _1t+\alpha _1)\cos (\zeta _2t+\alpha _2) \sin (p_1t+\beta
_1)\sin (p_2t+\beta _2)\sin (p_3t+\beta _3)]$ \\
$i_3[- \cos (\zeta _1t+\alpha _1)\sin (p_1t+\beta _1)\sin
(p_2t+\beta _2)\sin (p_3t+\beta _3) + \sin (\zeta _1t+\alpha _1)\sin
(\zeta _2t+\alpha _2)\sin (\zeta _3t+\alpha _3)\cos (p_1t+\alpha _1)
- \sin (\zeta _1t+\alpha _1)\cos (\zeta _2t+\alpha _2) \sin
(p_1t+\beta _1) \sin (p_2t+\beta _2)\cos (p_3t+\beta _3) + \sin
(\zeta _1t+\alpha _1) \sin (\zeta _2t+\alpha _2)\cos (\zeta
_3t+\alpha _3)\sin (p_1t+\beta _1)\cos (p_2t+\beta _2)]$.
\par For each $c>0$ with:
\par $(9)$ $\int_0^{\infty }e^{-ct}(\prod_{s=1}^n\cos (\alpha
_st))dt = 2^{1-n}c\sum_{v_2,...,v_n\in \{ 1, 2 \} }[c^2+(\alpha
_1+(-1)^{v_2}\alpha _2+...+(-1)^{v_n}\alpha _n)^2]^{-1}$;
\par $(10)$ $\int_0^{\infty }e^{-ct}(\prod_{s=1}^{2m}\sin (\beta
_st))dt = 2^{1-2m}c \sum_{w_1,...,w_m, v_2,...,v_m\in \{ 1, 2 \} }
(-1)^{w_1+...+w_m+m} [c^2 $ \\  $+ ((\beta _1+(-1)^{w_1}\beta _2)
+(-1)^{v_2}(\beta _3+(-1)^{w_2}\beta _4)+...+(-1)^{v_m}(\beta
_{2m-1}+(-1)^{w_m}\beta _{2m}))^2]^{-1}$;
\par $(11)$ $\int_0^{\infty }e^{-ct}(\prod_{s=1}^{2m+1}\sin (\beta
_st))dt = 2^{-2m} \sum_{w_1,...,w_m, v_1,...,v_m\in \{ 1, 2 \} }
(-1)^{w_1+...+w_m+m} (\beta _{2m+1}$ \\  $+ (-1)^{v_1}(\beta
_1+(-1)^{w_1}\beta _2) +(-1)^{v_2}(\beta _3+(-1)^{w_2}\beta
_4)+...+(-1)^{v_m}(\beta _{2m-1}+(-1)^{w_m}\beta _{2m}) [c^2 +
(\beta _{2m+1}$ \\  $+ (-1)^{v_1}(\beta _1+(-1)^{w_1}\beta _2)
+(-1)^{v_2}(\beta _3+(-1)^{w_2}\beta _4)+...+(-1)^{v_m}(\beta
_{2m-1}+(-1)^{w_m}\beta _{2m}))^2]^{-1}$;
\par $(12)$ $\int_0^{\infty }e^{-ct}(\prod_{s=1}^n\cos (\alpha
_st))(\prod_{q=1}^{2m}\sin (\beta _qt))dt= 2^{1-n-2m}c
\sum_{v_2,...,v_n, w_1,...,w_m, u_1,...,u_m \in \{ 1, 2 \} }
[c^2+(\alpha _1+(-1)^{v_2}\alpha _2+...+(-1)^{v_n}\alpha
_n+(-1)^{u_1}(\beta _1+(-1)^{w_1}\beta _2)+...+(-1)^{u_m}(\beta
_{2m-1}+(-1)^{w_m}\beta _{2m}))^2]^{-1}$;
\par $(13)$ $\int_0^{\infty }e^{-ct}(\prod_{s=1}^n\cos (\alpha
_st))(\prod_{q=1}^{2m+1}\sin (\beta _qt))dt= 2^{-n-2m}c
\sum_{v_1,...,v_n, w_1,...,w_m, u_1,...,u_m \in \{ 1, 2 \} }$ \\
$(-1)^{w_1+...+w_m+m} ((-1)^{v_1}\alpha _1+(-1)^{v_2}\alpha
_2+...+(-1)^{v_n}\alpha _n+\beta _{2m+1}+ (-1)^{u_1}(\beta
_1+(-1)^{w_1}\beta _2)+...+(-1)^{u_m}(\beta _{2m-1}+(-1)^{w_m}\beta
_{2m})) [c^2+((-1)^{v_1}\alpha _1+(-1)^{v_2}\alpha
_2+...+(-1)^{v_n}\alpha
_n$\\
$+\beta _{2m+1}+ (-1)^{u_1}(\beta _1+(-1)^{w_1}\beta
_2)+...+(-1)^{u_m}(\beta _{2m-1}+(-1)^{w_m}\beta _{2m}))^2]^{-1}$.
\par Then from Formula $(8-13)$ it follows, that the image
${\cal F}(f,u; p;0)$ of the function-original $f(t)=\exp (u(\zeta
,t;0))Ch_{[0,\infty )}(t)$ with $u(p,t;0) := p_0t + M(p,t;0)$ is:
\par $(14)$ ${\cal F}(f,u; p;0)= \{ (c/2) \sum_{v_1\in \{ 1,2 \} }
[c^2+(p_1+(-1)^{v_1}\zeta _1)^2]^{-1}+ (c/8)\sum_{v_1,w_1,u_1\in \{
1, 2 \} } (-1)^{w_1+1}[c^2+ (p_2+(-1)^{v_1}\zeta
_2+(-1)^{u_1}(p_1+(-1)^{w_1}\zeta _1))^2]^{-1}$ \\
$+(c/32) \sum_{v_1, w_1,w_2, u_1,u_2\in \{ 1,2 \}
}(-1)^{w_1+w_2}[c^2+(p_3+(-1)^{v_1}\zeta
_3+(-1)^{u_1}(p_1+(-1)^{w_1}\zeta _1)+
(-1)^{u_2}(p_2+(-1)^{w_2}\zeta _2))^2]^{-1}$ \\
$+(c/32)\sum_{w_1,w_2,w_3,u_1,u_2\in \{ 1, 2 \}
}(-1)^{w_1+w_2+w_3+1}[c^2+((p_1+(-1)^{w_1}\zeta
_1)+(-1)^{u_1}(p_2+(-1)^{w_2}\zeta
_2)+(-1)^{u_2}(p_3+(-1)^{w_3}\zeta _3))^2]^{-1} \} $ \\
$+i_1 \{ -4^{-1}\sum_{v_1,v_2\in \{ 1, 2 \} } [(\zeta
_1+(-1)^{v_1}p_1+(-1)^{v_2}\zeta _2)[c^2+ (\zeta
_1+(-1)^{v_1}p_1+(-1)^{v_2}\zeta
_2)^2]^{-1}-(p_1+(-1)^{v_1}p_2+(-1)^{v_2}\zeta _1)[c^2+
(p_1+(-1)^{v_1}p_2+(-1)^{v_2}\zeta _1)^2]^{-1}] - (32)^{-1}
\sum_{v_1,w_1,w_2,u_1,u_2\in \{ 1, 2 \} } (-1)^{w_1+w_2}[ (\zeta
_3+(-1)^{v_1}p_3+(-1)^{u_1}(p_1+(-1)^{w_1}\zeta _1)+(-1)^{u_2}
(p_2+(-1)^{w_2}\zeta _2))[c^2+(\zeta
_3+(-1)^{v_1}p_3+(-1)^{u_1}(p_1+(-1)^{w_1}\zeta _1)+(-1)^{u_2}
(p_2+(-1)^{w_2}\zeta _2))^2]^{-1} -(p_3+(-1)^{v_1}\zeta
_3+(-1)^{u_1}(p_1+(-1)^{w_1}\zeta _1)+(-1)^{u_2}(p_2+(-1)^{w_2}\zeta
_2)) [c^2+(p_3+(-1)^{v_1}\zeta _3+(-1)^{u_1}(p_1+(-1)^{w_1}\zeta
_1)+(-1)^{u_2}(p_2+(-1)^{w_2}\zeta _2))^2]^{-1}] \} $ \\
$+i_2 \{ (c/8)\sum_{v_1,w_1,u_1\in \{ 1, 2 \} } (-1)^{w_1+1}
[[c^2+(p_1+(-1)^{v_1}\zeta _3+(-1)^{u_1}(\zeta _1+(-1)^{w_1}\zeta
_2))^2]^{-1} - [c^2+(\zeta
_1+(-1)^{v_2}p_3+(-1)^{u_1}(p_1+(-1)^{w_1}p_2))^2]^{-1}] +(c/16)
\sum_{w_1,w_2,u_1,u_2\in \{ 1, 2 \} } (-1)^{w_1+w_2}[
-[c^2+(p_2+(-1)^{u_1}(p_1+(-1)^{w_1}\zeta _1)+(-1)^{u_2}(\zeta
_2+(-1)^{w_2}\zeta _3))^2]^{-1} + [c^2+(\zeta
_2+(-1)^{u_1}(p_1+(-1)^{w_1}\zeta
_1)+(-1)^{u_2}(p_2+(-1)^{w_2}p_3))^2]^{-1}] \} $ \\
$+ i_3 \{ 8^{-1}\sum_{v_1,w_1,u_1\in \{ 1, 2 \} } (-1)^{w_1+1}[
(\zeta _3+(-1)^{v_1}p_1+(-1)^{u_1}(\zeta _1+(-1)^{w_1}\zeta _2))
[c^2+ (\zeta _3+(-1)^{v_1}p_1+(-1)^{u_1}(\zeta _1+(-1)^{w_1}\zeta
_2))^2]^{-1} -(p_3+(-1)^{v_1}\zeta _1+(-1)^{u_1}(p_1+(-1)^{w_1}p_2))
[c^2+(p_3+(-1)^{v_1}\zeta _1+(-1)^{u_1}(p_1+(-1)^{w_1}p_2))^2]^{-1}]
+(16)^{-1}\sum_{v_1,v_2,w_1,u_1\in \{ 1, 2 \} } (-1)^{w_1}[
(p_2+(-1)^{v_1}\zeta _2+(-1)^{v_2}p_3+(-1)^{u_1}(p_1+(-1)^{w_1}\zeta
_1))[c^2+(p_2+(-1)^{v_1}\zeta
_2+(-1)^{v_2}p_3+(-1)^{u_1}(p_1+(-1)^{w_1}\zeta _1))^2]^{-1} -
(\zeta _2+(-1)^{v_1}p_2+(-1)^{v_2}\zeta
_3+(-1)^{u_1}(p_1+(-1)^{w_1}\zeta _1)) [c^2+(\zeta
_2+(-1)^{v_1}p_2+(-1)^{v_2}\zeta _3+(-1)^{u_1}(p_1+(-1)^{w_1}\zeta
_1))^2]^{-1}] \} $, \\
where $c=\zeta _0-p_0$. For non zero intial phases the formula for
${\cal F}(f,u; p;\beta )$ the function-original $f(t)=\exp (u(\zeta
,t;\alpha ))Ch_{[0,\infty )}(t)$ при $u(p,t;\alpha ) := p_0t +
M(p,t;\alpha )+\alpha _0$ is analogous, but it is more complicated
and it can be written with the use of the equalities:
\par $(15)$ $\int_0^{\infty }e^{-ct}\cos (at+b)dt=
[c\cos (b) -a\sin (b)] [c^2+a^2]^{-1}$,\\
\par $(16)$ $\int_0^{\infty }e^{-ct}\sin (at+b)dt=
[c\sin (b) + a \cos (b)] [c^2+a^2]^{-1}$\\
for each $c>0$, $a, b\in \bf R$, as well as Formulas $(8-13)$ also.
\par {\bf 31. Theorem.} {\it Let $f(t)$ be a function-original such that
$f'(t)$ also satisfies Conditions 1(1-3), $u(p,t) := p_0t +
M(p,t;\zeta )+\zeta _0$ over ${\cal A}_r$ с $2\le r<\infty $ (see
Definition 3). Then
\par $(1)$  ${\cal
F}(f'(t)Ch_{[0,\infty )}(t),u; p;\zeta ) = -f(0)+p_0{\cal
F}(f(t)Ch_{[0,\infty )}(t),u; p;\zeta ) $ \\  $+p_1{\cal
F}(f(t)Ch_{[0,\infty )}(t),u; p;\zeta - i_1\pi
/2)+...+p_{2^r-1}{\cal
F}(f(t)Ch_{[0,\infty )}(t),u; p;\zeta - i_{2^r-1}\pi /2)$ \\
in the domain $s_0<Re (p)$ for $s_0=\max (s_0(f),s_0(f'))$.}
\par {\bf Proof.} From Equations 30(6,7) it follows, that
the equality is satisfied:
\par $(2)$  $\partial \exp (-u(p,t;\zeta ))/\partial t= - p_0\exp
(-u(p,t;\zeta )) - p_1\exp (-u(p,t;\zeta -i_1\pi /2)) -...$ \\
$-p_{2^r-1} \exp (-u(p,t;\zeta -i_{2^r-1}\pi /2))$, since
\par $\exp
(-u(p,t;\zeta ))=\exp (-p_0t-\zeta _0)\exp (-M(p,t;\zeta ))$,
\par $\partial \exp (-p_0t-\zeta _0)/\partial t= -p_0\exp (-p_0t-\zeta
_0)$, \par  $\partial [\cos (p_jt+\zeta _j)-\sin (p_jt+\zeta
_j)i_j]/\partial t =\partial \exp (-(p_jt+\zeta _j)i_j)/
\partial t$  \\  $ = -p_ji_j\exp (-(p_jt+\zeta _j)i_j)=
-p_j\exp (-(p_jt+\zeta _j-\pi /2)i_j)$ \\  $=-p_j[\cos (p_jt+\zeta
_j-\pi /2) - \sin (p_jt+\zeta _j-\pi /2)i_j]$,\\ with this the
integration by parts gives \par $\int_0^{\infty }f'(t)\exp
(-u(p,t;\zeta ))dt = f(t)\exp (-u(p,t;\zeta ))|_0^{\infty } -
\int_0^{\infty }[f(t)\partial \exp (-u(p,t;\zeta ))/\partial t]dt$.
\par {\bf 32. Theorem.} {\it Let $f(t)$ be a function-original,
$u(p,t) := p_0t + M(p,t;\zeta )+\zeta _0$ над ${\cal A}_r$ with
$2\le r<\infty $ (see Definition 3). Then
\par $(1)$  $(\partial {\cal F}(f(t)Ch_{[0,\infty )}(t),u; p;\zeta )/
\partial p).h =
- {\cal F}(f(t)tCh_{[0,\infty )}(t),u;p;\zeta )h_0 $ \\  $- {\cal
F}(f(t)tCh_{[0,\infty )}(t),u; p;\zeta - i_1\pi /2)h_1 -...-{\cal
F}(f(t)tCh_{[0,\infty )}(t),u; p;
\zeta - i_{2^r-1}\pi /2)h_{2^r-1}$\\
for each $h=h_0i_0+...+h_{2^r-1}i_{2^r-1}\in {\cal A}_r$, where
$h_0,...,h_{2^r-1}\in {\bf R}$, $Re (p)>s_0$.}
\par {\bf Proof.} From Equations 30(6,7) it follows, that
\par $(2)$  $(\partial \exp (-u(p,t;\zeta ))/\partial p).h= - p_0\exp
(-u(p,t;\zeta ))h_0 - \exp (-u(p,t;\zeta -i_1\pi /2))h_1 -... -\exp
(-u(p,t;\zeta -i_{2^r-1}\pi /2))h_{2^r-1}$. \\
Since ${\cal F}(f(t)Ch_{[0,\infty )}(t),u; p;\zeta )$ is the
holomorphic function by $p$ for $Re (p)>s_0$, where $f(t)$ satisfies
conditions of Definition 1, $|\int_0^{\infty }e^{-ct}t^ndt|<\infty $
for each $c>0$ and $n=0,1,2,...$, then it can be differentiated
under the sign of the integral:
\par $(\partial (\int_0^{\infty }f(t)\exp (-u(p,t;\zeta
))dt)/\partial p).h =\int_0^{\infty }f(t)(\partial \exp
(-u(p,t;\zeta ))/\partial p).hdt$.
\par {\bf 33. Example.} Let $f_n(t) = t^nCh_{[0,\infty )}(t)$ be
the function-originals, where $n=0,1,2,...$, for $u(p,t) := p_0t +
M(p,t;\zeta )+\zeta _0$, then
\par $(1)$  ${\cal F}(f_n(t),u; p;\zeta
) = \int_0^{\infty }t^n\exp (-p_0t-\zeta _0) \cos (p_1t+\zeta
_1)dt]$
\\ $+\sum_{k=1}^{2^r-2} [\int_0^{\infty }t^n\exp (-p_0t-\zeta _0)
\sin(p_1t+\zeta _1)...\sin(p_kt+\zeta _k) \cos (p_{k+1}t+\zeta
_{k+1})dt]i_k$ \\ $+ [\int_0^{\infty }t^n\exp (-p_0t-\zeta
_0)\sin(p_1t+\zeta _1)...\sin(p_{2^r-2}t+\zeta _{2^r-2}) \sin
(p_{2^r-1}t+\zeta _{2^r-1})dt]i_{2^r-1}$, \\
where $r=2$ for ${\bf K}=\bf H$ and $r=3$ for ${\bf K}=\bf O$, $ \{
i_0, i_1,...,i_{2^r-1} \} $ are the standard generators of the
algebra ${\cal A}_r$, $2\le r\in \bf N$. Consider the integral:
\par $(2)$ $J_n(\alpha ) := \int_0^{\infty }t^ne^{-\alpha t}dt$ \\
for each $\alpha \in \bf C$ with $Re (\alpha )>0$. In particular,
\par $(3)$ $J_0(\alpha )=- \exp (-\alpha t) (\alpha
)^{-1}|_0^{\infty }= \alpha ^{-1}$. Therefore,
\par $(4)$  $J_n(\alpha ) = (-1)^n dJ_0 (\alpha )/d\alpha
^n = n!\alpha ^{-n-1}$, consequently,
\par $(5)$ $\int_0^{\infty }t^ne^{-\alpha t- b}dt=e^{-b}J_n(\alpha
)$ for each $b\in \bf C$, from which it follows, that
\par $(6)$ $\int_0^{\infty }t^n\exp ^{-\alpha _0t}\cos (\alpha
_1t+\beta )dt = Re (\exp (-{\bf i}\beta )J_n(\alpha )) = n! Re (\exp
(-{\bf i}\beta ) \alpha ^{-n-1})$,
\par $(7)$ $\int_0^{\infty }t^n\exp ^{-\alpha _0t}\sin (\alpha
_1t+\beta )dt = - Im (\exp ( - {\bf i}\beta ) J_n(\alpha )) =
n! Im (\exp ( - {\bf i}\beta ) \alpha ^{-n-1})$ \\
for each $\beta \in \bf R$, where $Re (\alpha )=\alpha _0=(\alpha +
{\bar {\alpha }})/2$, $Im (\alpha )=\alpha _1=(\alpha -{\bar {\alpha
}})/(2{\bf i})$, ${\bf i}=i_1=(-1)^{1/2}$, $\alpha =\alpha _0+\alpha
_1{\bf i}$. From the equations
\par $(8)$  $\alpha ^{-k} = {\bar {\alpha
}}^k |\alpha |^{-2k}$ and
\par $(9)$  $\alpha ^{-k} = |\alpha |^{-2k}
\sum_{q=0}^{[k/2]} {k\choose {2q}} (-1)^q\alpha _0^{k-2q}\alpha
_1^{2q} - {\bf i}\sum_{q=0}^{[(k-1)/2]} {k\choose {2q+1}}
(-1)^q\alpha _0^{k-2q-1}\alpha _1^{2q+1}$, \\ where ${k\choose q} =
k!/[(k-q)!q!]$ are the binomial coefficients for each $0\le q\le
k\in \bf Z$ it follows, that
\par $(10)$  $Re (\exp (-{\bf i}\beta )J_n(\alpha )) =
n! [\cos (\beta )\sum_q {{n+1}\choose {2q}} (-1)^q\alpha _0^{n-2q+1}
\alpha _1^{2q}$ \\  $ - \sin (\beta )\sum_q {{n+1}\choose {2q+1}}
(-1)^{q}\alpha _0^{n-2q} \alpha _1^{2q+1}] (\alpha _0^2+\alpha
_1^2)^{-n-1} =: T_n(\alpha _0,\alpha _1,\beta )$,\\
\par $(11)$  $Im (\exp (-{\bf i}\beta )J_n(\alpha )) =
- (n!) [\sin (\beta )\sum_q {{n+1}\choose {2q}} (-1)^q\alpha
_0^{n-2q+1} \alpha _1^{2q}$ \\  $+ \cos (\beta )\sum_q {{n+1}\choose
{2q+1}} (-1)^{q}\alpha _0^{n-2q} \alpha _1^{2q+1}] (\alpha
_0^2+\alpha _1^2)^{-n-1} =: - S_n(\alpha _0,\alpha _1,\beta )$. \\
Then from Equations $(1-11)$ it follows, that for ${\bf K}=\bf H$
\par $(12)$ ${\cal F}(f_n(t)Ch_{[0,\infty )}(t),u; p;\zeta
) = e^{-\zeta _0} \{ T_n(p_0,p_1,\zeta _1)$ \\  $ -
(i_1/2)\sum_{v_1\in \{
1,2 \} } S_n(p_0,p_1+(-1)^{v_1}p_2,\zeta _1+(-1)^{v_1}\zeta _2)$ \\
$-(i_2/4)\sum_{v_1,u_1\in \{ 1,2 \} } (-1)^{v_1}
T_n(p_0,p_3+(-1)^{u_1}(p_1+(-1)^{v_1}p_2),\zeta _3 +
(-1)^{u_1}(\zeta _1 + (-1)^{v_1}\zeta _2))$ \\
$- (i_3/4)\sum_{v_1, u_1\in \{ 1,2 \} } (-1)^{v_1+1}S_n(p_0,
p_3+(-1)^{u_1}(p_1+(-1)^{v_1}p_2), \zeta _3+(-1)^{u_1}(\zeta
_1+(-1)^{v_1}\zeta _2)) \} $, \\
as well as over the algebra ${\cal A}_r$ with $r\ge 3$ the image is
given by the formula:
\par $(13)$ ${\cal F}(f_n(t)Ch_{[0,\infty )}(t),u; p;\zeta
) = e^{-\zeta _0} \{ T_n(p_0,p_1,\zeta _1)$ \\  $ -
\sum_{s=1}^{2^{r-1}-1} i_{2s} 2^{-2s} \sum_{v_1,...,v_s,
u_1,...,u_s\in \{ 1,2 \} } (-1)^{v_1+...+v_s+s}
T_n(p_0,p_{2s+1}+(-1)^{u_1} (p_1+(-1)^{v_1}p_2)+...$ \\
$+(-1)^{u_s}(p_{2s-1}+(-1)^{v_s}p_{2s}), \zeta
_{2s+1}+(-1)^{u_1}(\zeta _1+(-1)^{v_1}\zeta _2)+...+
(-1)^{u_s}(\zeta _{2s-1}+(-1)^{v_s}\zeta _{2s}))$ \\
$- \sum_{s=0}^{2^{r-1}-2}(i_{2s+1}2^{-2s-1})
\sum_{w,v_1,...,v_s,u_1,...,u_s\in \{ 1,2 \} } (-1)^{v_1+...+v_s+s}
S_n(p_0,p_{2s+1}+(-1)^wp_{2s+2}+$
\\ $(-1)^{u_1}(p_1+(-1)^{v_1}p_2+...+
(-1)^{u_s}(p_{2s-1}+(-1)^{v_s}p_{2s}),\zeta _{2s+1}+(-1)^w\zeta
_{2s+2}+(-1)^{u_1}(\zeta _1+(-1)^{v_1}\zeta _2+...+ (-1)^{u_s}(\zeta
_{2s-1}+(-1)^{v_s}\zeta _{2s}))$ \\
$- (i_{2^r-1}2^{-2^r+2})\sum_{v_1,v_2,...,v_{2^{r-1}-1},
u_1,u_2,..., u_{2^{r-1}-1} \in \{ 1,2 \} }
(-1)^{v_1+v_2+v_{2^{r-1}-1}+1}S_n(p_0,p_{2^r-1}+(-1)^{u_1}(p_1+$ \\
$(-1)^{v_1}p_2)+(-1)^{u_2}(p_3+(-1)^{v_2}p_4)+
...+(-1)^{u_{2^{r-1}-1}}(p_{2^r-3}+(-1)^{v_{2^{r-1}-1}} p_{2^r-2}),
\zeta _{2^r-1}+ (-1)^{u_1}(\zeta _1+(-1)^{v_1}\zeta
_2)+(-1)^{u_2}(\zeta _3+ (-1)^{v_2}\zeta
_4)+...+(-1)^{u_{2^{r-1}-1}}
(\zeta _{2^r-3}+(-1)^{v_{2^{r-1}-1}}\zeta _{2^r-2})) \} $, \\
where for the convenience of the notation in the Formula for $s=0$
the sum is not accomplished by $v_1,...,v_s, u_1,...,u_s$ and the
value  $(-1)^{v_1+...+v_s+s}=1$ is taken.
\par {\bf 34. Theorem.} {\it Let $f(t)$ be a function-original
with values in ${\cal A}_r$ with $2\le r<\infty $, $u=pt$, $g(t) :=
\int_0^t f(x)dx$, then
\par ${\cal F}(g(t)Ch_{[0,\infty )}(t),pt;p;0)p =
{\cal F}(f(t)Ch_{[0,\infty )}(t),u;p;0)$ \\
in the domain $Re (p)>\max (s_0,0)$, in particular,
\par ${\cal F}(g(t)Ch_{[0,\infty )}(t),pt;p;0) =
{\cal F}(f(t)Ch_{[0,\infty )}(t),u;p;0)p^{-1}$ over the algebra
${\bf K}=\bf H$ or ${\bf K}=\bf O$.}
\par {\bf Proof.} In view of Theorem 28
\par ${\cal F} (g'(t)Ch_{[0,\infty )}(t),pt;p;0) = {\cal F}
(g(t)Ch_{[0,\infty )}(t),pt;p;0)p$ \\  in the same domain $Re
(p)>\max (s_0,0)$, since $g(0)=0$ and $g(t)$ also satisfies
conditions of Definition 1, moreover, $s_0(g)<\max (s_0(f),0)+b$ for
each $b>0$, where $s_0\in \bf R$. From the alternativity of the
algebra $\bf K$ it follows, that
\par $ {\cal F} (f(t)Ch_{[0,\infty )}(t),pt;p;0)p^{-1}
= {\cal F} (g(t)Ch_{[0,\infty )}(t),pt;p;0)$,\\
since $g'(t)=f(t)$ for each $t>0$ and $g(t)=0$ for each $t\le 0$.
\par {\bf 35. Theorem.} {\it If $F(p)$ is an image of a
function $f(t)$ for $u=pt$ in the half space $\{ p\in {\cal A}_r: Re
(p)>s_0 \} $ with $2\le r<\infty $, the integral $\int_p^{\infty }
F(z)dz$ converges and Condition 23(3) is satisfied, then
\par ${\cal F} (f(t)Ch_{[0,\infty )}(t)/t,pt;p;0) =
\int_p^{\infty }F(z)dz$.}
\par {\bf Proof.} Let a path of an integration belong to the half space
$Re (p)\ge a$ for some constant $a>s_0$, then $|\int_0^{\infty
}f(t)\exp (-pt)dt|\le C\int_0^{\infty }\exp (-(p_0-s_0)t)dt < \infty
$ converges, where $C=const>0$, $p_0\ge a$. It can be supposed
$t>0$, then conditions of Lemma 23 are satisfied, where
$\int_p^{\infty }F(z)dz=\lim_{0<\theta \to 0}\int_p^{\infty
}F(z)\exp (-\theta z)dz$, since the integral $\int_{a-S\infty
}^{a+S\infty }F(p)dp$ is absolutely converging and $\lim_{\theta \to
0}\exp (-\theta z)=1$ uniformly by $z$ on each compact subset in
${\cal A}_r$. Therefore, in the integral
\par $\int_p^{\infty
} F(z)dz = \int_p^{\infty } (\int_0^{\infty }f(t)\exp (-pt)dt)dz $ \\
it can be changed the order of the integration: \\
$\int_p^{\infty } F(z)dz = \int_0^{\infty } f(t)(\int_p^{\infty }
\exp (-zt)dz)dt = \int_0^{\infty }f(t)[- e^{-zt}/t]|_p^{\infty }dt =
\int_0^{\infty }f(t)t^{-1}\exp (-pt)dt$,
\\  since in view of Lemma 23 it can be taken the argument $z-z_0$
varying along the straight line parallel to $p-p_0$ with $z_0$
tending to $+\infty $, where $p_0:=Re (p)$.
\par {\bf 36. Examples.} 1. For $u=pt$ we have ${\cal F}((e^{bt}-e^{ct})
Ch_{[0,\infty )}(t),pt; p;0) = (p-b)^{-1} - (p-c)^{-1}$ for each $b,
c \in \bf R$. Then in view of Example 26.6 and Theorem 35 ${\cal
F}((e^{bt}-e^{ct})Ch_{[0,\infty )}(t)/t,pt; p;0) = \int_p^{\infty
}((z-b)^{-1} - (z-c)^{-1})dz = Ln [(p-b)^{-1}(p-c)]$ (see also
Corollaries 3.5 and 3.6 and Notes 3.7, 3.8 in \cite{ludfov}).
\par 2. In view of Example 26 and Theorem 35
${\cal F} (\sin (t)Ch_{[0,\infty )}(t)/t,pt; p;0) = \int_p^{\infty
}((1+z^2)^{-1}dz = (\pi /2)- arctg (p)=arcctg (p)$.
\par Then from Theorem 34 it follows, that
${\cal F} (si (t)Ch_{[0,\infty )}(t),pt; p;0)
 = p^{-1} arcctg (p)$, where $si (t) := \int_0^t t^{-1} \sin (t)dt$
and it denotes the integral sine.
\par 3. Calculate an image $F(p)$ of the function-original
$f(t)=Ch_{[0,\infty )}(t) (e^{bt} - e^{at})/t$, where $t\in \bf R$,
parameters $a$ and $b$ belong to the algebra ${\cal A}_r$, $u=pt$.
Then
\par $(i)$ $F(p) = \int_0^{\infty } t^{-1} [\exp (b_0t) (\cos
(|b'|t) +\sin (|b'|t) b'/|b'|) - \exp (a_0t) (\cos (|a'|t) + \sin
(|a'|t)a'/|a'|)]$ \\  $\exp (-p_0t) (\cos (|p'|t) -
\sin (|p'|t)p'/|p'|) dt$, \\
where $a = a_0+a'$, $a_0 := Re (a)$, $b=b_0+b'$, $b_0:=Re (b)$,
$p=p_0+p'$, $p_0:=Re (p)$. Let $a_0\le b_0<p_0$. The addition and
subtraction from the right part of Formula $(i)$ of the integral
$\int_0^{\infty }t^{-1}\exp (b_0t)(\cos (|a' |t)+\sin
(|a'|t)a'/|a'|)\exp (-p_0t) (\cos (|p'|t) - \sin (|p'|t)p'/|p'|)dt$,
the grouping of terms and the using of Formulas 30.3(i-iv) give:
\par $(ii)$  $F(p)= - \int_0^{\infty } t^{-1}\exp ((b_0-p_0)t)
[\sin (((|b'| + |a'|)2^{-1}+|p'|)t) + \sin (((|b'|
+|a'|)2^{-1}-|p'|)t) + (\cos ((|b'|+|a'|)2^{-1}+|p'|)t) - \cos ((
|b'|+|a'|)2^{-1}-|p'|)t))p']\sin ((|b'|-|a'|)2^{-1}t)]dt +
\int_0^{\infty }t^{-1}\exp ((b_0-p_0)t) (\sin (|b'|t)\cos
(|p'|t)b'/|b'| -\sin (|b'|t)\sin (|p'|t)b'p'|b'|^{-1}|p'|^{-1}$ \\
$- \sin (|a'|t)\cos (|p'|t)a'/|a'| + \sin (|a'|t)\sin
(|p'|t)a'p'|a'|^{-1}|p'|^{-1})dt + 2^{-1}\int_0^{\infty }t^{-1}\exp
((b_0-p_0)t)(1-\exp ((a_0-b_0)t)(\cos ((|a'|+|p'|)2^{-1}t) + \cos
((|a'|-|p'|)2^{-1}t))dt$ $+\int_0^{\infty }t^{-1}[\exp ((b_0-p_0)t)
-\exp ((a_0- p_0)t] [ -\cos (|a'|t)\sin (|p'|t)p'/|p'| +\cos
(|p'|t)\sin (|a'|t)a'/|a'| - \sin (|p'|t) \sin (|a'|t)
a'p'|a'|^{-1}|p'|^{-1})dt$. \\
Now use Formulas for the appearing integrals with the parameters $v,
w, s\in \bf C$:
\par $(iii)$  $J_1(v,s) := \int_0^{\infty }\exp (-vx)x^{-1}\sin
(sx)dx= \arctan (s/v)$ for $Re (v)>|Im (s)|$;
\par $(iv)$  $J_2(v,w,s):=\int_0^{\infty }x^{-1}\exp (-vx)\sin (wx)
\sin (sx)dx= 4^{-1} ln [(v^2+(w+s)^2)(v^2+(w-s)^2)^{-1}]$, \\
for $Re (p)>|Im (w)|+|Im (s)|$;
\par $(v)$  $J_3(v,w,s):=\int_0^{\infty }x^{-1}\exp (-v x)\sin (w x)
\cos (s x)dx=2^{-1} \arctan (2wv(v^2+s^2-w^2)^{-1}) +
\{ {0\choose {\pi /2}} \} $, \\
where $Re (v)> |Im (w)|+|Im (s)|,$ $\{ {{v^2+s^2\ge w^2} \choose
{v^2+s^2<w^2}} \} $ and depending of the satisfaction of the upper
or the lower inequality it is taken $0$ or $\pi /2$ in $\{ {0\choose
{\pi /2}} \} $;
\par $(vi)$ $J_4(v,w,s):= \int_0^{\infty }\exp (-vx)(1-\exp
(-w x))\cos (s x)dx= (2w)^{-1} [B(2,(v-{\bf i}s)w^{-1})+ B(2,(v+{\bf
i}s)w^{-1})]$,\\  where $Re (w)>0$, $Re (v)>|Im (s)|$ (see Formulas
11 on page 447, 8, 10 on page 449 and 8 on page 450 in \cite{prud}),
$B(v,w) := \int_0^1\tau ^{v-1}(1-\tau )^{w-1}d\tau =\Gamma (v)\Gamma
(w)/\Gamma (v+w)$ and $\Gamma (v):=\int_0^{\infty }e^{-t}t^{v-1}dt$
при $Re (v)>0$, $Re (w)>0$. Therefore from Formula $(ii-vi)$ it
follows, that
\par $(vii)$ $F(p)=-J_2(p_0-b_0,(|b'|-|a'|)2^{-1},(|a'|+|b'|)
2^{-1}+|p'|) -
J_2(p_0-b_0,(|b'|-|a'|)2^{-1}, (|a'|+|b'|)2^{-1}-|p'|)-
J_3(p_0-b_0,(|b'|-|a'|)2^{-1}, (|b'|+|a'|)2^{-1}+|p'|)p' +
J_3(p_0-b_0,(|b'|-|a'|)2^{-1},(|b'|+|a'|)2^{-1}-|p'|)p' +
J_3(p_0-b_0,|b'|,|p'|)b'/|b'| - J_3(p_0-b_0,|a'|,|p'|) a'/|a'|$ \\
$-J_2(p_0-b_0,|b'|,|p'|)b'p'|b'|^{-1}|p'|^{-1} +
J_2(p_0-b_0,|a'|,|p'|)a'p'|a'|^{-1}|p'|^{-1} + 2^{-1}J_4(p_0-b_0,
b_0-a_0,(|a'|+|p'|)/2)+ 2^{-1}J_4(p_0-b_0, b_0-a_0,(|a'|-|p'|)/2)
-J_3(p_0-b_0,|p'|,|a'|)p'/|p'| + J_3(p_0-b_0,|a'|,|p'|)a'/|a'| -
J_2(p_0-b_0,|a'|,|p'|)a'p'|a'|^{-1}|p'|^{-1} +
J_3(p_0-a_0,|p'|,|a'|)p'/|p'| - J_3(p_0-a_0,|a'|,|p'|)a'/|a'| +
J_2(p_0-a_0,|a'|,|p'|)a'p'|a'|^{-1}|p'|^{-1}$, \\
where ${\bf i}=(-1)^{1/2}=i_1$. Moreover, \par ${\cal F}((e^{bt} -
e^{at})/t,pt;p;0)=\int_p^{\infty }({\cal F}(e^{bt},pt;z;0)- {\cal
F}(e^{at},pt;z;0))dz$, \\
where $z$ tends to the infinity in the
angle $|Arg (z)|<\pi /2-\delta $ for some $0<\delta <\pi /2$.
\par {\bf 37. Theorem.} {\it If $f(t)$ is a function-original, then
\par $(1)$  ${\cal F}((fCh_{[0,\infty )}(t-\tau ),u;p;\zeta )=
{\cal F}((f Ch_{[0,\infty )})(t),u;p; \zeta +p\tau )$ for
$u(p,t;\zeta )= p_0t + \zeta _0 + M(p,t;\zeta )$ or $u(p,t;\zeta )=
pt + \zeta $ over ${\cal A}_r$ with $2\le r<\infty $,
\par $(2)$ ${\cal F}((fCh_{[0,\infty )})(t-\tau ),pt;p;0)={\cal F}((f
Ch_{[0,\infty )})(t),pt;p;0)e^{-p\tau }$ over ${\bf K}=\bf H$ or
${\bf K}=\bf O$ with $f({\bf R})\subset \bf K$, or over ${\cal A}_r$
with $4\le r$ for $f({\bf R})\subset \bf R$, in the half space $Re
(p)>s_0$.}
\par {\bf Proof.} For $p$  in the half space $Re (p)>s_0$
the equality is satisfied: \par ${\cal F}((fCh_{[0,\infty )})
(t-\tau ),u;p;\zeta ) = \int_{\tau }^{\infty }f(t-\tau )
e^{-u(p,t;\zeta )}dt = \int_0^{\infty }f(t_1)e^{-u(p,t_1;\zeta
+p\tau )}dt_1$
\par   $={\cal F}((fCh_{[0,\infty )})(t),u;p;\zeta +p\tau )$, in
view of Formula 3(1-5), since $p_jt+\zeta _j = p_jt_1+(\zeta
_j+p_j\tau )$ for each $j=0,1,...,2^r-1$, where $t=t_1+\tau $, in
particular, for $u=pt$, over ${\bf K}=\bf H$ or ${\bf K}=\bf O$ with
$f({\bf R})\subset \bf K$, or over ${\cal A}_r$ with $4\le r$ for
$f({\bf R})\subset \bf R$:
\par  ${\cal F}((fCh_{[0,\infty )})(t-\tau ),pt;p;0)=\int_{\tau
}^{\infty }f(t-\tau ) e^{-pt}dt = \int_0^{\infty
}f(t_1)e^{-p(t_1+\tau )}dt_1$ \\  $ = {\cal F}((fCh_{[0,\infty
)})(t), pt;p;0)e^{-p\tau }$, since $t, \tau \in \bf R$, the algebra
$\bf K$ is alternative and the center of the algebra ${\cal A}_r$ is
$\bf R$ .
\par {\bf 38. Note.} In view of the definition of the transformation
${\cal F}$ and $u(p,t;\zeta )$ and Theorem 37 there can be
interpreted $\zeta _1 i_1+...+\zeta _{2^r-1}i_{2^r-1}$ as the
initial phase of retarding.
\par {\bf 39. Theorem.} {\it If $f(t)$ is a function-original
with values in ${\cal A}_r$ for $2\le r<\infty $, $b\in \bf R$, then
${\cal F}(e^{bt}f(t)Ch_{[0,\infty )}(t),pt;p;\zeta )= {\cal
F}(f(t)Ch_{[0,\infty )}(t),pt;p-b;\zeta )$ for each $Re (p)>s_0+b$.}
\par {\bf Proof.} If $Re (p)>s_0+b$, then the integral
\par ${\cal F}(e^{bt}f(t)Ch_{[0,\infty )}(t),pt;p;\zeta )= \int_0^{\infty
}f(t)e^{bt}\exp (-pt-\zeta )dt$ \\   $= \int_0^{\infty } f(t)\exp
(-(p-b)t -\zeta )dt = {\cal F}(f(t)Ch_{[0,\infty )}(t),pt;p-b;\zeta
)$ converges.
\par {\bf 40. Examples.} 1. ${\cal F}(e^{-bt}\sin (at)Ch_{[0,\infty )}(t),
pt;p;0)=
a[(p+b)^2+a^2]^{-1}$, where $a, b\in \bf R$, $Re (p)>b$, $p\in {\cal
A}_r$.
\par 2. ${\cal F}(e^{-bt}\cos (at)Ch_{[0,\infty )}(t),pt;p;0)=
(p+b)[(p+b)^2+a^2]^{-1}$, where $a, b\in \bf R$, $Re (p)>b$.
\par 3. ${\cal F}(e^{-bt}t^nCh_{[0,\infty )}(t),pt;p;0)=
n!(p+b)^{-n-1}$, where $b\in \bf R$, $Re (p)>b$ (see examples 26.2,
26.3 and 30.1).
\par {\bf 41. Theorem.} {\it If functions $f(t)$ and $g(t)$ are
originals, where $g$ is with real values for each $t$, then
\par ${\cal F}(\int_0^tf(\tau )g(t-\tau )d\tau , pt;p;0) = {\cal
F}(f(t)Ch_{[0,\infty )}(t),pt;p;0) {\cal F}(g(t)
Ch_{[0,\infty )}(t),pt;p;0)$ \\
for each $p\in {\cal A}_r$ с $Re (p)>s_0$, where $s_0 = \max
(s_0(f), s_0(g))$ and $2\le r<\infty $.}
\par {\bf Proof.} The convolution of functions
$q(t) := \int_0^tf(\tau )g(t-\tau )d\tau $ satisfies Conditions 1
and 2 of Definition 1. Moreover, $|\int_0^tf(\tau )g(t-\tau )d\tau
|<C |\int_0^t \exp (s_0\tau )\exp (s_0(t-\tau ))d\tau |= Ct\exp
(s_0t) <C_1\exp ((s_0+\epsilon )t)$, where $C, C_1$ are positive
constants, $\epsilon >0$ can be chosen an arbitrary small number for
the corresponding constant $C_1$. Hence, Property 3 of Definition 1
is satisfied, that is, the function $q(t)$ is original. Write the
functions $f$ and $g$ in the form: $f = \sum_v f_vi_v$, $g = \sum_v
g_vi_v$, where $f_v$ and $g_v$ are functions with values in $\bf R$,
$\{ i_v: v=0,1,...,2^r-1 \} $ are generators of the algebra ${\cal
A}_r$. Then in view of the Fubini Theorem
\par ${\cal F}(\int_0^tf(\tau )g(t-\tau )d\tau , pt;p;0) =
\sum_{w,v=0}^{2^r-1} (i_wi_v)\int_0^{\infty }
(\int_0^tf_w(\tau )g_v(t-\tau )d\tau ) \exp (-pt)dt$ \\
$=\sum_{w,v=0}^{2^r-1} (i_wi_v)(\int_0^{\infty }f_w(\tau )\exp
(-p\tau
)d\tau ) (\int_0^{\infty }g_v(t_1)\exp (-pt_1)dt_1)$, \\
for each $Re (p)>s_0$, where $t_1:=t-\tau $. Since  $g({\bf
R})\subset \bf R$, then $g=g_0$, that gives the Formula of this
Theorem.
\par {\bf 42. Theorem.} {\it  Let a function $f(t)$ be a real valued
original, \\ $F(p) = {\cal F}(f(t)Ch_{[0,\infty )}(t);u;p;0)$, let
also $G(p)$ and $q(p)$ be analytic functions such that
\par ${\cal F}(g(t,\tau )Ch_{[0,\infty )}(t);u;p;0) = G(p)
\exp (-u(q(p),\tau ;0))$ \\  for $u=pt$ or $u=p_0t+ M(p,t;0)$, then
\par ${\cal F}(\int_0^{\infty }g(t,\tau )f(\tau )d\tau ;u;p;0) =
G(p)F(q(p))$ \\  for each $p\in {\cal A}_r$ with $Re (p)>s_0$ and
$Re (q(p))>s_0$, where $s_0=\max (s_0(f),s_0(g))$, $2\le r<\infty
$.}
\par {\bf Proof.} Let $Re (p)>s_0$ and $Re (q(p))>s_0$,
where $s_0=\max (s_0(f),s_0(g))$, then in view of the Fubini Theorem
and the theorem's conditions the order of integration's change gives
the equalities:
\par $\int_0^{\infty } (\int_0^{\infty } g(t,\tau ) f(\tau )d\tau
)\exp (-u(p,t;0))dt = \int_0^{\infty }
(\int_0^{\infty } g(t,\tau ) \exp (-u(p,t;0))dt )f(\tau )d\tau $ \\
$= \int_0^{\infty }G(p)\exp (- u(q(p),\tau ;0))f(\tau )d\tau
=G(p)\int_0^{\infty }f(\tau )\exp (-u(q(p),\tau ;0)) d\tau =
G(p)F(q(p))$, \\
since $t, \tau \in \bf R$ and the center of the algebra ${\cal A}_r$
is $\bf R$.
\par {\bf 43. Examples.} Consider $G(p)=p^{-1/2}$ and $q(p)=p^{1/2}$,
where by the inversion's formula of Theorem 19 we have $g(t)=(2\pi
i_1)^{-1}\int_{a-i_1\infty }^{a+i_1\infty }p^{-1/2}\exp (-\tau
p^{1/2}+pt)dp$, \\
since $p^{s}$ and $p$ commute for $s\in \bf R$ and $p\in {\cal A}_r$
in view of the polar form $p = |p|\exp (\theta M)$ of the
Cayley-Dickson numbers, in particular, quaternions and octonions,
where $\theta \in \bf R$, $M\in {\cal A}_r$, $Re (M)=0$, $|M|=1$ and
the exponential function's Formula $\exp (\theta M) = \cos (\theta )
+ \sin (\theta )M$ \cite{ludfov,ludoyst}. \par We consider at first
values  $t>0$. Let $0<v<R<\infty $ and consider a loop (closed
curve) $\gamma $ composed of the segment $\{ z\in {\cal A}_r:
z=a+i_1 c, c\in [-b,b] \} $, where $a^2+b^2=R^2$, arcs $S_1(R)$ and
$S_2(R)$ of the circle $S(R) = \{ z\in {\cal A}_r: z=z_0+i_1z_1,
|z|=R, z_0, z_1\in {\bf R} \} $, where the arc $S_1(R)$ joins the
point $(a+i_1b)$ with $(-R+i_1\epsilon )$ while going clockwise,
$S_2(R)$ joins the point $(-R-i_1\epsilon )$ with $(a-i_1b)$,
$\epsilon >0$ is a small number. The points $(-v+i_1\epsilon )$ and
$(-R+i_1\epsilon )$ are joined by the segment $[-v+i_1\epsilon ,
-R+i_1\epsilon ]$ of the straight line, also $(-R-i_1\epsilon )$ and
$(-v-i_1\epsilon )$ are joined by the straight line's segment
$[-R-i_1\epsilon , -v-i_1\epsilon ]$. The rest of $\gamma $ is
composed of the arc $S_3(v)$ of the circle $S(v)$ from the point
$-v-i_1\epsilon $ to $-v+i_1\epsilon $ while going counterclockwise.
Inside this contour we put $-\pi <\theta <\pi $, where $M=i_1$. Then
in view of Theorem  2.11 \cite{ludfov,ludoyst} $\int_{\gamma
}p^{-1/2}\exp (-\tau p^{1/2}+pt)dp=0$, consequently,
\par $\int_{a-i_1\infty }^{a+i_1\infty }p^{-1/2}\exp (-\tau
p^{1/2}+pt)dp= $ \\  $\int_{S_2(R)\cup [-R-i_1\epsilon ,
-v-i_1\epsilon ]\cup S_3(v)\cup [-v+i_1\epsilon , -R+i_1\epsilon
]\cup S_1(R)} p^{-1/2}\exp (-\tau p^{1/2}+pt)dp$. \\ In view of
Lemma 12.1 $\lim_{R\to \infty }\int_{S_2(R)\cup S_1(R)} p^{-1/2}\exp
(-\tau p^{1/2}+pt)dp=0$, consequently, \par $g(t;\tau )=(2\pi
i_1)^{-1}\int_{ [-R-i_1\epsilon , -v-i_1\epsilon ]\cup S_3(v)\cup
[-v+i_1\epsilon , -R+i_1\epsilon ]} p^{-1/2}\exp (-\tau
p^{1/2}+pt)dp$. \\ Along the segment $[-R-i_1\epsilon ,
-v-i_1\epsilon ]$ we have $p=x\exp (-i_1\pi )-i_1\epsilon $ and
$(p+i_1\epsilon )^{1/2}=-i_1x^{1/2}$ and along the segment
$[-v+i_1\epsilon , -R+i_1\epsilon ]$ we get $p=x\exp (i_1\pi
)+i_1\epsilon $ and $(p-i_1\epsilon )^{1/2}=i_1x^{1/2}$, where
$x>0$, $x^{1/2}$ is the positive branch of the square root
(arithmetical value), consequently, \par $\lim_{\epsilon \to 0}
\int_{ [-R-i_1\epsilon , -v-i_1\epsilon ]} p^{-1/2}\exp (-\tau
p^{1/2}+pt)dp=\int_R^v{\tilde i}_1 x^{-1/2}\exp (i_1\tau
x^{1/2}-xt)dx$ and \par $\lim_{\epsilon \to 0} \int_{
[-v+i_1\epsilon , -R+i_1\epsilon ]} p^{-1/2}\exp (-\tau
p^{1/2}+pt)dp=- \int_v^R{\tilde i}_1x^{-1/2}\exp (-i_1\tau
x^{1/2}-xt)dx$, also \par  $\lim_{v\to 0}\int_{S_3(v)} p^{-1/2}\exp
(-\tau p^{1/2}+pt)dp=0$, since $|\int_{S_3(v)} p^{-1/2}\exp (-\tau
p^{1/2}+pt)dp|\le C 2\pi v^{1/2}$, where $C=const>0$. Make the
change of the variable $x=y^2$ and use the known Poisson integral
$\int_0^{\infty }\exp (-ax^2)\cos (bx)dx=(\pi /a)^{1/2}\exp
(-b^2/(4a))/2$, where $a>0$, $b\in \bf R$, then $g(t;\tau )=\pi
^{-1}\int_0^{\infty }x^{1/2}\cos (\tau x^{1/2})\exp (-xt)dx= (2/\pi
)\int_0^{\infty }\exp (-ty^2)\cos (\tau y)dy=(\pi t)^{-1/2} \exp
(-\tau ^2/(4t))$. Further we get, that $\lim_{R\to \infty
}\int_{S_4(R)} p^{-1/2}\exp (-\tau p^{1/2}+pt)dp=0$ for $t<0$, where
$S_4(R)$ is the arc of the circle $S(R)$ joining points $(a-i_1b)$
and $(a+i_1b)$ while going counterclockwise with $Re (p)>0$,
consequently, the analogous reasoning to the given above shows, that
$g(t;\tau )=0$ for $t<0$. Then
\par $(1)$ $g(t;\tau )=(\pi t)^{-1/2} \exp (-\tau ^2/(4t))$ for
$t>0$. \\  Thus, in view of Theorem 42:
\par $(1)$ ${\cal F}((\pi t)^{-1/2}\int_0^{\infty }\exp (-\tau ^2/(4t))
f(\tau )d\tau ;u;p;0) = p^{-1/2}F(p^{1/2})$, \\
where $f(t)$ is the original of the function $F(p)$. \par In
particular, take $F(p)=p^{-1}\exp (-ap)$, where $0<a\in \bf R$. In
accordance with Formula 26(2) $f(t)=Ch_{[0,\infty )}(t-a)$. Then in
view of Formula (1) we get:
\par $(2)$  $p^{-1}\exp (-ap^{1/2})={\cal F}((\pi t)^{-1/2}
\int_a^{\infty }\exp (-\tau ^2/(4t))d\tau ;u;p;0)$, \\
where $\int_a^{\infty }\exp (-\tau ^2/(4t))d\tau = 2(\pi )^{-1/2}
\int_{t^{-1/2}a/2}^{\infty }\exp (-x^2)dx$, $x := t^{-1/2}\tau /2$.
Denoting in the standard manner $erf (y):=\int_0^y\exp (-x^2)dx$,
$Erf (y):=1-erf (y)$ we write Formula (2) in the form:
\par $(3)$  ${\cal F}(Erf(t^{-1/2}a/2);u;p;0)=p^{-1}\exp
(-ap^{1/2})$ for $u=pt$ and $Re (p)>0$, $p\in {\cal A}_r$, $2\le
r\in \bf N$.
\par {\bf 44. Theorem.} {\it Let a function $F(p)$ be holomorphic in
the half space $Re (p)>s_0$, where $p\in {\cal A}_r$, $2\le r\in \bf
N$, $s_0\in \bf R$. Moreover, $F(p)$ is regular at the infinite
point and we have in its neighborhood $\{ p\in {\cal A}_r: |p|\ge R
\} $ the Loran decomposition:
\par $(1)$  $F(p) = \sum_{l=1}^{\infty }c_lg_l(p)$, \\
where $c_l\in {\cal A}_r$, $g_l(p):=p^{-l}$ for $u=pt$, $g_l(p)$
with $l=n-1$ given by Formulas (12,13) of Example 33 for
$u=p_0t+\zeta _0+M(p,t;\zeta )$, $0<R<\infty $, then the original
$f(t)$ is the function:
\par $(2)$  $f(t)=Ch_{[0,\infty )}(t)
\sum_{l=1}^{\infty }c_lt^{l-1}/(l-1)!$.}
\par {\bf Proof.} Put $q=p^{-1}$ and denote
$G(q):=F(1/q)$. Therefore, the function $G(q)$ is analytic in the
ball $ \{ q: |q| \le R^{-1} \} $. For $u=pt$ each function is equal
to $g_l(1/q)=q^l$, for $u=p_0t+\zeta _0+M(p,t;\zeta )$ the function
$g_l(1/q)$ is holomorphic by $q$, since in Example 33 and the
following examples there is used the expression of $p_j$ through
$p$. With the help of generators of the algebra ${\cal A}_r$ (see
Formulas 3(6)), that gives the local analyticity of $F(p;\zeta )$ by
$p$ and $\zeta $. In view of Formulas 33(4,10,11) asymptotically
$|q_l(1/q)|$ behaves itself like $|q_l(1/q)|\le C'|q|^l$ for $|q|\to
\infty $ in the plane ${\bf R}\oplus i_1{\bf R}$, where $C'=const
>0$ is independent from $l=n-1$. Then Inequality (4) from Theorem 2.7
and Theorem 3.21 \cite{ludfov,ludoyst} give $|c_l|< C_0C_1\exp
(C_2R^n) R^l$ for each $l\in \bf N$, where $n= 2^r+2$, $2\le r\in
\bf N$, in particular, $r=2$ for ${\bf K}=\bf H$, $r=3$ for ${\bf
K}=\bf O$, $C_1$ and $C_2$ are positive constants independent from
functions, $C_0=\max_{|q|\le 1/R} |G(q)|$. Therefore,
\par $|f(t)|\le \sum_{l=1}^{\infty } |c_l| t^{l-1}/(l-1)! \le
C_0C_1\exp (C_2R^n)\sum_{l=0}^{\infty }(R|t|)^l/l!=C_0C_1\exp
(C_2R^n+R|t|)$, consequently,  \par  $|f(t)|\le C\exp (Rt)$ for each
$t\ge 0$, where $C=C_0C_1\exp (C_2R^n)$. Thus, the function
$f(t)Ch_{[0,\infty )}(t)$ is original, since $t\in \bf R$, the
center of the algebra ${\cal A}_r$ is  $\bf R$. In view of uniform
convergence of Series (2) it can be integrated termwise with the
multiplier $\exp (u(p,t;\zeta ))$ by $t$ from $0$ to $\infty $ when
$Re (p)>R$. It is necessary to note, that $|{\cal
F}(t^{l-1}/(l-1)!,u;p;\zeta )| \le [(l-1)!]^{-1} (\int_0^{\infty
}t^{l-1}\exp (-p_0t)dt)\le p_0^{-l}$ for $u=p_0t+\zeta
_0+M(p,t;\zeta )$. Since ${\cal F}(t^n;u;p,\zeta )=g_l(p)$, where
$\zeta =0$ for $u=pt$, $\zeta \in {\cal A}_r$ can be different from
zero for $u=p_0t+\zeta _0+M(p,t;\zeta )$, then from this and
Examples 26 and 33 it follows the statement of this Theorem.
\par {\bf 45. Definition.} A function $F: U\to {\cal A}_r$ from a
domain $U$ into the algebra ${\cal A}_r$ is called meromorphic, if
it is holomorphic on $U\setminus P$, where a set $P$ of point poles
of a function $F$ is countable, moreover, at each bounded subdomain
$V$ in $U$ a subset $V\cap P$ of poles is finite, in particular, it
can be $P=\emptyset $.
\par {\bf  46. Theorem.} {\it Let a function $F(p)$ satisfy conditions:
\par $(1)$ $F(p)$ is meromorphic in the half space
$W:= \{ p\in {\cal A}_r: Re (p)>s_0 \} $ and all its poles may be
only of finite orders;
\par $(2)$ moreover, for each $a>s_0$ there exist constants $C_a >0$ and
$\epsilon _a >0$ such that $|F(p)|\le C_a\exp (-\epsilon _a |p|)$
for each $p\in {\cal A}_r$ with $Re (p)\ge a$, where $s_0$ is fixed,
$2\le r<\infty $, $f({\bf R})\subset \bf R$ for $r\ge 4$;
\par $(3)$ the integral $\int_{a+i_1\theta , \theta \in \bf R} F(p)dp$
converges absolutely for each $a>s_0$, where $i_0, i_1,...,$ \\
$i_{2^r-1}\in {\cal A}_r$ are the standard generators of the algebra
${\cal A}_r$. Then an original for $F(p)$ is the function
\par $(4)$ $f(t)= (i_1)^{-1}
\sum_{p_k} res (p_k,F(p)\exp (u(p,t;0)).i_1$ \\
for  $r=2, 3$; and for $4\le r\in \bf N$, if $f({\bf R})\subset \bf
R$, where the sum of residues is taken over all singular points
$p_k$ of the function $F(p)$ in the order of non decreasing their
absolute values.}
\par {\bf  Proof.} In view of Theorem 22 over ${\bf K}=\bf H$ or
${\bf K}=\bf O$ with $f({\bf R})\subset \bf K$,  also for the
arbitrary algebra ${\cal A}_r$ with $f({\bf R})\subset \bf R$ for
$4\le r\in \bf N$, a function $F(p)$ is an image
\par $(5)$
$f(t)=(2\pi )^{-1}{\tilde i_1} \int_{a-i_1\infty }^{a+i_1\infty
}F(p)\exp (u(p,t;0))dp$. From Definition 34 it follows, that there
exists a sequence of radii $0<R_n<R_{n+1}$ for each $n\in \bf N$,
such that $\lim_{n\to \infty }R_n = \infty $, the sphere $S_n :=
S({\cal A}_r,0,R_n)$ does not contain any pole $p_k\in P$,
consequently, $\min_{p\in P, z\in S_n} |p-z| =: \delta _n >0$ in
view of the compactness of the sphere $S_n$ and the ball $B({\cal
A}_r,0,R_{n+1})$, where $S({\cal A}_r,y,R) := \{ z\in {\cal A}_r:
|z-y|=R \} $, $y\in {\cal A}_r$, $0<R<\infty $, $B({\cal A}_r,y,R)
:= \{ z\in {\cal A}_r: |z|\le R \} $. Consider the section of the
sphere $S_n$ by the plane ${\bf R}\oplus i_1{\bf R}$ and the part
$\gamma _n$ of the circle $\psi _n := S_n\cap ({\bf R}\oplus i_1{\bf
R})$ displayed to the right of the hyperplane $\{ p\in {\cal A}_r:
Re (p)=a \} $ over $\bf R$ in ${\cal A}_r$, that is, in the half
space $Re (p)>a$. Denote by $a-i_1b_n$ and $a+i_1b_n$ points of the
circle's intersection $\psi _n$ with the hyperplane $\{ p\in {\cal
A}_r: Re (p)=a \} $, also by $w_n$ the loop passed counterclockwise
and compose of $\gamma _n$ and the segment $[a-i_1b_n, a+i_1b_n]$.
In view of Lemma 12.1
\par $(6)$ $\lim_{n\to \infty }\int_{\gamma _n} F(p)\exp
(u(p,t))dp=0$, consequently,
\par $(7)$  $f(t)=(2\pi )^{-1}{\tilde i_1}
\lim_{n\to \infty } \int_{w_n}F(p)\exp (u(p,t;0))dp$. \\ Then from
Theorem 3.23 \cite{ludfov,ludoyst} and Formulas $(6,7)$ above
Formula $(4)$ of this theorem follows.
\par {\bf 47. Note.} Lemma 23 is also true under the condition:
\par $(2')$ there exists a sequence of hyper spheres
$S_n=S({\cal A}_r,0,R_n)$ over $\bf R$ embedded into ${\cal A}_r$
with the center at zero and radii $R_n$ with $R_n<R_{n+1}$ for each
$n\in \bf N$, $\lim_{n\to \infty }R_n = \infty $, such that
$\lim_{n\to \infty }
\sup_{p\in {S_n\cap W}} |{\hat F}(p)| =0$ \\
instead of Condition $(2)$ while others conditions of Lemma 23 are
satisfied.
\par Indeed, it is sufficient to prove it for some sequence of arcs
$\gamma _n$ of circles $\psi _n$ contained in the plane ${\bf
R}\oplus N\bf R$, where $N\in {\cal A}_r$, $Re (N)=0$, $|N|=1$. \par
If $F(p)$ is holomorphic in $W$, then due to Theorem 2.11
\cite{ludfov,ludoyst} $\int_{\gamma _n} F(p)\exp (-u(p,t;\zeta ))dp$
does not depend on a form of a curve and it is defined by an initial
and final points. If $V(\gamma _n)\le C_V R_n$ for each $n$, then it
is sufficient to prove the statement of the lemma for each
subsequence $R_{n(k)}$ with $R_{n(k+1)}\ge R_{n(k)}+1$ for each
$k\in \bf N$, that also can be seen from the estimate of integrals
below. Denote such subsequence by $R_n$. Each rectifiable curve can
be approximated by a converging sequence of rectifiable polygonal
lines, composed of arcs of circles. If a curve is displayed on a
sphere, then these circles can be taken on this sphere with the same
center. Since Condition $(2')$ is satisfied uniformly relative to
directrix $N$ and there can be accomplished a diffeomorphism $g$ in
${\cal A}_r$, such that $g(W)=W$, $g(S_n)=S_n$ for each $n\in \bf N$
and an image of a $C^1$ curve from $W$ is an arc of a circle, since
$0<R_n+1<R_{n+1}$ for each $n\in \bf N$ and $\lim_{n\to \infty
}R_n=\infty $. The functional $(F,\gamma )\mapsto \int_{\gamma
}F(p)dp$ from $C^0(V,{\cal A}_r)\times \Gamma $ into ${\cal A}_r$ is
continuous, where $V$ is a compact domain in ${\cal A}_r$, $\Gamma $
is a family of rectifiable curves in $V$ with the metric $\rho
(v,w):= \max (\sup_{z\in v}\inf_{\zeta \in w}|z-\zeta |, \sup_{z\in
w}\inf_{\zeta \in v}|z-\zeta |)$ (see Theorem 2.7 in
\cite{ludfov,ludoyst}). In addition a rectifiable curve is an
uniform limit of $C^1$ curves, since each rectifiable curve is
continuous, where the space of all $C^1$ functions is dense in the
space of all continuous functions $C^0$ relative to the compact-open
topology, each $C^1$ curve is rectifiable, that is, consider $\gamma
_n=\psi _n\cap \{ p\in {\cal A}_r: Re (p)>a \} \cap W$. Using an
automorphism of the algebra ${\cal A}_r$ it is sufficient to prove
the equality
\par $\lim_{n\to \infty }\int_{\gamma _n} F(p)\exp
(-u(p,t))dp=0$ \\
for each $t>0$ for $u=pt+\zeta $ in accordance with Lemmas 8.1, 8.2
and the proof of Theorem 19. Moreover, $\exp (tNz)dz=(tN)^{-1}d\exp
(tNz)$ in the plane ${\bf R}\oplus N\bf R$. Make the substitution of
the variable $z=(p+\zeta /t)N$, then \par $\lim_{n\to \infty
}\int_{\gamma _n} F(p)\exp (- u(p,t;\zeta ))dp = \lim_{n\to \infty
}[\int_{\eta _n}F({\tilde N}z)\exp (tNz)dz]{\tilde N} $, \\
where $\eta _n=\gamma _nN$, а $z=x+Ny$, where $x, y\in \bf R$, $z\in
{\bf R}\oplus N\bf R$, $b_n:= \sup_{p\in {S_n\cap W}} |{\hat
F}_n(p)|$, $\phi _n := \arcsin (a/R_n)$. Then \par $|\int_{\eta _n;
-a\le y\le 0}F({\tilde N}z)\exp (tNz)dz|\le b_n\exp (at)\phi _nR_n$.
\\ Since $\sin (\phi )\ge 2\phi /\pi $ for each $0\le \phi \le
\pi /2$, then $|\exp (tNz)| = \exp (-tR_n\sin (\phi )) \le \exp
(-2tR_n\phi /\pi )$, consequently, \\ $|\int_{\eta _n; 0\le y\le
R_n} F({\tilde N}z)\exp (tNz)dz| \le 2b_nR_n\int_0^{\pi /2} \exp
(-2tR_n\phi /\pi )d\phi =2b_n\pi (1-\exp (-tR_n))/(2t)$. \\
Thus, $\lim_{n\to \infty }\int_{\eta _n} F({\tilde N}z)\exp
(-u({\tilde N}z,t))dz=0$ for each $t>0$. From this it follows
Equality 12.1$(4)$.
\par  Then Theorem 22 is accomplished under imposing Condition $(2')$
instead of Condition $(i)$ while satisfaction all others conditions
of Theorem 22, hence in such variant under Condition $(2')$ instead
of $(2)$ Theorem 46 is true.
\par {\bf 48. Theorem.} {\it If a function $f(t)$ and its derivative
$f'(t)$ are originals, $F(p)$ is the image function of $f(t)$ for
$u=pt$ over ${\cal A}_r$ with $2\le r\in \bf N$, then
\par $(1)$  $\lim_{p\to \infty }F(p)p =f(0)$, where $p\to \infty $
inside the angle $|Arg (p)|<\pi /2 - \delta $ for some $0<\delta
<\pi /2$ and $f(0)=\lim_{t\to +0}f(t)$; if in addition there exists
the limit $\lim_{t\to \infty }f(t)=f(\infty )$, then
\par $(2)$  $\lim_{p\to 0}F(p)p=f(\infty )$, where $p\to 0$ inside the
same angle.}
\par {\bf Proof.} If $z\in {\cal A}_r$, then it has the polar
decomposition $z=|z|\exp (M)$, where $M\in {\cal A}_r$, $Re (M)=0$,
$Arg (z):=M$ (see Corollary 3.6 in \cite{ludfov}). In view of
Theorem 28 there is satisfied the equality ${\cal
F}(f'(t)Ch_{[0,\infty )}(t),pt;p;0)=F(p)p-f(0)$. Then in accordance
with Note 8 $\lim_{p\to \infty , |Arg (p)|<\pi /2-\delta }
(F(p)p-f(0))=0$. If there exists $\lim_{t\to \infty }f(t)=f(\infty
)$, where $|f(\infty )|<\infty $, then $f$ is bounded on $\bf R$,
consequently, its rate of growth is $s_0\le 0$, therefore, the
function $F(p)$ is defined for each $Re (p)>0$. From the formula
$\int_0^{\infty }f'(t)\exp (-pt)dt=F(p)p-f(0)$ it follows, that for
$p=0$ the integral $\int_0^{\infty }f'(t)dt$ converges. Then in the
angle $|Arg (p)|<\pi /2 -\delta $ the integral $\int_0^{\infty
}f'(t)\exp (-pt)dt$ converges uniformly by $p$. Therefore, there can
be going to the limit for $p\to 0$ in the same angle, that is,
$\int_0^{\infty } f'(t)dt = \lim_{p\to 0} F(p)p - f(0)=f(\infty
)-f(0)$, consequently, $\lim_{p\to 0} F(p)p - f(0)=f(\infty ).$
\par {\bf 49. Theorem.} {\it  If a function $f(t)$ is original
together with its derivative, where $F(p;\zeta )$ is an image
function of $f(t)$ over the Cayley-Dickson algebra ${\cal A}_r$ with
$2\le r\in \bf N$ for $u=p_0t+ \zeta _0 + M(p,t;\zeta )$, then
\par  $\lim_{p\to
\infty } p_0F(p;\zeta )+p_1F(p;\zeta -i_1\pi /2)+...+p_{2^r-1}
F(p;\zeta -i_{2^r-1}\pi /2)=f(0)$,\\
where $f(0)=\lim_{t\to +0} f(t)$,  $p$ tends to the infinity inside
the angle $|Arg (p)|<\pi /2-\delta $ for some $0<\delta <\pi /2$. If
there also exists  $\lim_{t\to \infty }f(t)=f(\infty )$, then
$\lim_{p\to 0} \{ p_0F(p;\zeta )+p_1F(p;\zeta -i_1\pi /2)+...+
p_{2^r-1}F(p;\zeta -i_{2^r-1}\pi /2) \} =f(\infty )$, where $p\to 0$
inside the same angle.}
\par {\bf Proof.} In view of Theorem 31 the equality is satisfied:
\par ${\cal F}(f'(t)Ch_{[0,\infty )}(t),u;p;\zeta )=p_0F(p;\zeta
)+p_1F(p;\zeta -i_1\pi /2)+...+p_{2^r-1}F(p;\zeta -i_{2^r-1}\pi /2)-
f(0)$, for $u=p_0t+\zeta _0+M(p,t;\zeta )$, \\ where
$p=p_0+p_1i_1+...+p_{2^r-1}i_{2^r-1}$, $p_0,...,p_{2^r-1}\in {\cal
A}_r$, $ \{ i_0,...,i_{2^r-1} \} $ are the generators of the algebra
${\cal A}_r$. From Note 8 it follows, that $\lim_{p\to \infty , |Arg
(p)|<\pi /2-\delta } {\cal F}(f'Ch_{[0,\infty )},u;p;\zeta )=0$,
which gives the first statement of the  Theorem.
\par If there exists $\lim_{t\to \infty }f(t)=f(\infty )$, then
$f$ is bounded on $\bf R$ and $s_0\le 0$, where $F(p;\zeta )$ is
defined for each $Re (p)>0$. Therefore, there exists the limit
\par $\lim_{p\to 0, |Arg (p)|<\pi /2 -\delta }\int_0^{\infty } f'(t)\exp
(-p_t-\zeta _0 -M(p,t;\zeta ))=\int_0^{\infty } f'(t)dt=f(\infty
)-f(0)=\lim_{p\to 0, |Arg (p)|<\pi /2 -\delta } \{ p_0F(p;\zeta
)+p_1F(p;\zeta -i_1\pi /2)+...+ p_{2^r-1}F(p;\zeta -i_{2^r-1}\pi /2)
-f(0) \} $, \\
from which the second statement of this theorem follows.
\par {\bf 50. Example.} Consider images of fractional powers.
The Euler gamma-function is given by the integral $\Gamma (a+1) :=
\int_0^{\infty } t^ae^{-t}dt$ for each $Re (a)>-1$. In the polar
form a number $p\in {\cal A}_r$ has the form: $p=\rho \exp (S\alpha
)$, where $\rho =|p|$, $\alpha \in \bf R$, $S\in {\cal A}_r$, $Re
(S)=0$, $|S|=1$, $Arg (p)=S\alpha $ (see Corollary 3.6
\cite{ludfov}). For $Re (p)>0$ take $-\pi /2 <\alpha <\pi /2$.
Introduce new variable $q:=tp^{-1}$, then $\Gamma
(a+1)=p^{a+1}\int_L q^a \exp (-pq)dq$ (see Definition 4.1 and
Proposition 4.2 \cite{ludfov}), where the integral is taken along
the ray $L$ characterized by the condition $Arg (q)=-\alpha S$. For
points of the arc $\gamma _R := \{ q\in {\cal A}_r: |q|=R, -\alpha <
(Arg (q))S^*<0 \} $ put $q = R\exp (S\phi )$. Then $|\int_{\gamma
_R}q^a\exp (-pq)dq| \le R^a\int_{-a}^0 \exp (-\rho R\cos (\phi
+\alpha ))Rd\phi $. Since $0<\alpha +\phi <\alpha $, then $\cos
(\alpha +\phi )\ge \delta _0>0$, where $\delta _0=const
>0$. Therefore, there exists the limit
\par  $\lim_{R\to \infty }\int_{\gamma _R}q^a\exp (-pq)dq =0$, since
$0<\alpha +\phi <\alpha $ and $\cos (\alpha +\phi )>c_0$, where
$c_0=conts >0$. For each $S\in {\cal A}_r$ with $|S|=1$, $Re (S)=0$
between the ray $L$ and the real axis $i_0{\bf R}$ in ${\cal A}_r$
there are not singular points of the function under the integral.
Thus, the integral $\int_L q^a \exp (-pq)dq$ can be substituted on
$\int_0^{\infty } q^a \exp (-pq)dq$, consequently, $p^{-a-1}\Gamma
(a+1)=\int_0^{\infty }t^a\exp (-pt)dt$ over $\bf H$ or $\bf O$ in
view of the alternativity of the octonion algebra and the
associativity of the quaternion skew field; or $\Gamma
(a+1)=p^{a+1}\int_0^{\infty }t^a\exp (-pt)dt$ in the case of the
algebra ${\cal A}_r$ with $4\le r\in \bf N$, where the variable of
the integration is denoted anew by $t$. Thus,
\par $(i)$  ${\cal F}(t^aCh_{[0,\infty )}(t),pt;p;0)=
p^{-a-1}\Gamma (a+1)$ over $\bf H$ or $\bf O$; or
\par $(i')$  $p^{a+1}{\cal F}(t^aCh_{[0,\infty )}(t),pt;p;0)=
\Gamma (a+1)$ over ${\cal A}_r$ with $4\le r\in \bf N$. For $Re
(a)\ge 0$ the function $f(t)$ is original, while for $-1<Re (a)<0$
the function $f(t)=t^a$ is increasing unboundedly for $t$ tending to
zero does not satisfy conditions imposed on an original. But for the
latter interval of values of the parameter $a$ the integral
converges and Formula $(i)$ or $(i')$ respectively is satisfied. In
this situation it is possible to say, that $t^a$ is the exceptional
original, while $p^{-a-1}\Gamma (a+1)$ is the exceptional image over
$\bf H$ or $\bf O$; or $p^{a+1}t^a$ is the exceptional original for
the exceptional image $\Gamma (a+1)$ over the Cayley-Dickson algebra
${\cal A}_r$ with $4\le r\in \bf N$.
\par In the particular case $a=-1/2$ we get: $p^{-1/2} ={\cal F}((\pi
t)^{-1/2},pt;p;0)$ over $\bf H$ or $\bf O$; or $1 =p^{1/2}{\cal
F}((\pi t)^{-1/2},pt;p;0)$ over ${\cal A}_r$ with $4\le r\in \bf N$.
\par The case $u=p_0t+\zeta _0+M(p,t;\zeta )$ with
$\Gamma _M(a+1;\zeta ) := \int_0^{\infty }t^a\exp
(-M((1,0,...,0),t;\zeta )-t)dt$ for each $Re (a) > -1$ reduces to
that of considered above, since $M(g,t;\zeta )=M(0,0;\zeta )$ for
$g=Re (g)$. But in the general case, $\int_0^{\infty }t^a\exp
(-p_0t- M(p,t;0)) dt$ can be expressed through integrals (see
Formula 4 on page 446 in \cite{prud}):
\par $(1)$ $T_a(v):=\int_0^{\infty }t^a\exp (-vx)\cos (bx)dx=
\Gamma (a+1)(b^2+v^2)^{-(a+1)/2}\cos (c)$,
\par $(2)$ $S_a(v):=\int_0^{\infty }t^a\exp (-vx)\sin (bx)dx=
\Gamma (a+1)(b^2+v^2)^{-(a+1)/2}\sin (c)$, \\
where $c=(a+1) \arctan (b/p)$, $Re (a)>-1$, $Re (v)>|Im (b)|$, like
as it was done above in Example 33. That is, it is sufficient to use
Formulas of Example 33 with the substitution of integrals $(6,7)$ on
$T_a$ and $S_a$ respectively.
\par Let $S\in {\cal A}_r$, $Re (S)=0$, $|S|=1$, $\psi $
is a path in the plane ${\bf R}\oplus S\bf R$ embedded into ${\cal
A}_r$, such that $\psi $ goes along the two shores cut of the plane
${\bf R}\oplus S\bf R$ from $-\infty $ to $-\delta $ from the side
$z=a+Sb$ with $b<0$, then passes the circle of radius $\delta
>0$ with the center at zero clockwise and goes along the straight line
from $-\delta $ to $-\infty $ from the side $z=a+Sb$ with positive
$b>0$, where $2\le r\in \bf N$, $a, b\in \bf R$. Then denote
\par $(3)$ $Re (N_1S^*)\int_{\psi }\zeta ^a\exp
(u(\zeta ,1;0))d\zeta =:(2\pi S) (\Gamma _u(-a))^{-1}$,\\
where $u(p,t;\zeta ) = p_0t+\zeta _0+M(p,t;\zeta )$, $a, p, \zeta
\in {\cal A}_r$, $Re (N_1S^*)\ne 0$. In the particular case of
complex numbers and $u=pt$ Formula $(3)$ is known as the Hankel
integral representation of the gamma-function. For $u=pt$ and $p\in
{\cal A}_r$ the function $\Gamma _u(-a)$ coincides with the usual
gamma-function $\Gamma (-a)$ of the variable $a\in {\cal A}_r$ (see
Definition 4.1 and Propositions 4.2, 4.12 and Corollary 4.13 in
\cite{ludfov}).
\par {\bf 51. Theorem.} {\it Let a holomorphic function  $F(p)$
satisfy Conditions 23$(1-3)$ or 47$(2')$ instead of 23$(2)$ for
$p\to \infty $ in the domain $Re (p)<a$, $F(p)=\sum_{j=0}^{\infty }
p^{j\beta +\alpha }c_j$ for $Re (p)<a$, $F(p)$ has not singular
points in ${\cal A}_r$ may be besides zero, which is the branching
point of a finite order, while the integral $\int_{(-\infty ,
-\delta ]\cup S({\cal A }_r,0,\delta )} |c_j| |p^{j\beta +\alpha
}|dt $ converges for some $\delta >0$, where $2\le r\in \bf N$,
$\alpha \in \bf R$ and  $\beta $ is a positive rational number,
$c_j\in \bf R$ for each $j$ for $4\le r\in \bf N$, $c_j\in \bf K$
for ${\bf K}=\bf H$ or ${\bf K}=\bf O$. Then the original function
is \par $(i)$ $f(t)=Ch_{[0,\infty )}(t) t^{-\alpha -1}
\sum_{j=0}^{\infty
} (\Gamma _u(-\alpha -j\beta ))^{-1}t^{-j\beta }c_j$, \\
satisfying Properties 1(1,2).}
\par {\bf Proof.} Let a function $F(p)$ be having a decomposition
into a generalized power series $F(p)=\sum_{j=0}^{\infty } p^{j\beta
+\alpha }c_j$, where $0<\beta \in \bf Q$. Consider the loop $\gamma
$ consisting of the segment $ \gamma _1 := \{ a+\theta S, |\theta
|\le b \} $ passed from above to bottom, where $b>0$, $S\in {\cal
A}_r$, $|S|=1$, $Re (S)=0$, $\gamma _2$ is the part of the circle in
the plane ${\bf R}\oplus S\bf R$ from $-R$ to $-a+bS$ with the
center at zero, $R=(a^2+b^2)^{1/2}$, the arc $\gamma _3$ of the
circle from $a-bS$ to $-R$, two shores cut $\gamma _4\cup \gamma _5$
along the axis $\bf R$ from $-R$ to $-\delta $, $\gamma _4= [-
\delta ,-R ]+\epsilon S$, $\gamma _5= [-R,-\delta ]+\epsilon S$,
$\gamma _6$ is the arc of the circle of radius $\rho =(\delta
^2+\epsilon ^2)^{1/2}$ in the plane ${\bf R}\oplus S\bf R$ passed
clockwise from $(- \delta - \epsilon S)$ to $-\delta +\epsilon S$,
where $\epsilon >0$ and $\delta
>0$ are arbitrary small numbers, $\gamma =\bigcup_{k=1}^6\gamma _k$.
For $z$ in the plane ${\bf R}\oplus S\bf R$ we put $z=|z|\exp (S\phi
)$, where $-\pi /2<\phi <\pi /2$. In view of Theorem 2.15
\cite{ludfov,ludoyst} $\int_{a-bS}^{a+bS}F(p)\exp
(u(p,t;0)dp=\int_{\bigcup_{k=2}^6 \gamma _k}F(p)\exp (u(p,t;0))dp$.
From Lemma 23 and Note 47 it follows, that $\lim_{R\to \infty
}\int_{\gamma _2\cup \gamma _3} F(p)\exp (u(p,t;0))dp=0$ for each
$t>0$. The contour of integration $\psi $ in \S 50 is symmetrical
relative to the axis $ox$, such that while reflection in the plane
containing the axis $ox$ and perpendicular to $S$ relative to the
scalar product $(z,\eta ):=Re (z\eta ^*)$, the direction of its
circuit changes on the opposite one, therefore, $\Gamma _u(w)\in \bf
R$ for each $w\in \bf R$ like usual gamma-function.
\par From Theorem 22 it follows, that the function-original is given
by Formula:
\par $Sf(t)=(2\pi )^{-1}\int_{\gamma _4\cup \gamma _5\cup \gamma
_6}F(p)\exp (u(p,t;0)) dp$.
\par Then $ Sf(t)=\sum_{j=0}^{\infty }(2\pi )^{-1}
Re (N_1S^*) [\int_{\gamma _4\cup \gamma _5\cup \gamma _6}p^{\alpha
+j\beta }\exp (u(p,t;0)) dp] c_j$. \\
In the case, when all $c_j\in \bf R$, then $f({\bf R})\subset \bf
R$, but $Sa=b$ for $a\in \bf R$ and $b\in {\cal A}_r$ means that
$a=S^*b$, since $\bf R$ is the center of the algebra ${\cal A}_r$.
We make the substitution of the variable $\zeta =pt$, then
\par $(2\pi )^{-1}\int_{\gamma _4\cup \gamma _5\cup \gamma
_6}p^{\alpha +j\beta }\exp (u(p,t;0))dp = (2\pi )^{-1} t^{-\alpha -
j\beta -1}\int_{(\gamma _4\cup \gamma _5\cup \gamma _6)t}\zeta
^{\alpha +j\beta } \exp (u(\zeta ,1;0))d\zeta = t^{-\alpha - j\beta
-1} (\Gamma _u(-\alpha -j\beta ))^{-1}$. Consequently, for the
function $F(p)$ there exists the original given by Formula $(i)$,
since $t\in \bf R$ commutes with each $c_j\in {\cal A}_r$. If
$\alpha + j \beta $ is the integer number, then the integral of
$p^{\alpha +j\beta }\exp (u(p,t;0))$ is equal to zero, that is, from
the decomposition $(i)$ it is necessary to eliminate all terms with
integer $\alpha +j\beta $. \par It is necessary to note, that this
original in the general case may not be satisfying Property 1(3),
that is, $f(t)$ is the exceptional original.
\par {\bf 52. Example.} Consider $F(p) :=
p^{-3/2}\exp (p^{1/2}(S-1))$, where $S\in \bf K$, ${\bf K}=\bf H$ or
${\bf K}=\bf O$, $Re (S)=0$, $|S|=1$. Let $p=R\exp (S\phi )$,
$S-1=2^{1/2}\exp (3 S \pi /4)$, then $Re
(p^{1/2}(S-1)=(2R)^{1/2}\cos (\phi /2 +3\pi /4)<0$, since $\pi /2<
\phi <3\pi /2$ and $\pi <\phi /2 +3\pi /4<3\pi /2$. If $|p'|$ is
large, where $p'=p-p_0$, $p_0:=Re (p)$, and $0<p_0<a$, then $\phi $
is approximately equal to $\pi /2$ or $3\pi /2$, therefore, $p^{1/2}
(S-1)$ is approximately equal to $- 2^{1/2}R$ or $-S2^{1/2}R$. If
$p\in {\bf R}\oplus S\bf R$, then $F(p)= \sum_{j=0}^{\infty }
2^{j/2}(j!)^{-1}\exp (3\pi jS/4)p^{(j-3)/2}$. Thus, $F(p)\to 0$ for
$p\to \infty $ and $Re (p)<a$, that is, $F(p)$ satisfies Conditions
23(1,2',3) and conditions of Theorem 51, since $a>0$ and $p\ne 0$
along the straight line $\{ z=a+S\theta : \theta \in {\bf R} \}$. In
view of Theorem 51: $f(t)=\sum_{j=0}^{\infty }2^j (-S)^j
((2j)!\Gamma _u(3/2 -j))^{-1} t^{-j-1/2}$.
\par {\bf 53. Note.} As it is seen from the previous example
there exist originals not satisfying Condition 1(3), but for which
there exist images. On the other hand, not for all images there
exist originals, satisfying Conditions 1(1-3). For example,
functions $1, p, p^2,...$ do not have usual originals. To extend
possibilities of the Laplace transform there are used generalized
functions, that is, functionals.
\par {\bf 54. Definitions.} A space of test function
$\cal D$ consists of all infinite differentiable functions $f: {\bf
R}\to {\cal A}_r$ on $\bf R$ with compact supports, that is, $f$ is
zero outside a segment depending on a function. A sequence of
functions $f_n\in \cal D$ tends to zero, if all $f_n$ are zero
outside some segment $[a,b]$, while on it the sequence $f^{(k)}_n$
converges to zero uniformly for each $k=0,1,2,...$, where
$f^{(k)}(t):=d^kf(t)/dt^k$, $f^{(0)}=f$. Such convergence defines
closed subsets in $\cal D$, their complements by the definition are
open, that gives the topology on $\cal D$.
\par By a generalized function of class ${\cal D}'$ is called continuous
$\bf R$-linear ${\cal A}_r$-additive function $g: {\cal D} \to {\cal
A}_r$. The set of all such functionals is denoted by ${\cal D}'$.
That is, $g$ is continuous, if for each sequence $f_n\in \cal D$,
converging to zero, a sequence of numbers $g(f_n)=:(g,f_n)\in {\cal
A}_r$ converges to zero for $n\to \infty $. A generalized function
$g$ is zero on an open subset $V$ in $\bf R$, if $(g,f)=0$ for each
$f\in \cal D$ equal to zero outside $V$. By a support of a
generalized function $g$ is called the family, denoted by $supp
(g)$, of all points $t\in \bf R$ such that in each neighborhood of
each point $t\in supp (g)$ the functional $g$ is different from
zero. The addition of generalized functions $g, h$ is given by the
Formula: $(g+h,f):=(g,f)+(h,f)$. The multiplication $g\in {\cal D}'$
on an infinite differentiable function $w$ is given by the Formula:
$(gw,f)=(g, wf)$ for each test function $f\in \cal D$ with a real
image $f({\bf R})\subset \bf R$. By a derivative  $g'$ of a
generalized function $g$ is called a generalized function $g'$ given
by the Formula: $(g',f):= - (g,f')$.
\par The space $\cal B$ of test functions consists of all
infinite differentiable functions  $f: {\bf R}\to {\cal A}_r$ such
that there exists $\lim_{t\to +\infty } t^mf^{(j)}(t)=0$ for each
$m, j=0,1,2,...$. A sequence $f_n\in \cal B$ is called converging to
zero, if the sequence $t^mf_n^{(j)}(t)$ converges to zero uniformly
on the segment $[a,\infty )$ for each $m, j=0,1,2,...$ and each
$-\infty <a< + \infty $. The family of all $\bf R$-linear and ${\cal
A}_r$-additive functionals on $\cal B$ is denoted by ${\cal B}'$.
\par A generalized function $f\in {\cal A}'$ with a support contained in
$[0,\infty )$ we call a generalized original, if there exists a real
number $s_0$ such that for each $s>s_0$ the generalized function is
$f\exp (-st)\in {\cal B}'$. An image of such original we call a
function
\par $(1)$ ${\cal F}(f,u;p;\zeta ):=(f,\exp (-u(p,t;\zeta )))$
of the variable $p\in {\cal A}_r$, defined in the half space $Re
(p)>s_0$ by the following rule. For a given $p\in {\cal A}_r$ with
$Re (p)=s>s_0$ choose $s_1\in \bf R$ such that $s_0<s_1<s$, then
\par $(2)$  $(f,\exp (-u(p,t;\zeta
)) := (f\exp (-s_1t),\exp (-[u(p,t;\zeta )-s_1t]))$, \\
since $\exp (-[u(p,t;\zeta )-s_1t])\in {\cal B}$, by the condition
$f \exp (-s_1t)\in {\cal B}'$.
\par {\bf 55. Note and Examples.} It is evident that
${\cal F}(f,u;p;\zeta )$ does not depend on a choice of $s_1$, since
$(f\exp (-s_1t),\exp (-[u(p,t;\zeta )-s_1t]))=(f\exp (-s_1t-bt),\exp
(-[u(p,t;\zeta )-s_1t-bt]))$ for each $b\in \bf R$ such  that
$s_0<s_1+b<s$, since $\exp (-bt)\in \bf R$, while $\bf R$ is the
center of the Cayley-Dickson algebra ${\cal A}_r$, where $2\le r\in
\bf N$.
\par  Let $\delta $ be the Dirac delta function, defined by the
equation $(\delta (t),\phi (t)):=\phi (0)$ for each $\phi \in {\cal
B}$. Then
\par $(i)$ ${\cal F}(\delta ^{(j)}(t-\tau ),u;p;\zeta )
=(\delta ^{(j)}\exp (-s_1t),\exp (-[u(p,t;\zeta )-s_1t]))$ \\  $=
(-1)^j(d^j \exp (-[u(p,t;\zeta )]/ dt^j)|_{t=\tau }$, \\
since it is possible to take $-\infty <s_0<0$ and $s_1=0$, where
$\tau \in \bf R$ is the parameter. In particular, for $j=0$ we have
${\cal F}(\delta (t-\tau ),u;p;\zeta ) = \exp (-u(p,\tau ;\zeta ))$,
while for $u=pt$ we have ${\cal F}(\delta ^{(j)}(t-\tau ),pt;p;\zeta
)= p^n\exp (-p\tau )$.
In the general case: \\
\par ${\cal F}(\delta ^{(j)}(t),u;p;\zeta )=$ \par
$\sum_{n_0,n_1,...,n_{2^r-1}=j} p_0^{n_0}p_1^{n_1}...
p_{2^r-1}^{n_{2^r-1}}\exp (- \zeta _0 - M(p,0;\zeta - (i_1n_1+...
+i_{2^r-1}n_{2^r-1})\pi /2))$,\\
moreover, $M(p,0;\zeta )=M(0,0;\zeta )$, where $n_0, n_1,...,
n_{2^r-1}$ are nonnegative numbers.
\par The transformation ${\cal F}$ of the generalized function is
a holomorphic function by $p\in {\cal A}_r$ with $Re (p)>s_0$ and by
$\zeta \in {\cal A}_r$, since the right side of Equation 54(2) is
holomorphic by $p$ with $Re (p)>s_0$ and by $\zeta $ in view of
Theorem 7. From Equation 54(2) it follows, that Theorems 27, 28 and
29 are accomplished also for generalized  functions.
\section{Two-sided noncommutative integral transformations.}
\par {\bf 1. Definition.} Consider function-originals, satisfying
conditions $(1-3)$ below:
\par $(1)$ $f(t)$ satisfies the H\"older condition: $|f(t+h)-f(t)|
\le A |h|^{\alpha }$ for each $|h|<\delta $ (where $0<\alpha \le 1$,
$A=const >0$, $\delta >0$ are constants for a given $t$) everywhere
on $\bf R$ may be besides points of discontinuity of the first kind.
On each finite interval in $\bf R$ a function $f$ may have only a
finite number of points of discontinuity and of the first kind only.
\par $(2)$ $|f(t)| < C_1 \exp (-s_1t)$ for each $t<0$,
where $C_1=const >0$, $s_1=s_1(f)=const \in \bf R$.
\par $(3)$ $|f(t)|<C_2 \exp (s_0t)$ for each $t\ge 0$, that is,
$f(t)$ is growing not faster, than the exponential function, where
$C_2=const>0$, $s_0=s_0(f)\in \bf R$.
\par The two-sided Laplace transformation over the Cayley-Dickson
algebras ${\cal A}_r$ with $2\le r\in \bf N$ is defined analogously
to that of one-sided, but with the substitution of the lower
integration limit on $-\infty $ instead of zero:
\par $(4)$ ${\cal F}^s(f,u;p;\zeta ) :=
\int_{-\infty }^{\infty }f(t)\exp (-u(p,t;\zeta ))dt$ \\
for all numbers $p\in {\cal A}_r$, for which the integral exists.
Denote for short ${\cal F}^s(f,u;p;\zeta )$ through $F^s_u(p;\zeta
)$ (see also Definitions 2.1 and 2.3). For a basis of generators $\{
N_0,...,N_{2^r-1} \} $ in ${\cal A}_r$ we shall write in more
details $\mbox{ }_N{\cal F}^s(f,u;p;\zeta )$ or $\mbox{
}_NF^s_u(p;\zeta )$ in the case of necessity.
\par {\bf 2. Note.} Naturally, that the two-sided Laplace integral
can be considered as the sum of two one-sided integrals
\par $(1)$  $\int_{-\infty }^{\infty }f(t)\exp (-u(p,t;\zeta ))dt=
\int_{-\infty }^0f(t)\exp (-u(p,t;\zeta ))dt +
\int_0^{\infty }f(t)\exp (-u(p,t;\zeta ))dt$ \\
$=\int_0^{\infty }f(-t)\exp (-u(p,-t;\zeta ))dt +
\int_0^{\infty }f(t)\exp (-u(p,t;\zeta ))dt$. \\
Let $u$ has the form 2.1 or 2.3(1-4). The second integral converges
for $Re (p)>s_0$. Since $u(-p,-t;\zeta ) = u(p,t;\zeta )$, then the
first integral converges for $Re (-p)> -s_1$, that is, for $Re (p)<
s_1$. Then there is a region of convergence $s_0 < Re (p) < s_1$ of
the two-sided Laplace integral. For $s_1=s_0$ the region of
convergence reduces to the vertical hyperplane in ${\cal A}_r$ over
$\bf R$. For $s_1<s_0$ there is no any common domain of convergence
and $f(t)$ can not be transformed with the help of the two-sided
transformation 1(4).
\par {\bf 3. Example.} ${\cal F}^s(\exp (-\alpha t^2),pt;p;0)=
\int_{-\infty }^{\infty } \exp (-\alpha t^2-pt)dt= (\pi/\alpha
)^{1/2}\exp (p^2/(4\alpha ))$, where $\alpha >0$, since
$\int_{-\infty }^{\infty }\exp (-t^2)dt=(\pi )^{1/2}$. For a
comparison the one-sided Laplace transformation gives:
\par ${\cal F}(\exp (-\alpha t^2)Ch_{[0,\infty )},pt;p;0)=
\int_0^{\infty } \exp (-\alpha t^2-pt)dt$ \par  $= (\alpha
)^{-1/2}\exp (p^2/(4\alpha ))\int_{p/(2(\alpha )^{1/2})}^{\infty
}\exp (-t^2)dt=2^{-1}(\pi /\alpha )^{1/2}\exp (p^2/(4\alpha )) Erf
(p/(2(\alpha )^{1/2}))$ (see also \S 1.43).
\par The application of Theorem 1.7 to
$\int_0^{\infty }f(-t)\exp (-u(p,-t))dt$ and $\int_0^{\infty
}f(t)\exp (-u(p,t))dt$ gives.
\par {\bf 4. Theorem.} {\it If an original $f(t)$ satisfies
Conditions 1(1-3), and moreover, $s_0<s_1$, then its image ${\cal
F}^s(f,u;p;\zeta )$ is holomorphic by $p$ in the domain $\{ z\in
{\cal A}_r: s_0< Re (z)<s_1 \} $, as well as by $\zeta \in {\cal
A}_r$, where $2\le r\in \bf N$.}
\par {\bf 5. Proposition.} {\it If $\mbox{ }_NF^s_u(p;\zeta )$ and
$\mbox{ }_NG^s_u(p;\zeta )$ are images of the original functions
$f(t)$ and $g(t)$ in the domains $s_0(f)< Re (p)< s_1(f)$ and
$s_0(g)<Re (p)<s_1(g)$ with values in ${\cal A}_r$, then for each
$\alpha , \beta \in {\cal A}_r$ in the case ${\bf K}=\bf H$; as well
as $f$ and $g$ with values in $\bf R$ and each $\alpha , \beta \in
{\cal A}_r$ or $f$ and $g$ with values in ${\cal A}_r$ and each
$\alpha , \beta \in \bf R$ in the case of ${\cal A}_r$ the function
$\alpha \mbox{ }_NF^s_u(p;\zeta ) + \beta \mbox{ }_NG^s_u(p;\zeta )
$ is the image of the function $\alpha f(t) +\beta g(t)$ in the
domain $\max (s_0(f),s_0(g))< Re (p)< \min (s_1(f),s_1(g))$.}
\par {\bf Proof} follows from Remark 2 applying Proposition 2.25 to each
of the integrals \par $\alpha \int_0^{\infty }f(-t)\exp
(-u(p,-t;\zeta
))dt + \beta \int_0^{\infty }g(-t)\exp (-u(p,-t;\zeta ))dt$ and \\
$\alpha \int_0^{\infty }f(t)\exp (-u(p,t;\zeta ))dt + \beta
\int_0^{\infty }g(t)\exp (-u(p,t;\zeta ))dt$, since their common
domain of convergence is the domain $\max (s_0(f),s_0(g))< Re (p)<
\min (s_1(f),s_1(g))$, where $p\in {\cal A}_r$.
\par {\bf 6. Examples.} 1. There may be cases, when a domain of
convergence for a sum is greater, then for each additive. For
example, ${\cal F}^s(\exp (at) U(t),pt;p;0)=(p-a)^{-1}$, also ${\cal
F}^s((\exp (at) -1) U(t),pt;p;0)=ap^{-1}(p-a)^{-1}$ for $Re (p)>a$
in both cases, when $a\in \bf R$, $U(t):=1$ for $t>0$, $U(0)=1/2$,
while $U(t)=0$ for $t<0$. But $(p-a)^{-1}+ ap^{-1}(p-a)^{-1}=p^{-1}$
and ${\cal F}^s(U(t),pt;p;0)=p^{-1}$ for each $Re (p)>0$.
\par It is necessary to note, that the two-sided Laplace
transformation of the function $t^n$ does not exist, but the
one-sided transformation was elucidated in examples 2.30.1 and 2.33.
\par 2. ${\cal F}^s(\exp (-\alpha |t|)/2,pt;p;0)=\alpha
(\alpha ^2-p^2)^{-1}$ in the domain $-\alpha <Re (p)<\alpha $ for
$\alpha >0$.
\par 3. Consider the two-sided transformation
\par ${\cal F}^s((e^t+1)^{-1},pt;p;0) = \int_{-\infty }^{\infty }
(e^t+1)^{-1}\exp (-pt)dt$ \\ in the domain $-1<Re (p)<0$. Make the
substitution $v=(e^t+1)^{-1}$, then the integral reduces to the
Euler integral of the first kind $\int_0^1v^p(1-v)^{-p-1}dv = -\pi
/\sin (\pi p)$ (see Proposition 4.6, Definition 4.14 and Theorem
4.17 in \cite{ludfov}).
\par {\bf 7. Theorem.} {\it If ${\cal F}^s(f(t),u;p;\zeta )$
is defined in the domain $s_0< Re (p) <s_1$, then
\par ${\cal F}^s(f(bt),u;p;\zeta )={\cal F}^s(f(t),u;p/b;\zeta )/b$
for $b>0$ in the domain $bs_0< Re (p)<bs_1$,\par  ${\cal
F}^s(f(bt),u;p;\zeta )= - {\cal F}^s(f(t),u;p/b;\zeta )/b$ for $b<0$
in the domain $bs_1< Re (p)<bs_0$.}
\par {\bf Proof} follows from the definition of the two-sided
transformation with the subsequent substitution $v=bt$ of the
variable of integration, since $p_jt=(p_j/b)(bt)$ for each $j$ and
$u(p,t;\zeta ) = u(p/b,bt;\zeta )$ as for $u=pt+\zeta $, and for $u=
p_0t+M(p,t;\zeta )+\zeta _0$ as well (see \S 2.3).
\par {\bf 8. Theorem.} {\it If $f(t)$ is a function-original with
parameters $s_0, s_1$, then
\par $(1)$  ${\cal F}^s((f(t-\tau ),u;p;\zeta )={\cal F}^s(
(f(t),u;p; \zeta +p\tau )$ for $u(p,t;\zeta )= p_0t + \zeta _0 +
M(p,t;\zeta )$ or $u(p,t;\zeta )= pt + \zeta $ over ${\cal A}_r$
with $2\le r<\infty $,
\par $(2)$ ${\cal F}^s((f(t-\tau ),pt;p;0)={\cal F}^s((f
(t),pt;p;0)e^{-p\tau }$ over either ${\bf K}=\bf H$ or ${\bf K}=\bf
O$ with $f({\bf R})\subset \bf K$, or over ${\cal A}_r$ with $4\le
r$ for $f({\bf R})\subset \bf R$, in the domain $s_0< Re (p)< s_1$.}
\par {\bf 9. Theorem.} {\it If $f(t)$ is a function-original
with values in ${\cal A}_r$ for $2\le r<\infty $, $b\in \bf R$, then
${\cal F}^s(e^{bt}f(t),u;p;\zeta )= {\cal F}^s(f(t),u;p-b;\zeta )$
for each $s_0+b<Re (p)<s_1+b$, where $u=pt+\zeta $ or $u=p_0t+\zeta
_0 + M(p,t;\zeta )$.}
\par {\bf Proof} follows by application of Theorem 2.37 or
2.39 respectively to each integral of the form $\int_0^{\infty
}f(-t)\exp (-u(p,-t;\zeta ))dt$ or $\int_0^{\infty }f(t)\exp
(-u(p,t;\zeta ))dt$ (see Note 2), since $s_0(e^{bt}f(t))=s_0(f)+b$ и
$s_1(e^{bt}f(t))=s_1+b$.
\par {\bf 10. Theorem.} {\it Let  $f: [a,b]\to {\cal A}_r$
and $\alpha : [a,b]\to \bf R$ be functions of bounded variations,
moreover, let $\alpha $ be continuous, where $a<b\in \bf R$. Then
\par $\int_a^bfd\alpha = f(b)\alpha (b) -f(a)\alpha (a)
- \int_a^b\alpha df$, \\
where $\int_a^bfd\alpha $ is the Riemann-Stieltjes integral.}
\par {\bf Proof.} Write the function $f$ in the form $f=
\sum_{p=0}^{2^r-1}f_pi_p$, where each function $f_p$ takes only real
values. Then it is sufficient to prove the theorem for each $f_p$,
since $\alpha $ commutes with $i_p$ for each $p$. Since $f_p$ is of
bounded variation, then there exists the Riemann-Stieltjes integral
$\int_a^b\alpha df_p$ for each $p$ and $\sum_p(\int_a^b\alpha
df_p)i_p=\int_a^b\alpha df$. For an arbitrary partition $P = \{ x_0,
x_1,..., x_n \} $ of the segment $[a,b]$ we choose points
$t_0,...,t_{n+1}$ such that $t_0=a$, $t_{n+1}=b$, $x_{j-1}\le t_j\le
x_j$ for each $j=1,...,n$, then we get the partition $Q = \{
t_0,...,t_{n+1} \} $ of the segment $[a,b]$. The summation by parts
gives the equality: \par $S(P,\alpha ,f_p) := \sum_{j=1}^n \alpha
(t_j)[f_p(x_j)-f_p(x_{j-1})] = f_p(b)\alpha (b) -f_p(a)\alpha (a) -
\sum_{j=1}^{n+1} f_p(x_{j-1})[\alpha (t_j)-\alpha (t_{j-1})]=
f_p(b)\alpha (b) -f_p(a)\alpha (a) - S(Q,f_p,\alpha )$, \\
where $S(Q,f_p,\alpha ) := \sum_{j=1}^{n+1} f_p(x_{j-1})[\alpha
(t_j)-\alpha (t_{j-1})]$, since $t_{j-1}\le x_{j-1}\le x_j$. If the
diameter of the partition $diam (P)\to 0$ tends to zero, then $diam
(Q)\to 0$ also, where $diam (P) := \max_{1\le j\le n}
|x_{j-1}-x_j|$. Therefore,
\par $\lim_{diam (P)\to 0}S(P,\alpha ,f_p) = \int_a^b\alpha df_p$, and
\par $\lim_{diam (Q)\to 0}S(Q,f_p,\alpha ) = \int_a^bf_pd\alpha $ \\
due to Theorem 6.14 \cite{rudin}, which states that if there exists
$\lim_{diam (Q)\to 0} S(Q,g,\alpha )$, then $g$ is Riemann-Stieltjes
integrable relative to $\alpha $ and $\lim_{diam (Q)\to 0}
S(Q,g,\alpha ) = \int_a^b g d\alpha $.
\par {\bf 11. Theorem.} {\it If $f(t)$ is an original function,
$F(p)$ is its image, $f^{(n)}(t)$ is an original for $n\ge
1$, then \par $(i)$ ${\cal F}^s(f^{(n)}(t),pt;p;0) = F(p)p^n$ \\
is the image of the function $f^{(n)}(t)$ over ${\cal A}_r$ for
$u=pt$, if $f$ is real-valued for $r\ge 4$, $f$ is with values in
$\bf K$ for ${\bf K}=\bf H$ or ${\bf K}=\bf O$;
\par $(ii)$  ${\cal F}^s(f'(t),u; p;\zeta ) =
p_0{\cal F}^s(f(t),u; p;\zeta ) $ \\
$+p_1{\cal F}^s(f(t),u; p;\zeta - i_1\pi /2)+...+p_{2^r-1}{\cal
F}^s(f(t),u; p;\zeta - i_{2^r-1}\pi /2)$ \\
for $u(p,t) := p_0t + M(p,t;\zeta )+\zeta _0$ over ${\cal A}_r$ with
$2\le r<\infty $. Domains, where Formulas $(i,ii)$ are true may be
different from a domain of the two-sided Laplace transformation for
$f$, but they are satisfied in the domain $s_0< Re (p)<s_1$, where
$s_0= \max (s_0(f),s_0(f'),...,s_0(f^{(n)}))$, $s_1=\min
(s_0(f),s_0(f'),...,s_0(f^{(n)}))$, if $s_0<s_1$.}
\par {\bf Proof.} In view of Theorem 10 the integration by parts
on the segment $a\to -\infty $, $b\to \infty $ gives:
\par $\int_{-\infty }^{\infty }f'(t)e^{-pt}dt = [f(t)e^{-pt}]|_{
-\infty }^{\infty } + (\int_{-\infty }^{\infty }f(t)e^{-pt}dt)p$ \\
due to the alternativity of ${\bf K}=\bf H$ and ${\bf K}=\bf O$ for
$f$ with values in $\bf K$, also for $r\ge 4$ due to real-valuedness
of $f$. For $s_0<Re (p)<s_1$ the first additive is zero. The
application of the last formula $n$ times gives the general Formula
$(i)$.
\par For $u= p_0t+\zeta _0+M(p,t; \zeta )$ in view of Theorem 10
\par $\int_{-\infty }^{\infty }f'(t)e^{-u(p,t;\zeta )}dt =
[f(t)e^{-u(p,t;\zeta )}]|_{-\infty }^{\infty } + (\int_{-\infty
}^{\infty }f(t)[\partial e^{-u(p,t;\zeta )}/\partial t] dt)$. \\
Again for $s_0<Re (p)<s_1$ the first additive is zero, while the
second integral converts with the help of Formula 2.31(2), from
which Formula (ii) follows.
\par For a derivation it can be used Theorem 8, since
\par $\lim_{\tau \to 0}[{\cal F}^s(f(t),u;p;\zeta ) -
{\cal F}^s(f(t-\tau ),u;p;\zeta )]/\tau = \lim_{\tau \to 0}[{\cal
F}^s(f(t),u;p;\zeta ) - {\cal F}^s(f(t),u;p;\zeta +p\tau )]/\tau
=\lim_{\tau \to 0}\int_{-\infty }^{\infty }f'(t)[e^{-u(p,t;\zeta )}
- e^{-u(p,t;\zeta +p\tau )}]\tau ^{-1}dt$. \\
The change of places of $\lim_{\tau \to 0} $ and $\int_{-\infty
}^{\infty }$ may change a domain of convergence, but in the
indicated in the theorem domain $s_0< Re (p)<s_1$, when it is non
void, the application of Statement 6 from \S XVII.2.3 \cite{zorich}
and Note 2 gives Formula (i,ii) of the given theorem, since if there
exists an original $f^{(n)}(t)$, then $f^{(n-1)}(t)$ is continuous
for $n\ge 1$, where $f^0:=f$.
\par {\bf 12. Theorem.} {\it If $f(t)$ is an original function
for ${\cal A}_r$ with $2\le r<\infty $, then
\par $(1)$  $({\cal F}^s(f,u;p;0))^{(n)}(p).(\mbox{ }_
1h,...,\mbox{ }_nh) = {\cal F}^s((-t)^nf(t)\mbox{ }_1h...\mbox{
}_nh;u;p;0)$ for each $n\in \bf N$ and $\mbox{ }_1h,...,\mbox{
}_nh\in {\bf R}\oplus p'{\bf R}\subset {\cal A}_r$ for $u=pt$, where
$p' := p-Re (p)$. If $u(p,t) := p_0t + M(p,t;\zeta )+\zeta _0$, then
\par $(2)$  $(\partial {\cal F}^s(f(t),u; p;\zeta )/
\partial p).h = - {\cal F}^s(f(t)t,u;p;\zeta )h_0 $ \\
$- {\cal F}^s(f(t)t,u;p;\zeta - i_1\pi /2)h_1 -...-{\cal F}(f(t)t,u;
p; \zeta - i_{2^r-1}\pi /2)h_{2^r-1}$ \\
for each $h=h_0i_0+...+h_{2^r-1}i_{2^r-1}\in {\cal A}_r$, where
$h_0,...,h_{2^r-1}\in {\bf R}$. Both formulas are true in the domain
$s_0(f)<Re (p)<s_1(f)$.}
\par The {\bf proof} is quite analogous to the proof of Theorems
2.29 and 2.32 with the substitution of $\int_0^{\infty }$ on
$\int_{-\infty }^{\infty }$ and with the use of Note 2, since
$s_0(f)<Re (p)<s_1(f)$ is equivalent to the inequalities
$s_0(f(t)t)<Re (p)<s_1(f(t)t)$, since $\lim_{t\to + \infty }\exp
(-bt)t=0$ for each $b>0$, also $s_0(f)\le s_0(f(t)t)<s_0(f)+b$ for
$t>0$ and $s_1(f)\le s_1(f(t)t)<s_1(f)+b$ for $t<0$ for each $b>0$.
Indeed, ${\cal F}^s(f(t),u; p;\zeta )$ is the holomorphic function
by $p$ for $s_0(f)<Re (p)<s_1(f)$, also $|\int_0^{\infty
}e^{-ct}t^ndt|<\infty $ for each $c>0$ and $n=0,1,2,...$, then it is
possible to differentiate under the sign of the integral:
\par $(\partial (\int_{-\infty }^{\infty }f(t)\exp (-u(p,t;\zeta
))dt)/\partial p).h =$ \\ $(\partial (\int_{-\infty }^0f(t)\exp
(-u(p,t;\zeta ))dt)/\partial p).h + (\partial
(\int_0^{\infty }f(t)\exp (-u(p,t;\zeta ))dt)/\partial p).h=$ \\
$\int_{-\infty }^0f(t)(\partial \exp (-u(p,t;\zeta ))/\partial
p).hdt + \int_0^{\infty }f(t)(\partial \exp (-u(p,t;\zeta
))/\partial p).hdt$ \\ $=
 \int_{-\infty }^{\infty }f(t)(\partial \exp (-u(p,t;\zeta
))/\partial p).hdt$.
\par {\bf 13. Theorem.} {\it If $f(t)$ and
$g(t)$ are original functions, where $g$ is real-valued for each
$t$, also the convolution $f*g(t):= \int_{-\infty }^{\infty }f(\tau
)g(t-\tau )dx$ is the original and
\par $(1)$ $f$, $g$ and $f*g(t)$ have the converging two-sided
Laplace transformations in a common domain of convergence or
\par $(2)$ $f$ or $g$ has an absolutely converging two-sided Laplace
integrals, then
\par ${\cal F}^s(\int_{-\infty }^{\infty }
f(\tau )g(t-\tau )d\tau , pt;p;0) = {\cal F}^s(f(t),pt;p;0) {\cal
F}^s(g(t),pt;p;0)$ \\
for each $p\in {\cal A}_r$ with $s_0<Re (p)<s_1$, where $s_0 = \max
(s_0(f), s_0(g))$, $s_1=\min (s_1(f),s_1(g))$, and $2\le r<\infty
$.}
\par {\bf Proof.} Transform the product with the help of the
definition: \par  ${\cal F}^s(f(t),pt;p;0) {\cal F}^s(g(t),pt;p;0)
=\int_{-\infty }^{\infty }{\cal F}^s(f(t),pt;p;0)g(\tau ) \exp (-pt)
d\tau $ in the domain $s_0<Re (p)<s_1$. \\
Since $g(\tau )$ and $\exp (-p\tau )$ commute, then in view view of
Theorem 8
\par ${\cal F}^s(f(t),pt;p;0)\exp (-p\tau ) = {\cal
F}^s(f(t-\tau ),pt;p;0)$, consequently,
\par ${\cal F}^s(f(t),pt;p;0) {\cal F}^s(g(t),pt;p;0)
=\int_{-\infty }^{\infty }{\cal F}^s(f(t-\tau ),pt;p;0)g(\tau )d\tau
$ \par $=\int_{-\infty }^{\infty }(\int_{-\infty }^{\infty }f(t-\tau
)\exp (-pt)dt)g(\tau )d\tau $. \\
If there is accomplished Condition $(1)$ or $(2)$, then using the
expression of a ${\cal A}_r$ valued function in the form
$f=\sum_vf_vi_v$ the order of integration can be changed in the last
integral in view of the Fubini theorem (see also Proposition 2.18
above), from which the statement of this theorem follows
\par ${\cal F}^s(f(t),pt;p;0) {\cal F}^s(g(t),pt;p;0) =
\int_{-\infty }^{\infty }(\int_{-\infty }^{\infty }f(t-\tau )g(\tau
)d\tau )\exp (-pt)dt$.
\par {\bf 14. Theorem.} {\it Let $f(t)$ be an original function
with values in ${\cal A}_r$ for $2\le r<\infty $, $u=pt$, moreover,
$s_0(f)\le 0\le s_1(f)$, then
\par $(1)$ ${\cal F}^s(g(t),pt;p;0)p =
{\cal F}^s(f(t),u;p;0)$ \\
for $g(t) := \int_{-\infty }^t f(x)dx$ in the domain $\max
(s_0(f),0)<Re (p)<s_1(f)$, in particular,
\par ${\cal F}^s(g(t),pt;p;0) =
{\cal F}^s(f(t),u;p;0)p^{-1}$ over the algebra ${\bf K}=\bf H$ or
${\bf K}=\bf O$;
\par $(2)$ ${\cal F}^s(g(t),pt;p;0)p =
{\cal F}^s(f(t),u;p;0)$ \\
for $g(t) := \int_{+\infty }^t f(x)dx$ in the domain $s_0(f)<Re
(p)<\min (s_1(f),0)$, in particular,
\par ${\cal F}^s(g(t),pt;p;0) =
{\cal F}^s(f(t),u;p;0)p^{-1}$ over the algebra ${\bf K}=\bf H$ or
${\bf K}=\bf O$, when these domains are non void.}
\par {\bf Proof.} For $s>0$ we have the equality
$\int_0^T\exp (st)dt=(\exp (sT) -1)/s$, also for $s\le 0$ we have
the inequality $\int_0^T\exp (st)dt\le T$. Moreover, $\int_{-\infty
}^0 \exp (ct)dt$ converges for $c<0$. Then $s_0(g)<\max
(s_0(f),0)+b$ for each $b>0$ when $g(t) := \int_{-\infty }^t
f(x)dx$.
\par For $s>0$ we have the equality $\int_{+\infty }^T
\exp (-st)dt= -\exp (-sT)/s$ for each $T\in \bf R$, in particular
for $T<0$. Moreover, $\int_{+\infty }^0 \exp (-ct)dt$ converges for
$c>0$. Then $s_1(g)>\min (s_1(f),0)-b$ for each $b>0$ when $g(t) :=
\int_{+\infty }^t f(x)dx$. Therefore, the first and the second
statements of this theorem follow from Theorem 11 above. In
particular, in the algebra ${\bf K}=\bf H$ or ${\bf K}=\bf O$ the
equation $zp=y$ is resolvable relative to $z$ for $p\ne 0$ such that
$z=yp^{-1}$, then Formulas $(1,2)$ converts to the forms indicated
above.
\par {\bf 15. Theorem.} {\it If a function $f(t)$ is an original
such that \par $\mbox{ }_N{\cal F}^s(f,u;p;\zeta ) :=
\sum_{j=0}^{2^r-1} \mbox{ }_NF^s_{u,j}(p;\zeta )N_j$ is its image,
where a function $f$ is written in the form \par $f(t) =
\sum_{j=0}^{2^r-1} f_j(t)N_j$, $f_j: {\bf R}\to \bf R$ for each
$j=0,1,...,2^r-1$, $f({\bf R})\subset \bf K$ for ${\bf K}=\bf H$ or
${\bf K}=\bf O$, $f({\bf R})\subset \bf R$ over the Cayley-Dickson
algebra ${\cal A}_r$ with $4\le r\in \bf N$,
\par $\mbox{ }_NF^s_{u,j}(p;\zeta ) := \int_{-\infty }^{\infty
}f_j(t)\exp (-u(p,t;\zeta ))dt$. Then at each point $t$, where
$f(t)$ satisfies the H\"older condition there is true the equality:
\par $(i)$ $f(t) = (2\pi N_1)^{-1} Re (S{\tilde N}_1) \sum_{j=0}^{2^r-1}
\int_{a-S\infty }^{a+S\infty }\mbox{ }_NF^s_{u,j}(p;\zeta )
\exp (u(p,t;\zeta ))dp)N_j$ \\
in the domain $s_0(f)< Re (p) <s_1(f)$, where $u(p,t;\zeta
)=pt+\zeta $ or $u(p,t;\zeta )=p_0t+M_N(p,t;\zeta )+\zeta _0$ and
the integral is taken along the straight line $p(\tau )=a+S\tau \in
{\cal A}_r$, $\tau \in \bf R$, $S\in {\cal A}_r$, $Re (S)=0$,
$|S|=1$, while the integral is understood in the sense of the
principal value.}
\par {\bf Proof.} The two-sided transformation in the basis of generators
\\ $N = \{ N_0,N_1,...,N_{2^r-1} \} $ can be written in the form
\par ${\cal F}^s(f,u;p;\zeta ) := \int_{-\infty }^{\infty }f(t)\exp
(-u(p,t))dt = {\cal F}^s(fU(t),u;p;\zeta ) + {\cal
F}^s(f(1-U)(t),u;p;\zeta ) $, \\ where $u=u(p,t;\zeta )$, the index
$N$ is omitted, also
\par ${\cal F}^s(f(1-U)(t),u(p,t;\zeta );p;\zeta ) = \int_0^{\infty
}f(-t)U(-t)\exp (-u(p,-t;\zeta ))dt$ \\  $= {\cal
F}^s(f(-t),u(p,-t;\zeta );p;\zeta )$, \\ where $|f(-t)|\le
C_1\exp (s_1t)$ for each $t>0$. The common domain of the existence
\\ ${\cal F}^s(f(-t),u(p,-t;\zeta );p;\zeta )$ and ${\cal
F}^s(fU(t),u;p;\zeta )$ is $s_0(f)<Re (p)<s_1(f)$, since $p_j(-t) =
-p_jt$, also the inequality $Re (-p)> - s_1(f)$ is equivalent to the
inequality $Re (p)<s_1(f)$. Then the application of Theorem 2.19
twice to $f(t)U(t)$ and to $f(-t)U(-t)$ gives the statement of this
theorem.
\par {\bf 16. Theorem.} {\it If a function $\mbox{
}_NF^s_u(p)$ is analytic by the variable $p\in {\cal A}_r$ in the
domain $W := \{ p\in {\cal A}_r: s_0< Re (p) <s_1 \} $, where $2\le
r\in \bf N$, $f({\bf R})\subset \bf K$ for ${\bf K}=\bf H$ or ${\bf
K}=\bf O$, also for an arbitrary algebra ${\cal A}_r$ with $f({\bf
R})\subset \bf R$ for $4\le r<\infty $, $u(p,t;\zeta )=pt+\zeta $ or
$u(p,t;\zeta ) := p_0t + M(p,t;\zeta )+\zeta _0$. Let also $\mbox{
}_NF^s_u(p)$ can be written in the form $\mbox{ }_NF^s_u(p)=\mbox{
}_NF^{s,0}_u(p) + \mbox{ }_NF^{s,1}_u(p)$, where $\mbox{
}_NF^{s,0}_u(p)$ is holomorphic by $p$ in the domain $s_0<Re (p)$,
also $\mbox{ }_NF^{s,1}_u(p)$ is holomorphic by $p$ in the domain
$Re (p)<s_1$, moreover, for each $a>s_0$ and $b<s_1$ there exists
constants $C_a>0$, $C_b>0$ and $\epsilon _a
>0$ and $\epsilon _b>0$ such that
\par $(i)$ $|\mbox{ }_NF^{s,0}_u(p)|\le C_a\exp (-\epsilon _a |p|)$
for each $p\in {\cal A}_r$ with $Re (p)\ge a$,
\par $(ii)$ $|\mbox{ }_NF^{s,1}_u(p)|\le C_b\exp (-\epsilon _b |p|)$
for each $p\in {\cal A}_r$ with $Re (p)\le b$, where $s_0$ and $s_1$
are fixed, also the integral
\par $(iii)$ $\int_{w-S\infty }^{w+S\infty }\mbox{
}_NF^{s,k}_u(p)dp$ \\ converges absolutely for $k=0$ and $k=1$ for
$s_0<w<s_1$, then $\mbox{ }_NF^s_u(p)$ is the image of the function
\par $(iv)$ $f(t)=(2\pi )^{-1}{\tilde S}\int_{w-S\infty
}^{w+S\infty }\mbox{ }_NF^s_u(p)\exp (u(p,t;0))dp$.}
\par {\bf Proof.} For the function
$\mbox{ }_NF^{s,1}_u(p)$ we consider the substitution of the
variable $p=-g$, $-s_1<Re (g)$. In view of Theorem 2.22 there exists
originals $f^0$ and $f^1$ for functions $\mbox{ }_NF^{s,0}_u(p)$ and
$\mbox{ }_NF^{s,1}_u(p)$ while a choice of $w\in \bf R$ in the
common domain $s_0<Re (p)<s_1$, that is, $s_0<w<s_1$. At the same
time the supports of the functions $f^0$ and $f^1$ are contained in
$[0,\infty )$ and $(-\infty ,0]$ respectively. Then $f=f^0+f^1$ is
the original for $\mbox{ }_NF^s_u(p)$ while $\zeta =0$, since
\par $f(t)=f^0(t)+f^1(t)=(2\pi )^{-1}{\tilde S}\int_{w-S\infty
}^{w+S\infty }\mbox{ }_NF^{s,0}_u(p)\exp (u(p,t;0))dp+ $\\  $(2\pi
)^{-1}{\tilde S}\int_{w-S\infty }^{w+S\infty }\mbox{
}_NF^{s,1}_u(p)\exp (u(p,t;0))dp= (2\pi )^{-1}{\tilde
S}\int_{w-S\infty }^{w+S\infty }\mbox{ }_NF^s_u(p)\exp (u(p,t;0))dp$
\\ due to the distributivity of the multiplication in the algebra
${\cal A}_r$.
\par {\bf 17. Note.} While the definition of the one- and two-sided
Laplace transformations over the Cayley-Dickson algebras above the
Riemann integral of the real variable was used, while for the
inverse transformation the noncommutative integral along paths over
${\cal A}_r$ from the works \cite{ludfov,ludoyst}. It can be
considered also a generalization of the direct transformation with
the Riemann-Stieltjes integral as the starting point. For a function
$\alpha (t)$ with values in ${\cal A}_r$ of the variable $t\in \bf
R$ such that $\alpha (t)$ has a bounded variation on each finite
segment $[a,b]\subset \bf R$, we consider the Stieltjes integral
\par $\int_{-\infty }^{\infty }(d\alpha (t))\exp (-u(p,t;\zeta )):=$
\\ $\lim_{b\to \infty } \int_0^b(d\alpha (t))\exp (-u(p,t;\zeta )) +
\lim_{b\to \infty } \int_{-b}^0(d\alpha (t))\exp (-u(p,t;\zeta )) $,
where \par $\int_a^b(d\alpha
(t))f(t)=\sum_{v,w}(\int_a^bf_v(t)d\alpha
_w(t)) (i_wi_v)$, \\
$i_0,i_1,...,i_{2^r-1}$ are generators of the Cayley-Dickson algebra
${\cal A}_r$, $f_v$ and $\alpha _w$ are real-valued functions such
that $\alpha =\sum_w\alpha _wi_w$ and $f=\sum_vf_vi_v$, also
$\int_a^b f_v(t)d\alpha _w(t)$ is the usual Stieltjes integral over
the field of real numbers on a finite segment $[a,b]$. Under
imposing the condition $\alpha (t)\exp (-p_0t)|_{-\infty }^{\infty
}=0$ the integration by parts gives the relation
\par $\int_{-\infty }^{\infty }(d\alpha (t))\exp (-u(p,t;\zeta ))=
\int_{-\infty }^{\infty }\alpha (t)d[\exp (-u(p,t;\zeta ))] $. For
$u=pt+\zeta $ over $\bf H$ and $\bf O$ in view of their
alternativity it gives:
\par $(1)$ $\int_{-\infty }^{\infty }(d\alpha (t))\exp (-u(p,t;0))=
{\cal F}^s(\alpha (t)p,u;p;0)$. \\
This relation is also true for real valued function $\alpha $ over
the general Cayley-Dickson algebra. In general in view of Formula
31(2):
\par $(2)$  $\int_{-\infty }^{\infty }(d\alpha (t))\exp (-u(p,t;0))
=p_0{\cal F}^s(\alpha (t),u; p;\zeta ) $ \\
$+p_1{\cal F}^s(\alpha (t),u; p;\zeta - i_1\pi
/2)+...+p_{2^r-1}{\cal
F}^s(\alpha (t),u; p;\zeta - i_{2^r-1}\pi /2)$ \\
for $u(p,t) := p_0t + M(p,t;\zeta )+\zeta _0$ over ${\cal A}_r$ with
$2\le r<\infty $.
\par Thus, the Laplace transformation over ${\cal A}_r$
can be spread on a more general class of originals. For this it is
used instead of an ordinary notion of convergence of improper
integrals their convergence by Cesaro of order $p>0$:
\par $(C,p) \int_0^{\infty }f(t)dt := \lim_{b\to \infty }
\int_0^bf(t)(1-t/b)^pdt$. \\
If this integral converges by Cesaro for some $p>0$, then it
converges for each $q>p$, moreover,
\par $(C,p) \int_0^{\infty }f(t)dt := (C,q)\int_0^{\infty }f(t)dt$. \\
That is, with the growth or the order $p$ a family of functions
enlarges for which an improper integral converges. The limit case of
the limit by Cesaro is the Cauchy limit:
\par $(C) \int_0^{\infty }f(t)dt := \lim_{\epsilon \to +0}
\int_0^{\infty }f(t)\exp (-\epsilon t)dt$. \\
For two-sided integrals convergence of improper integrals by Cesaro
of order $p$ is defined by the equality:
\par $(C,p) \int_{-\infty }^{\infty }f(t)dt := \lim_{b\to \infty }
\int_{-b}^bf(t)(1-|t|/b)^pdt$ \\
and by Cauchy:
\par $(C) \int_{-\infty }^{\infty }f(t)dt := \lim_{\epsilon \to +0}
\int_{-b}^bf(t)\exp (-\epsilon |t|)dt$, \\
when these limits exist.

\par {\bf 18. An application of the noncommutative Laplace transformation
to (super)differential equations.}
\par Let a differential equation be given
\par $(1)$ $L[x]=f(t)$, где $L[x] := (d^nx(t)/dt^n) a_0+... +
(dx/dt) a_{n-1}+x a_n,$ $a_0\ne 0$ \\
and an initial conditions are imposed:
\par $(2)$ $x(0)=x_0,$ $x'(0)=x_1$,..., $x^{(n-1)}(0)=x_{n-1}$. \\
Suppose that a function $f(t)$ and a solution $x(t)$ and its
derivatives $x^{(1)}(t),...,x^{(n)}(t)$ are originals. Then
$X(p)(p^na_0+p^{n-1}a_1+...+a_n) = F(p)+ x_0
(p^{n-1}a_0+...+pa_{n-2}+ a_{n-1}) + x_1
(p^{n-2}a_0+...+pa_{n-3}+a_{n-2}) +...+ x_{n-1} a_0$, where
$X(p):={\cal F}(x,pt;p;0)$, $F(p)={\cal F}(f,pt;p;0)$, $a_j\in \bf
H$ in the case of $\bf H$, $a_j\in \bf R$ in the case of ${\cal
A}_r$ with $3\le r\in \bf N$. Since ${\cal F}(f'(t),u;p;\zeta )=
F(p;\zeta )\Delta _p - f(0)$ for $u=p_0t+\zeta _0+M(p,t;\zeta )$,
where $F(p;\zeta )\Delta _p := p_0F(p;\zeta )+p_1F(p;\zeta -i_1\pi
/2) +...+p_{2^r-1}F(p;\zeta -i_{2^r-1}\pi /2)$, then $f^{(n)}(t)=
F(p;\zeta )\Delta _p^n-f(0)p^{n-1}-...-f^{(n-2)}(0)p-f^{(n-1)}(0)$,
since $\exp (i_k\pi /2)=i_k$ for each $k=1,2,...,2^r-1$, while
$c\Delta _p=cp$ for $c=const \in {\cal A}_r$ in view of Formula
2.44(i), because it is possible to suppose that $\zeta =0$ for
$c=const $. Thus,
\par $(3)$ $X(p)A(p)=F(p)+B(p)$, where
\par $(4)$ $A(p)=(p^na_0+p^{n-1}a_1+...+a_n)$ and
\par $(5)$ $B(p)=x_0(p^{n-1}a_0+...+pa_{n-2}+
a_{n-1}) + x_1(p^{n-2}a_0+...+pa_{n-3}+a_{n-2}) +...+ x_{n-1}a_0$,
in the case $u=pt$, $\zeta =0$. For $u=p_0t+\zeta _0+M(p,t;\zeta )$
we have the equation
\par $(6)$ $X(p)A(\Delta _p)=F(p;\zeta )+B(p)$, \\
where $A$ is already the polynomial of the operator $\Delta _p$,
\par $(7)$ $B(p)=[(x_0\Delta _p^{n-1})a_0+...+(x_0\Delta _p)a_{n-2}+
a_{n-1}] + [(x_1\Delta _p^{n-2})a_0+...+(x_1\Delta
_p)a_{n-3}+a_{n-2}] +...+ x_{n-1}a_0$. In particular, for the skew
field of the quaternions ${\bf K}=\bf H$ or an algebra of octonions
${\bf K}=\bf O$, due to the alternativity of the algebra  $\bf K$
there follows the operator solution
\par $(8)$ $X(p)=(F(p;\zeta )+B(p;\zeta ))A^{-1}(p)$,\\  where
$A=A(p)$ or $A=A(\Delta _p)$ respectively, we put $\zeta =0$ in the
case $u=pt$, $F(p)=F(p;\zeta )$, $X(p)=X(p;\zeta )$, $B(p)=B(p;\zeta
)$, $A(\Delta _p)=A(\Delta _p;\zeta )$, $\zeta $ is the parameter of
the initial phases. In the given case it is more useful to write the
argument on the left from the operator. If Equation $(1)$ for the
initial Condition $(2)$ has a solution, $x(t)$,
$x^{(1)}(t)$,...,$x^{(n)}(t)$ satisfy conditions imposed on
originals, then $x(t)$ satisfies Equation $(7)$. The operator
$\Delta _p$ can be written in the form:
\par $F(p;\zeta )\Delta _p = Re <p,F(p;\zeta )L>$, \\
where $F(p;\zeta )L:=i_0F(p;\zeta )+ i_1F(p;\zeta - i_1\pi /2)+
...+i_{2^r-1}F(p;\zeta -i_{2^r-1}\pi /2)$, also
\par  $(f)\Delta _p={\cal F}(f,u;p;\zeta )+\delta df(t)/dt$, \\
where $u=p_0t+\zeta _0 +M(p,t;\zeta )$, $\delta $ is the generalized
function from Example 2.55, $<a,b> := \sum_{l=1}^n\mbox{ }^l{\tilde
a}\mbox{ }^lb$, where $a=(\mbox{ }^1a,...,\mbox{ }^na)$, $a, b\in
{\cal A}_r^n$, $n\in \bf N$.
\par {\bf 19. Example.} 1. Consider the equation $x^{(2)}+xa^2=b\sin
(at)$, where $b\in \bf H$ in the case of the algebra ${\cal A}_2=\bf
H$, $b\in \bf R$ in the case of the algebra ${\cal A}_r$ with $3\le
r\in \bf N$, $a\in \bf R$. Then the operator equation for $u=pt$ has
the form: $X(p)(a^2+p^2)= =b(a^2+p^2)^{-1}a+ x_0p+x_1$. For algebras
$\bf H$ and $\bf O$ the solution has the form:
$X(p)=ab(p^2+a^2)^{-2}+ (x_0p) (p^2+a^2)^{-1}+ x_1(p^2+a^2)^{-1}$.
Example 2.26.2, Theorem 2.27 about integration of originals 2.34
give: $ab(p^2+a^2)^{-1}= {\cal F}([2^{-1}b\int_0^bt\sin
(at)dt];pt;p;\zeta )= {\cal F}(b(2a^2)^{-1}(\sin (at) - at \cos
(at);pt;p;\zeta )$. If $x_0\in \bf R$, then the solution can be
written in the form:
\par $x(t)=(x_1+b(2a)^{-1})a^{-1}\sin (at) +(x_0-bt(2a)^{-1})
\cos (at)$.
\par 2. The equation $x^{(3)}+x=1$ with zero initial conditions
has the operator solution $X(p)=p^{-1}(p^3+1)^{-1}$. By Theorem of
decomposition 2.46 we get the original solution
$x(t)=1-e^{-t}/3-(2/3) Re \exp [(1/2+3^{1/2}S/2)t]=
1-e^{-t}/3-(2/3)\exp (t/2) \cos (t3^{1/2}/2)$, where $S=p'/|p'|$ for
$p':=p-Re (p)\ne 0$, while for $p'=0$, $S$ is such that $S\in {\cal
A}_r$, $|S|=1$.

\section{The Mellin transformation over Cayley-Dickson algebras.}
\par The Mellin transformation is based on the two sided
transformation of the Laplace type over Cayley-Dickson algebras,
which was presented in section 3 above.
\par {\bf 1. Remark.} If $f$ is an original function of the two-sided
Laplace transformation over the Cayley-Dickson algebra ${\cal A}_r$
and $g(\tau )=f(\ln \tau )$ for each $0<\tau <\infty $, then
Conditions 3.1(1-3) for $f$ are equivalent to the following
Conditions M(1-3):
\par $M(1)$ $g(\tau )$ satisfies the H\"older condition:
$|g(\tau +h)-g(\tau )| \le A |h|^{\alpha }$ for each $|h|<\delta $
(where $0<\alpha \le 1$, $A=const >0$, $\delta >0$ are constants for
a given $\tau $) everywhere on $\bf R$ may be besides points of
discontinuity of the first kind. On each finite segment $[a,b]$ in
$(0, \infty )$ a function $g$ may have only a finite number of
points of discontinuity and of the first kind only.
\par $M(2)$ $|g(\tau )| < C_1 \tau ^{s_0}$ for each $0<\tau <1$,
where $C_1=const >0$, $s_0=s_0(g)=const \in \bf R$.
\par $M(3)$ $|g(\tau )|<C_2 \tau ^{-s_1}$ for each $\tau \ge 1$, that is,
$g(\tau )$ is growing not faster, than the power function, where
$C_2=const>0$, $s_1=s_1(g)\in \bf R$.
\par This is because of the fact that the logarithmic function
$\ln : (0,\infty )\to (-\infty ,\infty )$ is the diffeomorphism.
\par {\bf 2. Definition.} In the two-sided integral transformation
of the Laplace type substitute variables $p$ on $-p$ and $\zeta $ on
$-\zeta $ and $t$ on $\tau =e^t$, then the formula take the form:
\par $(1)$ ${\cal M}(g,u;p;\zeta ) := \int_0^{\infty }
f(\ln \tau )\exp (-u(-p,\ln \tau ,-\zeta ))\tau ^{-1}d\tau $, \\
where $f$ is an original function, $g(\tau )=f(\ln \tau )$ for each
$0<\tau <\infty $, a function   $u(p,t,\zeta )$ is given by Formulas
2.1 or 2.3(1-4). \par For a specified basis $ \{ N_0, N_1,...,
N_{2^r-1} \} $ of generators of the Cayley-Dickson algebra ${\cal
A}_r$ we can write the notation in more details $\mbox{ }_N{\cal
M}(g,u;p;\zeta )$ if necessary.
\par {\bf 3. Theorem.} {\it If an original function $g(\tau )$ satisfies
Conditions M(1-3), where $s_0<s_1$, then its image ${\cal
M}(g,u;p;\zeta )$ is holomorphic by $p$ in the domain $\{ z\in {\cal
A}_r: s_0< Re (z)<s_1 \} $, as well as by $\zeta \in {\cal A}_r$,
where $2\le r\in \bf N$.}
\par {\bf Proof.} The application of Theorem 3.4 to $g(\tau )=
f(\ln \tau )$ gives the statement of this theorem due to Formulas
2.3(6).
\par {\bf 4. Proposition.} {\it If $\mbox{ }_N{\cal M}(g,u;p;\zeta )$ and
$\mbox{ }_N{\cal M}(w,u;p;\zeta )$ are images of original functions
$g(\tau )$ and $w(\tau )$ in the domains $s_0(g)< Re (p)< s_1(g)$
and $s_0(w)<Re (p)<s_1(w)$ with values in ${\cal A}_r$, then for
each $\alpha , \beta \in {\cal A}_r$ in the case ${\bf K}=\bf H$;
and also either $g$ and $w$ with values in $\bf R$ and each $\alpha
, \beta \in {\cal A}_r$ or $g$ and $w$ with values in ${\cal A}_r$
and each $\alpha , \beta \in \bf R$ in the case of ${\cal A}_r$ the
function $\alpha \mbox{ }_N{\cal M}(g,u;p;\zeta ) + \beta \mbox{
}_N{\cal M}(w,u;p;\zeta ) $ is the image of the function $\alpha
g(\tau ) +\beta w(\tau )$ in the domain $\max (s_0(g),s_0(w))< Re
(p)< \min (s_1(g),s_1(w))$.}
\par {\bf Proof.} Applying Proposition 3.5 we get the statement of this
proposition.
\par {\bf 5. Theorem.} {\it If $g(\tau )$ is an original function with
parameters $s_0, s_1$, then for each $\alpha >0$:
\par $(1)$  ${\cal M}((g(\alpha \tau ),u;p;\zeta )={\cal M}(
(g(\tau ),u;p; \zeta - p\ln \alpha )$ \\
for $u(p,t;\zeta )= p_0t + \zeta _0 + M(p,t;\zeta )$ or $u(p,t;\zeta
)= pt + \zeta $ over ${\cal A}_r$ with $2\le r<\infty $,
\par $(2)$ ${\cal M}((g(\alpha \tau ),u;p;0)={\cal M}((g
(\tau ),u;p;0)\alpha ^{-p}$ \\
for $u(p;t;0)=pt$ over either ${\bf K}=\bf H$ or ${\bf K}=\bf O$
with $f({\bf R})\subset \bf K$, or over ${\cal A}_r$ with $4\le r$
for $f({\bf R})\subset \bf R$, in the domain $s_0(g)< Re (p)<
s_1(g)$.}
\par {\bf Proof.} Since $\ln (\alpha \tau )=\ln \alpha +\ln \tau $,
then application of Theorem 3.8 and Definition 2 gives the statement
of this theorem.
\par {\bf 6. Theorem.} {\it If $g(\tau )$ is an original function
with values in ${\cal A}_r$ for $2\le r<\infty $, $b\in \bf R$, then
${\cal M}(\tau ^bg(\tau ),u;p;\zeta )= {\cal M}(g(\tau ),u;p+b;\zeta
)$ for each $s_0(g)+b<Re (p)<s_1(g)+b$, where $u=pt+\zeta $ or
$u=p_0t+\zeta _0 + M(p,t;\zeta )$.}
\par {\bf Proof.} Since $\tau ^b=\exp (b \ln \tau )$ for each
$\tau >0$, then Theorem 3.9 gives the statement of this Theorem.
\par {\bf 7. Theorem.} {\it If ${\cal M}(g(\tau ),u;p;\zeta )$
is defined in the domain $s_0< Re (p) <s_1$, then
\par ${\cal M}(g(\tau ^b),u;p;\zeta )={\cal M}(g(\tau ),u;p/b;\zeta )/b$
for $b>0$ in the domain $bs_0< Re (p)<bs_1$,\par  ${\cal M}(g(\tau
^b),u;p;\zeta )= - {\cal M}(g(\tau ),u;p/b;\zeta )/b$ for $b<0$ in
the domain $bs_1< Re (p)<bs_0$.}
\par {\bf Proof.} Apply Theorem 3.7 to the original function
$f(bt)=g(\tau ^b)$ with $t=\ln \tau $ for each $\tau >0$, since
$s_0(g)=-s_1(f)$, $s_1(g)=-s_0(f)$.
\par {\bf 8. Theorem.} {\it Let $g, g',..., g^{(n)}$ be original
functions such that $g(\tau )\tau ^{p-1}|_0^{\infty }=0$,...,
$g^{(n-1)}\tau ^{p-n}|_0^{\infty }=0$, let also $g$ be real-valued
for $r\ge 4$, $g$ is with values in $\bf K$ for ${\bf K}=\bf H$ or
${\bf K}=\bf O$; then
\par $(1)$ ${\cal M}(g^{(n)}(\tau ),u;p;0) = (-1)^n{\cal M}(g(\tau ),
u;p-n;0)(p-1)...(p-n)$ \\
in the domain $s_0<Re (p)<s_1$ in ${\cal A}_r$, when $s_0<s_1$,
where $n\in \bf N$, $n\ge 1$, $s_0 := \max
(s_0(g),s_0(g'),...,s_0(g^{(n)}))$, $s_1 := \min
(s_1(g),s_1(g'),..., s_1(g^{(n)}))$, $u(p,t,0)=pt$. \par If $g(\tau
)$, $g'(\tau )\tau $,..., $g^{(n)}(\tau )\tau ^n$ are originals such
that $g(\tau )\tau ^p|_0^{\infty }=0$, $g'(\tau )\tau
^{p+1}|_0^{\infty }=0$,..., $g^{(n)}(\tau )\tau ^{p+n}|_0^{\infty
}=0$, then
\par $(2)$ ${\cal M}(g^{(n)}(\tau )\tau ^n,u;p;\zeta )=(-1)^n{\cal
M}(g(\tau ),u;p;\zeta )p(p+1)...(p+n-1)$ in the domain $s_0<Re
(p)<s_1$ in ${\cal A}_r$, when $s_0<s_1$, where $s_0 := \max
(s_0(g),s_0(g')-1,...,s_0(g^{(n)})-n)$, $s_1 := \min (s_1(g),
s_1(g')-1,...,s_1(g^{(n)})-n)$, $n\in \bf N$, $n\ge 1$, $u(p,t,\zeta
)=pt+\zeta $.}
\par {\bf Proof.} For $s_0<s_1$ and $n=1$ in view of Theorem 3.11
integration by parts gives:
\par $\int_0^{\infty }g'(\tau )\tau ^{p-1}d\tau =g(\tau )\tau
^{p-1}|_0^{\infty } - (\int_0^{\infty } g(\tau )\tau ^{p-2}d\tau)
(p-1) = -{\cal M}(g(\tau ),u,p-1;\zeta )(p-1)$, \\
since either $g$ is real-valued for $r\ge 4$, or $g$ is with values
in $\bf K$ for ${\bf K}=\bf H$ or ${\bf K}=\bf O$, but $\bf O$ is
the alternative algebra. Therefore, ${\cal M}(g'(\tau ),u;p;0) =
-{\cal M}(g(\tau ), u;p-1;0)(p-1)$. Then we apply this integration
by parts formula by induction and get Formula $(1)$.
\par If $|g(\tau )|<C_1\tau ^{s_0}$ for each $0<\tau \le 1$,
then $|g(\tau )\tau ^k|<C_1\tau ^{s_0+k}\le C_1\tau ^{s_0}\le
C_1\tau ^{s_0-k}$ for each $k\ge 0$ and $0<\tau \le 1$,
consequently, $s_0(g(\tau )\tau ^k)\le s_0(g)-k$. If $|g(\tau
)|<C_2\tau ^{-s_1}$ for each $\tau \ge 1$, then $|g(\tau )\tau
^k|<C_2\tau ^{-s_1+k}$ for each $k\ge 0$ and $\tau \ge 1$, hence
$s_1(g(\tau )\tau ^k)\ge s_1(g)-k$.
\par {\bf 9. Theorem.} {\it (1). Let $g(\tau )$ and $g'(\tau )$
be original functions for ${\cal A}_r$ with $2\le r<\infty $ such
that there exists \par $(a)$ $\lim_{\tau \to \tau _0}g(\tau )\tau
^{p-1}$ for values $\tau _0=0$ and $\tau _0=\infty $, then
\par $(b)$  ${\cal M}(g'(\tau ),u;p;\zeta )=-(p_0-1)
{\cal M}(g(\tau ),u;p-1;\zeta )-p_1{\cal M}(g(\tau ),u;p-1;\zeta
+i_1\pi /2) -...- p_{2^r-1} {\cal M}(g(\tau ),u;p-1;\zeta
+i_{2^r-1}\pi /2)$
\\ in the domain $s_0<Re (p)<s_1$ in ${\cal A}_r$, when $s_0<s_1$,
where $s_0:=\max (s_0(g),s_0(g'))$, $s_1:=\max (s_1(g),s_1(g'))$,
$u(p,t,\zeta )=p_0t+\zeta _0+M(p,t,\zeta )$.
\par $(2)$. If $g(\tau )$ and $g'(\tau )\tau $ are originals such
that there exists \par $(c)$ $\lim_{\tau \to \tau _0} g(\tau )\tau
^p=0$ for each value $\tau _0=0$ and $\tau _0=\infty $, then
\par $(d)$  ${\cal M}(g'(\tau )\tau ,u;p;\zeta )=-p_0 {\cal M}(g(\tau
),u;p;\zeta )-p_1{\cal M}(g(\tau ),u;p;\zeta +i_1\pi /2) -$\\ $...-
p_{2^r-1} {\cal M}(g(\tau ),u;p;\zeta +i_{2^r-1}\pi /2)$ \\
in the domain $s_0<Re (p)<s_1$ in ${\cal A}_r$, where $s_0:=\max
(s_0(g),s_0(g')-1)$, \\ $s_1:=\min (s_1(g),s_1(g')-1)$, $u(p,t;\zeta
)=p_0t+\zeta _0 +M(p,t;\zeta )$.}
\par {\bf Proof.} Integration by parts gives:
\par $\int_0^{\infty }g'(\tau )\exp (-u(-p,\ln \tau ;-\zeta ))\tau
^{-1}d\tau = g(\tau )e^{-u}\tau ^{-1}|_0^{\infty } -\int_0^{\infty }
g(\tau )(\partial e^{-u}\tau ^{-1}/\partial \tau )d\tau $, \\
where $\exp (-u(-p,\ln \tau ;-\zeta ))\tau ^{-1}= \exp (p_0\ln \tau
+\zeta _0)(\cos (p_1\ln \tau +\zeta _1) - i_1 \sin (p_1\ln \tau
+\zeta _1)$ \\ $\cos (p_2\ln \tau +\zeta _2) + i_2 \sin (p_1\ln \tau
+\zeta _1)\sin (p_2\ln \tau +\zeta _2)\cos (p_3\ln \tau +\zeta
_3)-...+i_{2^r-2}\sin (p_1\ln \tau +\zeta _1)...\cos (p_{2^r-1}\ln
\tau +\zeta _{2^r-1}) - i_{2^r-1}\sin (p_1\ln \tau +\zeta _1)...\sin
(p_{2^r-1}\ln \tau +\zeta _{2^r-1}))\tau ^{-1}$, \\
$\exp (p_0\ln \tau +\zeta _0)=\tau ^{p_0-1}\exp (\zeta _0)$, $\exp
(\zeta _0)=const$, where $\tau >0$. Therefore,
\par $\partial [\exp (-u(-p,\ln \tau ; -\zeta )\tau ^{-1}]/\partial
\tau = \partial [\exp (\zeta _0)\tau ^{p_0-1}\exp (-M(-p,\ln \tau
;-\zeta )]/\partial \tau $ \\ $= (p_0-1)\exp (\zeta _0)\tau
^{p_0-2}\exp (-M) +p_1 \exp (-u(-p,\ln \tau ;-\zeta -i_1\pi /2)/\tau
^2 +...$ \\ $+p_{2^r-1} \exp (-u(-p,\ln \tau
;-\zeta -i_{2^r-1}\pi /2))/\tau ^2$ \\
due to Formula 2.31(2), since $d \ln \tau /d\tau =1/\tau $, where
$\exp (\zeta _0)\tau ^{p_0-2}\exp (-M) = \tau ^{-1}\exp (-u(-p+1,\ln
\tau ;-\zeta ))$. From this Formula $(b)$ follows, since the phase
factor $\exp (\phi (\tau ))$ with $Re (\phi (\tau ))=0$ for each
$\tau >0$ does not influence on Limit $(a)$.
\par Consider now $g'(\tau )\tau $, then from the integration by parts
formula we get:
\par $\int_0^{\infty }g'(\tau )\tau \exp (-u(-p,\ln \tau ;-\zeta
))\tau ^{-1}d\tau = g(\tau )\exp (-u(-p,\ln \tau ;-\zeta
))|_0^{\infty } $ \\  $- p_0\int_0^{\infty } g(\tau )\exp (-u(-p,\ln
\tau ; -\zeta )\tau ^{-1}d\tau - p_1 \int_0^{\infty }g(\tau )\exp
(-u(-p, \ln \tau ; -\zeta -i_1\pi /2))\tau ^{-1}d\tau -...-
p_{2^r-1} \int_0^{\infty }g(\tau )\exp (-u(-p,\ln \tau ; -\zeta
-i_{2^r-1}\pi /2)) \tau ^{-1}d\tau $, \\
from which Formula $(d)$ follows, since the multiplier $\exp (\phi
(\tau ))$ with $Re (\phi (\tau ))=0$ has $|\exp (\phi (\tau ))|=1$
for each $\tau >0$ and it does not influence on Limit $(b)$.
\par {\bf 10. Theorem.} {\it If $g(\tau )$ is an original function
with values in ${\cal A}_r$ for $2\le r<\infty $, then
\par $(1)$ ${\cal M}(g,u;p;0)^{(n)}(p).(\mbox{ }_ 1h,...,\mbox{ }_nh)
= {\cal M}((\ln ^n\tau )g(\tau )\mbox{ }_1h...\mbox{ }_nh,u;p;0)$
for each $n\in \bf N$ and $\mbox{ }_1h,...,\mbox{ }_nh\in {\bf
R}\oplus p'{\bf R}\subset {\cal A}_r$ for $u=pt$, where $p' := p-Re
(p)$. If $u(p,t) := p_0t + M(p,t;\zeta )+\zeta _0$, then
\par $(2)$  $(\partial {\cal M}(g(\tau ),u; p;\zeta )/
\partial p).h = {\cal M}(g(\tau )\ln \tau ,u;p;\zeta )h_0 $ \\
$+ {\cal M}(g(\tau )\ln \tau ,u;p;\zeta + i_1\pi /2)h_1 +...+{\cal
M}(g(\tau )\ln \tau ,u; p; \zeta + i_{2^r-1}\pi /2)h_{2^r-1}$ \\
for each $h=h_0i_0+...+h_{2^r-1}i_{2^r-1}\in {\cal A}_r$, where
$h_0,...,h_{2^r-1}\in {\bf R}$. Both formulas are accomplished in
the domain $s_0(g)<Re (p)<s_1(g)$.}
\par The {\bf proof} follows from Theorem 3.12 applying change of
variables in accordance with Definition 2, since $\mbox{ }_1h,...,
\mbox{ }_nh\in {\bf R}\oplus p'{\bf R}$ commute with $p'+c$ for each
$c\in \bf R$, $\ln \tau \in \bf R$ for each $\tau >0$ and $\bf R$ is
the center of the Cayley-Dickson algebra ${\cal A}_r$ for $r\ge 2$.
\par {\bf 11. Definition.} Let $g$ and $w\in L^1((0,\infty ),
\lambda ,{\cal A}_r)$, where $\lambda $ is the Lebesgue measure on
$(0,\infty )$, then their convolution $g{\hat *}w$ is defined by the
formula:
\par $(g{\hat *}w)(b):=\int_0^{\infty }g(a)w(b/a)a^{-1}da$ \\
for each $b\in (0,\infty )$. We spread this definition on arbitrary
originals $g$ and $w$ whenever $g{\hat *}w$ exists as an original.
\par {\bf 12. Theorem.} {\it Let $g(\tau )$ and
$w(\tau )$ be original functions, where $w$ is real-valued for each
$\tau \in (0,\infty )$ and their convolution $g{\hat *}w$ is an
original also. If either
\par $(1)$ $g$, $w$ and $g{\hat *}w$ have converging Mellin
integrals over ${\cal A}_r$, $2\le r<\infty $, in their common
domain of convergence or
\par $(2)$ $g$ and $w$ have absolutely converging Mellin
integrals over ${\cal A}_r$, $2\le r<\infty $, in their common
domain of convergence, then
\par ${\cal M}((g{\hat *}w(\tau ),u;p;0)
= {\cal M}(g(\tau ),u;p;0) {\cal M}(w(\tau ),u;p;0)$ \\
for each $p\in {\cal A}_r$ with $s_0<Re (p)<s_1$, where
$u(p,t;0)=pt$, $s_0 = \max (s_0(g), s_0(g))$, $s_1=\min
(s_1(g),s_1(w))$ and $2\le r<\infty $.}
\par {\bf Proof.} Here the convolution is defined relative to
the multiplicative group $(0,\infty )$ apart from the convolution
relative to the additive group $(-\infty ,\infty )$ in Theorem 3.13.
Let $f(\ln a)=g(a)$ and $q(\ln b- \ln a)=w(b/a)$ for each $a>0$ and
$b>0$, then:
\par $\int_{-\infty }^{\infty }f(t)q(\ln b- t)dt=
\int_0^{\infty }g(a)w(b/a)a^{-1}da=(g{\hat *}w)(b)$. Therefore,
applying Theorem 3.13 we get the statement of this theorem.
\par {\bf 13. Theorem.} {\it Let $g(\tau )$ be an original function
with values in ${\cal A}_r$ for $2\le r<\infty $, $u(p,t;0)=pt$,
where $s_0(g)\le 0\le s_1(g)$, then
\par $(1)$ ${\cal M}(w(\tau ),u;p;0)p =
{\cal M}(g(\tau ),u;p;0)$ \\
for $w(\tau ) := \int_0^{\tau } g(a)a^{-1}da$ in the domain
$s_0(g)<Re (p)<\min (s_1(g),0)$, in particular,
\par ${\cal M}(w(\tau ),u;p;0) =
{\cal M}(g(\tau ),u;p;0)p^{-1}$ over the algebra ${\bf K}=\bf H$ or
${\bf K}=\bf O$;
\par $(2)$ ${\cal M}(w(\tau ),u;p;0)p =
{\cal M}(g(\tau ),u;p;0)$ \\
for $w(\tau ) := \int_{+\infty }^{\tau } g(a)a^{-1}da$ in the domain
$\max (s_0(g),0)<Re (p)<s_1(g)$, in particular,
\par ${\cal M}(w(\tau ),u;p;0) =
{\cal M}(g(\tau ),u;p;0)p^{-1}$ over the algebra ${\bf K}=\bf H$ or
${\bf K}=\bf O$, when these domains are non void.}
\par {\bf Proof.} Put $f(\ln \tau )=g(\tau )$, $x=\ln a$ and
$t=\ln \tau $, then $\int_0^{\tau }g(a)a^{-1}da=\int_{-\infty
}^tf(x)dx$ and $\int_{+\infty }^tf(x)dx=\int_{+\infty }^{\tau
}g(a)a^{-1}da$. Therefore, from Theorem 3.14 follows the statement
of this theorem, since $\min (a,b)= -\max (-a,-b)$ for each $a, b\in
\bf R$.
\par {\bf 14. Theorem.} {\it Let $g(\tau )$ be an original function
such that \par $\mbox{ }_N{\cal M}(g,u;p;\zeta ) :=
\sum_{j=0}^{2^r-1} \mbox{ }_NG_{u,j}(p;\zeta )N_j$ be its image,
where a function $g$ is written in the form \par $g(\tau ) =
\sum_{j=0}^{2^r-1} g_j(\tau )N_j$, $g_j: (0,\infty )\to \bf R$ for
each $j=0,1,...,2^r-1$, $g((0,\infty ))\subset \bf K$ for ${\bf
K}=\bf H$ or ${\bf K}=\bf O$, $g((0,\infty ))\subset \bf R$ over the
Cayley-Dickson algebra ${\cal A}_r$ with $4\le r\in \bf N$,
\par $\mbox{ }_NG_{u,j}(p;\zeta ) := \int_0^{\infty
}g_j(\tau )\exp (-u(-p,\ln \tau ; - \zeta ))\tau ^{-1}d\tau $. Then
at each point $\tau $, where $g(\tau )$ satisfies the H\"older
condition the equality is accomplished:
\par $(i)$ $g(\tau ) = (2\pi N_1)^{-1} Re (S{\tilde N}_1)
\sum_{j=0}^{2^r-1} \int_{a-S\infty }^{a+S\infty }\mbox{
}_NG_{u,j}(p;\zeta ) \exp (u(-p,\ln \tau ; - \zeta ))dp)N_j$ \\
in the domain $s_0(g)< Re (p) <s_1(g)$, where $u(p,t;\zeta
)=pt+\zeta $ or $u(p,t;\zeta )=p_0t+M_N(p,t;\zeta )+\zeta _0$ and
the integral is taken along the straight line $p(\theta )=a+S\theta
\in {\cal A}_r$, $\theta \in \bf R$, $S\in {\cal A}_r$, $Re (S)=0$,
$|S|=1$, while the integral is understood in the sense of the
principal value.}
\par {\bf Proof.} Putting $t=\ln \tau $ and substituting $p$ on $-p$
and applying Theorem 3.15 we get the statement of this theorem.
\par {\bf 15. Theorem.} {\it If a function $\mbox{
}_NG_u(p)$ is analytic by the variable $p\in {\cal A}_r$ in the
domain $W := \{ p\in {\cal A}_r: s_0< Re (p) <s_1 \} $, where $2\le
r\in \bf N$, $g((0,\infty ))\subset \bf K$ for ${\bf K}=\bf H$ or
${\bf K}=\bf O$, also for an arbitrary Cayley-Dickson algebra ${\cal
A}_r$ with $g((0,\infty ))\subset \bf R$ for $4\le r<\infty $,
$u(p,t;\zeta )=pt+\zeta $ or $u(p,t;\zeta ) := p_0t + M(p,t;\zeta
)+\zeta _0$. Let also $\mbox{ }_NG_u(p)$ can be written in the form
$\mbox{ }_NG_u(p)=\mbox{ }_NG^0_u(p) + \mbox{ }_NG^1_u(p)$, where
$\mbox{ }_NG^0_u(p)$ is holomorphic by $p$ in the domain $s_0<Re
(p)$, also $\mbox{ }_NG^1_u(p)$ is holomorphic by $p$ in the domain
$Re (p)<s_1$, moreover, for each $a>s_0$ and $b<s_1$ there exists
constants $C_a>0$, $C_b>0$ and $\epsilon _a
>0$ and $\epsilon _b>0$ such that
\par $(1)$ $|\mbox{ }_NG^0_u(p)|\le C_a\exp (-\epsilon _a |p|)$
for each $p\in {\cal A}_r$ with $Re (p)\ge a$,
\par $(2)$ $|\mbox{ }_NG^1_u(p)|\le C_b\exp (-\epsilon _b |p|)$
for each $p\in {\cal A}_r$ with $Re (p)\le b$, where $s_0$ and $s_1$
are fixed, while the integral
\par $(3)$ $\int_{w-S\infty }^{w+S\infty }\mbox{
}_NG^k_u(p)dp$ \\ converges absolutely for $k=0$ and $k=1$ for
$s_0<w<s_1$, then $\mbox{ }_NG_u(p)$ is the image of the function
\par $(4)$ $g(\tau )=(2\pi )^{-1}{\tilde S}\int_{w-S\infty
}^{w+S\infty }\mbox{ }_NG_u(p)\exp (u(-p,\ln \tau ;0))dp$.}
\par {\bf Proof.} The change of the variable $p$ on $-p$ and the
substitution $t=\ln \tau $ for $\tau >0$ with the help of Theorem
3.16 gives the statement of this Theorem.

\par Chair of Applied Mathematics,
\par Moscow State Technical University MIREA,
av. Vernadsky 78,
\par Moscow, Russia
\par e-mail: sludkowski@mail.ru
\end{document}